\documentclass[a4paper]{article}
\usepackage{mathrsfs}
\usepackage{latexsym,bm}
\usepackage{amssymb,amsmath,amsthm}
\usepackage[english,german]{babel}
\usepackage[colorlinks=true]{hyperref}
\usepackage{color}

\allowdisplaybreaks[4]

\newtheorem{theorem}{Theorem}[section]
\newtheorem{lemma}[theorem]{Lemma}
\newtheorem{remark}[theorem]{Remark}
\newtheorem{proposition}[theorem]{Proposition}
\newtheorem{definition}[theorem]{Definition}
\newtheorem{corollary}[theorem]{Corollary}
\numberwithin{equation}{section}
\hypersetup{linkcolor=blue,urlcolor=red,citecolor=red}
\title{Stabilizability with bounded feedback for analytic linear control systems\thanks{This work was supported by the National Natural Science Foundation of China
under grants 12171377, 12571483.}}
\author{Yaxing Ma{\footnote{School of
Mathematics and Statistics, Wuhan University, Wuhan 430072, China
(e-mail: yaxingma@yeah.net).}}\and Emmanuel Tr\'{e}lat{\footnote{Sorbonne Universit\'e, Universit\'e Paris Cit\'e, CNRS, Inria, Laboratoire Jacques-Louis Lions, LJLL, F-75005 Paris, France (email: emmanuel.trelat@sorbonne-universite.fr).}}\and Lijuan Wang{\footnote{School of
Mathematics and Statistics, Wuhan University, Wuhan 430072, China
(e-mail: ljwang.math@whu.edu.cn).}}\and Huaiqiang Yu{\footnote{School of Mathematics and KL-AAGDM, Tianjin University, Tianjin, 300350, China (email: huaiqiangyu@yeah.net).}}}
\date{}
\begin{document}
\selectlanguage{english}
\maketitle
\begin{abstract}
In this paper, we give sufficient conditions under which linear abstract control systems for which the semigroup is analytic are stabilizable with a bounded feedback. We obtain various characterizations of that property, which extend some earlier works.  We illustrate our findings with several examples.
\end{abstract}


\section{Introduction and main results}\label{yu-section-1}
\paragraph{Notations.}
Let $\mathbb{N}^+:=\mathbb{N}\setminus\{0\}$, $\mathbb{R}^+:=(0,+\infty)$, $\mathbb{C}_{\gamma}^+:=\{z\in\mathbb{C}:\mbox{Re}z>\gamma\}$ and $\mathbb{C}^{-}_{\gamma}:=\{z\in\mathbb{C}:\mbox{Re}z<\gamma\}$ for $\gamma\in\mathbb{R}$.  Given a Banach space $X$, we denote its norm and dual space as $\|\cdot\|_X$ and $X'$, respectively. For a Hilbert space $X$, we denote its inner product by $\langle\cdot,\cdot\rangle_X$. Given two Banach spaces $X_1$ and $X_2$, we denote
    the class of bounded linear operators from $X_1$ to $X_2$ as $\mathcal{L}(X_1;X_2)$ with the usual operator norm. We write $\mathcal{L}(X_1):=\mathcal{L}(X_1;X_2)$ for simplicity if $X_1=X_2$.
    Given an unbounded (or bounded) and linear operator $L$ from $X_1$ to $X_2$, we denote its domain as $D(L):=\{f\in X_1:Lf\in X_2\}$, its adjoint operator as $L^*$,  its resolvent set and spectrum  as $\rho(L)$ and $\sigma(L)$, respectively. If $\Lambda$ is a closed subspace of $X$, we denote its dimension as $\mbox{Dim}(\Lambda)$. We use $C(\cdots)$  to denote a positive constant that depends on what are enclosed in the parenthesis.
\subsection{Reminders and preliminaries}
     Let  $X$ and $U$ be two separable and complex Hilbert spaces with $X'=X$ and $U'=U$.
     We make the following assumption:
     \begin{enumerate}
      \item [$(A_1)$] The linear operator $A:D(A)\subset X\to X$ generates  an analytic semigroup $S(\cdot)$ on $X$.
     \end{enumerate}
       Let $\rho_0\in(\sup\{\mbox{Re}\lambda:\lambda\in \sigma(A)\},+\infty)$ be fixed. Under Assumption $(A_1)$, we define the following spaces:
\begin{equation}\label{yu-3-26-5}
    X_{\gamma}:=
\begin{cases}
    D((\rho_0I-A)^\gamma)&\mbox{if}\;\;\gamma\geq 0,\\
    D((\rho_0I-A^*)^{-\gamma})'&\mbox{if}\;\;\gamma<0,
\end{cases}
    \;\;\mbox{and}\;\;X_\gamma^*:=
    D((\rho_0I-A^*)^\gamma)\;\;\mbox{for}\;\;\gamma\geq 0.
\end{equation}
   Here and throughout this paper, given a Hilbert space $V$  densely and continuously embedded in $X$, if we do not emphasize other situations,  $V'$ stands for the dual space of $V$ with respect to (w.r.t.) the pivot space $X$.
\begin{remark}\label{yu-remark-6-22-1}
    The following facts are collected in \cite[Section 4.1.3]{Trelat}.
\begin{enumerate}
  \item [$(i)$] For any $\gamma\geq 0$, $X_{\gamma}$ and $X^*_{\gamma}$, endowed with norms $\|f\|_{X_{\gamma}}=\|(\rho_0I-A)^{\gamma}f\|_X$,
    $\|g\|_{X^*_\gamma}=\|(\rho_0I-A^*)^\gamma g\|_X$ for any $f \in X_{\gamma}$, $g\in X_\gamma^*$,
    respectively, are two separable Hilbert spaces. Moreover,
      by \eqref{yu-3-26-5}, one can easily check that $D(A)=X_1$, and for each $\gamma\geq 0$, $(\rho_0I-A^*)^{-\gamma}\in\mathcal{L}(X;X^*_\gamma)$, $X_{\gamma}=(\rho_0I-A)^{-\gamma} X$, $X^*_{\gamma}=(\rho_0I-A^*)^{-\gamma} X$ and
$X_{-\gamma}'=X_\gamma^*$.
  \item [$(ii)$]  The operator $A$ (with domain $X_1$) has a unique extension $\widetilde{A}\in \mathcal{L}(X;X_{-1})$ in  the sense: $\langle \widetilde{A}\varphi,\psi\rangle_{X_{-1},X^*_1}\\
      =\langle \varphi, A^*\psi\rangle_{X}$
    for any $\varphi\in X,\;\;\psi\in X^*_1$.
 For any $\gamma\in[0,+\infty)$, $(\rho_0I-\widetilde{A})^{-\gamma}\in \mathcal{L}(X_{-\gamma};X)$ and it is isomorphic. The norm of $X_{-\gamma}$ is given as follows:
\begin{equation}\label{yu-3-28-10}
    \|f\|_{X_{-\gamma}}=\|(\rho_0I-\widetilde{A})^{-\gamma}f\|_X\;\;\mbox{for any}\;\;f \in X_{-\gamma}.
\end{equation}
\item[$(iii)$]
    Let $\widetilde{S}(\cdot):=(\rho_0I-\widetilde{A})S(\cdot)(\rho_0I-\widetilde{A})^{-1}$, which is an extension of $S(\cdot)$ (from $X$ to $X_{-1}$).
    Then $\widetilde{S}(\cdot)$ is an analytic semigroup on $X_{-1}$ and $\widetilde{A}$ with domain $X$ is its generator (see \cite[Proposition 2.10.4, Section 2.10, Chapter 2]{Tucsnak-Weiss}). Moreover, for any $\gamma> -1$, $\widetilde{S}(\cdot)$ is also a $C_0$-semigroup on $X_\gamma$ and the corresponding generator is
\begin{equation}\label{yu-25-6-20-1}
    \widetilde{A}_\gamma:=(\rho_0I-\widetilde{A})^{-\gamma}A(\rho_0I-\widetilde{A})^{\gamma}\;\;
    \mbox{with domain}\;\;D(\widetilde{A}_\gamma):=X_{1+\gamma}.
\end{equation}
     Furthermore, $\widetilde{A}\varphi=\widetilde{A}_{\gamma}\varphi$ for any $\varphi\in D(\widetilde{A}_\gamma)$.
\end{enumerate}
\end{remark}
    Under Assumption $(A_1)$, we consider the following controlled system:
\begin{equation}\label{yu-4-9-1}
    y_t(t)=Ay(t)+Bu(t),\;\;t\in\mathbb{R}^+,
\end{equation}
    where $u\in L^2_{loc}(\mathbb{R}^+;U)$ and the control operator $B$ satisfies the assumption:
    \begin{enumerate}
      \item [$(A_2)$] There exists  $\gamma\in[0,1)$ so that
$B\in \mathcal{L}(U;X_{-\gamma})$.
      \end{enumerate}
   Throughout this paper, we denote the system \eqref{yu-4-9-1} by the pair $[A,B]$ for simplicity.
\begin{remark}\label{yu-remark-8-00-1}
 Under Assumption $(A_2)$, we have $B\in\mathcal{L}(U;X_{-1})$. Moreover, according to $(ii)$ in Remark \ref{yu-remark-6-22-1}, Assumption $(A_2)$ can be replaced by $(\rho_0I-\widetilde{A})^{-\gamma}B\in
  \mathcal{L}(U;X)$. Furthermore, by Assumption $(A_2)$ and $(i)$ in Remark \ref{yu-remark-6-22-1}, we have that $B^*\in \mathcal{L}(X^*_\gamma;U)$, i.e.,  $B^*(\rho_0I-A^*)^{-\gamma}\in\mathcal{L}(X;U)$.
  It follows from \cite[Theorem 6.13, Section 2.6, Chapter 2]{Pazy} that  if $\gamma$ (in Assumption $(A_2)$) belongs to $[0,\frac{1}{2})$, then for each $T>0$, there exists $C(T)>0$ so that
\begin{equation}\label{yu-25-6-22-1}
    \int_0^T\|B^*S^*(t)\varphi\|_U^2\mathrm{dt}\leq C(T)\|\varphi\|_X^2\;\;\mbox{for any}\;\;\varphi\in X_1^*.
\end{equation}
\end{remark}
    For the convenience of later statements, we recall definitions of admissible control operator and admissible observation operator (see \cite[Sections 4.2, 4.3 and 4.4, Chapter 4]{Tucsnak-Weiss}) as follows.
\begin{definition}\label{yu-definition-6-21-1}
   If $B\in\mathcal{L}(U;X_{-1})$ (equivalently, $B^*\in\mathcal{L}(X_1^*;U)$) and there exist $T>0$ and $C(T)>0$ so that \eqref{yu-25-6-22-1} holds, then $B$ is called an admissible control operator (w.r.t. $A$) in $X$, or equivalently, $B^*$ is   an admissible observation operator (w.r.t. $A^*$) in $X$.
\end{definition}
\begin{remark}\label{yu-remark-6-23-2}
    Two remarks on Definition \ref{yu-definition-6-21-1} are presented in order.
\begin{enumerate}
  \item [$(i)$]
     It is known that $B$ is  an admissible control operator (w.r.t. $A$) in $X$ if and only if for each $T>0$, there exists $C(T)>0$ so that \eqref{yu-25-6-22-1} holds for any $\varphi\in X_1$ (see \cite[Proposition 4.3.2 and Theorem 4.4.3, Chapter 4]{Tucsnak-Weiss}).
  \item[$(ii)$] If Assumption $(A_2)$ holds for $\gamma\in[0,\frac{1}{2})$, then it follows from
  \cite[Theorem 6.13, Section 6, Chapter 2]{Pazy} that $B$ is admissible (w.r.t. $A$) in $X$. This, along with $(i)$ above, implies that $\int_0^\cdot\widetilde{S}(\cdot-s)Bu(s)ds$ is in
    $C([0,+\infty);X)$ for any $u\in L^2_{loc}(\mathbb{R}^+;U)$. Then by \cite[Proposition 4.2.5, Section 4.2, Chapter 4]{Tucsnak-Weiss}, we observe that the pair
  $[A,B]$ is well-posed in $X$ in this case. When Assumption $(A_2)$ holds only for $\gamma\in[\frac{1}{2}, 1)$, the pair $[A,B]$ may  be
  ill-posed in $X$.
\end{enumerate}
\end{remark}
    Inspired by $(ii)$ in Remark \ref{yu-remark-6-23-2}, in order to incorporate the study of the pair $[A,B]$ into the framework in \cite{Kunisch-Wang-Yu, Liu-Wang-Xu-Yu}, we define
\begin{equation}\label{yu-25-6-23-5}
    \mathcal{A}:=
\begin{cases}
    A&\mbox{if}\;\;\gamma\in[0,\frac{1}{2}),\\
    \widetilde{A}_{-\frac{1}{2}}&\mbox{if}\;\;\gamma\in[\frac{1}{2},1),
\end{cases}\;\;\mbox{with}\;\;
    D(\mathcal{A})=
  \begin{cases}
   X_1&\mbox{if}\;\;\gamma\in[0,\frac{1}{2}),\\
    X_{\frac{1}{2}}&\mbox{if}\;\;\gamma\in[\frac{1}{2},1),
  \end{cases}
\end{equation}
   and
\begin{equation}\label{yu-25-6-23-5-b}
\mathcal{X}:=\begin{cases}
        X&\mbox{if}\;\;\gamma\in[0,\frac{1}{2}),\\
        X_{-\frac{1}{2}}&\mbox{if}\;\;\gamma\in[\frac{1}{2},1),
    \end{cases}
\end{equation}
    where $\widetilde{A}_{-\frac{1}{2}}$ is given in \eqref{yu-25-6-20-1}.  In Section \ref{yu-25-sec-1-2-1}, we will show  the following proposition.
\begin{proposition}\label{yu-proposition-6-27-1}
    Under Assumptions $(A_1)$ and $(A_2)$, the following statements hold:
\begin{enumerate}
  \item [$(i)$] The operator $\mathcal{A}$ with domain $D(\mathcal{A})$ generates an analytic semigroup $\mathcal{S}(\cdot)$ on $\mathcal{X}$.
  \item [$(ii)$] The control operator $B$ is admissible (w.r.t. $\mathcal{A}$)
    in $\mathcal{X}$.
  \item [$(iii)$] The pair $[\mathcal{A},B]$ is well-posed in $\mathcal{X}$, i.e., for any $\mathcal{Y}_0\in\mathcal{X}$ and $u\in L^2_{loc}(\mathbb{R}^+;U)$, the equation (in $\mathcal{X}$)
  \begin{equation}\label{yu-25-6-27-1}
  \begin{cases}
    \mathcal{Y}_t(t)=\mathcal{A}\mathcal{Y}(t)+Bu(t),\;\;t\in\mathbb{R}^+,\\
    \mathcal{Y}(0)=\mathcal{Y}_0
  \end{cases}
  \end{equation}
    has a unique solution $\mathcal{Y}(\cdot;\mathcal{Y}_0,u)\in C([0,+\infty);\mathcal{X})$.
\end{enumerate}
\end{proposition}
\par
    In the rest of this subsection, we always suppose that
    Assumptions $(A_1)$ and $(A_2)$ hold.
     Based on Proposition \ref{yu-proposition-6-27-1}, the exponential stabilizability property for the pair  $[\mathcal{A},B]$ (in $\mathcal{X}$) is defined as follows.
\begin{definition}\label{yu-definition-10-18-1}
    \begin{enumerate}
      \item [$(i)$] Let $\alpha\in\mathbb{R}^+$.  The pair $[\mathcal{A},B]$ is said to be
  $\alpha$-stabilizable in $\mathcal{X}$  if there exists a $C_0$-semigroup $\mathcal{T}(\cdot)$ on $\mathcal{X}$ (of generator $M: D(M)\subset \mathcal{X}\to \mathcal{X}$) and an operator
  $K\in \mathcal{L}(D(M);U)$ so that
\begin{enumerate}
  \item [$(a)$] there exists $C_1\geq 1$ so that $\|\mathcal{T}(t)\|_{\mathcal{L}(\mathcal{X})}
  \leq C_1e^{-\alpha t}$ for all  $t\in\mathbb{R}^+$;
\item[$(b)$] for any $x\in D(M)$,  $M x=\widetilde{\mathcal{A}}x+BKx$;
  \item [$(c)$] there exists $C_2\geq 0$ so that $\|K\mathcal{T}(\cdot)x\|_{L^2(\mathbb{R}^+;U)}
  \leq C_2\|x\|_\mathcal{X}$ for any $x\in D(M)$,
\end{enumerate}
   where $\widetilde{\mathcal{A}}$ is the extension of $\mathcal{A}$ from $\mathcal{X}$ to $\mathcal{X}_{-1}$
   \footnote{The space $\mathcal{X}_{-1}$ is the dual space of $D(\mathcal{A}^*):=\{\varphi\in\mathcal{X}:\exists\;C>0\;\mbox{s.t.}\;|\langle \varphi,\mathcal{A}\psi\rangle_{\mathcal{X}}|\leq C\|\psi\|_{\mathcal{X}}\;\forall\;\psi\in D(\mathcal{A})\}$  {\color{blue}w.r.t.} the pivot space $\mathcal{X}$. Here and in what follows, we denote by $\mathcal{A}^*$  the adjoint operator of $\mathcal{A}$ in the sense: $\langle\mathcal{A}\varphi,\psi\rangle_{\mathcal{X}}=
   \langle\varphi,\mathcal{A}^*\psi\rangle_{\mathcal{X}}$ for any $\varphi\in D(\mathcal{A})$ and $\psi\in D(\mathcal{A}^*)$ and let $\mathcal{X}_1^*:=D(\mathcal{A}^*)$.} (see
    $(ii)$ in Remark \ref{yu-remark-6-22-1}). The operator  $K$ is  called a feedback law to the $\alpha$-stabilizability property for the pair $[\mathcal{A},B]$ in $\mathcal{X}$.
      \item [$(ii)$] The pair $[\mathcal{A},B]$ is said to be stabilizable in $\mathcal{X}$ if there exists $\alpha\in\mathbb{R}^+$ so that $[\mathcal{A},B]$ is $\alpha$-stabilizable in $\mathcal{X}$. The pair $[\mathcal{A},B]$ is said to be rapidly stabilizable in $\mathcal{X}$ if $[\mathcal{A},B]$ is $\alpha$-stabilizable in $\mathcal{X}$ for any $\alpha\in\mathbb{R}^+$.
    \end{enumerate}
\end{definition}
\begin{remark}\label{yu-remark-11-20-111}
\begin{enumerate}
    \item[$(i)$] The concept of stabilizability given in Definition \ref{yu-definition-10-18-1}  was firstly given in \cite{Liu-Wang-Xu-Yu} (see also \cite{Kunisch-Wang-Yu, Ma-Wang-Yu}). The main idea to define the stabilizability by $(i)$ in Definition \ref{yu-definition-10-18-1} originally comes from the LQ-theory of the infinite-dimensional linear control systems established in \cite{Flandoli-Lasiecka, Lasiecka-Triggiani, Lasiecka-Triggiani-2000}, which shows that the finite cost
 condition of the standard LQ problem: $\int_{\mathbb{R}^+}[\|\mathcal{Y}(t;\mathcal{Y}_0,u)\|_{\mathcal{X}}^2+\|u(t)\|_U^2]dt<+\infty$
for all $(\mathcal{Y}_0,u)\in \mathcal{X}\times L^2(\mathbb{R}^+;U)$ implies the exponential stabilizability in the sense of Definition \ref{yu-definition-10-18-1}, where
$\mathcal{Y}(\cdot;\mathcal{Y}_0,u)$ is the solution of the equation \eqref{yu-25-6-27-1}. Indeed, they are equivalent (see \cite[Proposition 3.9]{Liu-Wang-Xu-Yu}). Moreover, if the pair $[A,B]$ is stabilizable, then the feedback law can be constructed by considering the standard LQ problem.
\item[$(ii)$] The admissibility of $B$ (w.r.t. $\mathcal{A}$) in $\mathcal{X}$ plays a key role in Definition \ref{yu-definition-10-18-1}. Indeed, as stated in $(i)$ above, Definition \ref{yu-definition-10-18-1} is based on the LQ problem studied in \cite{Flandoli-Lasiecka, Lasiecka-Triggiani}, which requires the controlled system to be well-posed in $C([0,+\infty);\mathcal{X})$ (see Assumption $(H_1)$ in \cite{Flandoli-Lasiecka}  or Assumption $(H_2)$ in \cite[Chapter 1]{Lasiecka-Triggiani}).
\item[$(iii)$] The feedback law in Definition \ref{yu-definition-10-18-1} may be unbounded (i.e., $K\notin \mathcal{L}(\mathcal{X};U)$). For example, we can use boundary damping to stabilize wave equations (see, for instance, \cite[Section 8.6, Chapter 8]{Komornik}).
    In order to establish the stabilizability property for the pair
$[\mathcal{A},B]$ within the framework of $(i)$ in Definition \ref{yu-definition-10-18-1}, we need  to find a  feedback law $K$,  a closed-loop system whose solution has a $C_0$-semigroup representation (i.e., the form $\mathcal{T}(\cdot)\mathcal{Y}_0$ for the initial data $\mathcal{Y}_0$), and then prove that $\mathcal{T}(\cdot)$ is exponentially stable (i.e., the condition $(a)$).
Once the semigroup $\mathcal{T}(\cdot)$ is obtained, its generator $M$ and $D(M)$ can be defined directly.
In order to test the condition $(b)$, $D(M)$ should be a subset of $D(K)$ (when the feedback law is  bounded, its domain is $X$, however, when it is unbounded, its domain is a
proper subset of $X$). In general, the domain of the feedback law and the domain of the generator are mutually determined. The condition $(c)$ is the admissibility condition of $K$ to $\mathcal{T}(\cdot)$. Indeed, the condition $(c)$ ensures that the stabilized closed-loop system can be also placed within  the framework of the pair $[\mathcal{A},B]$, where the admissible control set is $L^2_{loc}(\mathbb{R}^+;U)$.
\item[$(iv)$] In some cases, the pair $[\mathcal{A},B]$ is also said to be
  $\alpha$-stabilizable with feedback law $K$ in order to emphasize that $K$ is obtained and fixed.
  \item[$(v)$] In general, for unbounded control operator $B$, $D(M)$ is not equal to $D(\mathcal{A})$.
\end{enumerate}
\end{remark}
    We next define the exponential stabilizability property for the pair $[A,B]$ in $X$. Since Assumption $(A_2)$ does not imply the admissibility of $B$ (w.r.t. $A$) in $X$ (see $(ii)$ in Remark \ref{yu-remark-6-23-2}),
   the framework in  \cite{Kunisch-Wang-Yu, Liu-Wang-Xu-Yu, Ma-Wang-Yu} is no longer applicable to the pair $[A,B]$ under Assumption $(A_2)$. However, by the well-known Banach fixed-point theorem, one can easily show  the following lemma.
\begin{lemma}\label{yu-proposition-25-6-26-2}
    Suppose that Assumptions $(A_1)$ and $(A_2)$ hold. If the feedback law  $K$ is bounded (i.e., $K\in\mathcal{L}(X;U)$), then for any $y_0\in X$, the equation
\begin{equation}\label{yu-25-6-30-1}
\begin{cases}
    y_t(t)=Ay(t)+BKy(t),&t\in\mathbb{R}^+,\\
    y(0)=y_0
\end{cases}
\end{equation}
     has a unique solution in $C([0,+\infty);X)$.
\end{lemma}
    Inspired by Lemma \ref{yu-proposition-25-6-26-2}, we define the exponential/rapid stabilizability property for the pair $[A,B]$ in $X$ as follows.
\begin{definition}\label{yu-def-9-27-1}
\begin{enumerate}
      \item [$(i)$] Let $\alpha\in\mathbb{R}^+$. The pair $[A,B]$ is said to be  $\alpha$-stabilizable in $X$
       with bounded feedback law if there exists $K\in\mathcal{L}(X;U)$ so that for any $y_0\in X$, the solution
       $y_K(\cdot;y_0)$ of the closed-loop system \eqref{yu-25-6-30-1} is exponentially stable with decay rate $\alpha$, i.e., there exists $C(\alpha)>0$ so that $\|y_K(t;y_0)\|_X\leq C(\alpha)e^{-\alpha t}\|y_0\|_X$ for each $t\in\mathbb{R}^+$. The operator  $K$ is called a feedback law to the $\alpha$-stabilizability property for the  pair $[A,B]$ in $X$.

         When  the feedback law $K$ is given and fixed,  the pair $[A,B]$ is also said to be
  $\alpha$-stabilizable in $X$ with feedback law $K$.
      \item [$(ii)$] The pair $[A,B]$ is said to be stabilizable in $X$ with bounded feedback law if there exists $\alpha\in\mathbb{R}^+$ so that $[A,B]$ is $\alpha$-stabilizable in $X$ with bounded feedback law.
      \item [$(iii)$] The pair $[A,B]$ is said to be rapidly stabilizable in $X$ with bounded feedback law if $[A,B]$ is $\alpha$-stabilizable in $X$ with bounded feedback law for any $\alpha>0$.
    \end{enumerate}
\end{definition}

\subsection{Problems and motivation}\label{yu-sec-1-3-1}
\vskip 5pt
 In this paper, we investigate the following two issues:
\begin{enumerate}
  \item [$(P1)$] Provide a relationship between the  stabilizability property for the pair $[A,B]$ with bounded feedback law and the  stabilizability property for the pair $[\mathcal{A},B]$.
  \item [$(P2)$] Give a frequency-domain characterization on the  stabilizability  property for the pair $[A,B]$ with bounded feedback law.
\end{enumerate}
\vskip 5pt

   We first present the motivation to study the problem $(P1)$. As stated in Remark \ref{yu-remark-11-20-111}, a new definition (i.e., Definition \ref{yu-definition-10-18-1})  of exponential stabilizability of the control system within the framework of admissible control operators
   was given in \cite{Liu-Wang-Xu-Yu} and  many new results in this framework were obtained (see \cite{Liu-Wang-Xu-Yu, Kunisch-Wang-Yu, Ma-Wang-Yu}). In the setting of unbounded control operators, this framework is one of many existing frameworks.
      Indeed, to our knowledge, besides the framework in \cite{Liu-Wang-Xu-Yu}, there are many other frameworks used to study the stabilizability of control systems with unbounded control operators (see, for instance, \cite{Badra-Takahashi, Curtain-Weiss, Raymond, Staffans-1998, Weiss-Rebarber} and so on). In particular, the framework in \cite{Badra-Takahashi, Raymond}
      does not need the admissibility of control  operators. In practice, there are many examples in which control  operators are not admissible, but the systems are feedback stabilizable, for instance, the heat equation with Dirichlet boundary controls (see \cite{Badra-Takahashi} or Section \ref{yu-sec-25-7-1} in this paper). Therefore, it is natural to study \emph{whether so many frameworks are equivalent or not}.
      This is a very important problem in control theory of the distributed parameter systems. Unfortunately, we do not find literature touching on such issues.
      As the first purpose of this paper, we try to unify the framework of \cite{Badra-Takahashi, Raymond} and that of \cite{Liu-Wang-Xu-Yu} within an abstract setting (i.e., the problem $(P1)$). In this paper, we will show that under Assumptions $(A_1)$, $(A_2)$ and an almost exponential decay condition (see Definition \ref{yu-def-7-2-1} below)  for the state operator, the stabilizability  property for the pair $[\mathcal{A},B]$ is equivalent to that of the pair $[A,B]$ with bounded feedback law.

      \medskip We next give the motivation to study the problem $(P2)$. In \cite{Trelat-Wang-Xu}, the authors introduced a weak observability inequality (a kind of time-domain condition) to characterize the stabilizability property for infinite-dimensional control systems for bounded control  operators or under Assumptions $(A_1)$ and $(A_2)$ for $\gamma\in[0,\frac{1}{2})$. This result has been extended to the case of admissible  control operators (see \cite{Liu-Wang-Xu-Yu, Ma-Wang-Yu}). It is worth noting that the admissibility property plays an important role in the weak observability inequality. In fact, if the control  operator is not admissible, then the weak observability inequality is meaningless
      (see \cite[Sections 4.2 and 4.3, Chapter 4]{Tucsnak-Weiss}).
      However, as we mention in $(ii)$ of Remark \ref{yu-remark-6-23-2}, the control operator $B$ may be not admissible under Assumption $(A_2)$. Therefore, it is inappropriate to use the weak observability inequality to characterize the stabilizability  property for the pair $[A,B]$. Recently, in \cite{Kunisch-Wang-Yu}, the authors obtained the frequency-domain characterizations on the stabilizability property for infinite-dimensional control systems with admissible control operators under some assumptions on the state operator $A$, for instance, the almost exponential decay condition (see Definition \ref{yu-def-7-2-1} below). When we strengthen the requirements for the state operator $A$ (for example, requiring it to generate an analytic semigroup), whether we can generalize the results in \cite{Kunisch-Wang-Yu} or not (in particular,  Theorem 2 in \cite{Kunisch-Wang-Yu}) to the case that the control operator is not admissible  is an important
      and interesting problem.

\subsection{Main result and novelties}
\begin{definition}\label{yu-def-7-2-1}
    An operator $L$: $D(L)\to X$ is said to satisfy the almost exponential decay condition ((AEDC) for short) in $X$ if the following conditions hold:
\begin{enumerate}
  \item [$(i)$] The linear operator $L$ with domain $D(L)$ generates a $C_0$-semigroup $S_L(\cdot)$ in $X$.
  \item [$(ii)$] For each $\alpha>0$, there are two closed subspaces $Q_1:=Q_1(\alpha)$ and $Q_2:=Q_2(\alpha)$ of $X$ so that $(a)$ $X=Q_1\oplus Q_2$;  $(b)$ $Q_1$ and $Q_2$ are  invariant subspaces of $S_L(\cdot)$;
    $(c)$ $L|_{Q_1}$ (the restriction of $L$ on $Q_1$) is bounded and satisfies that  $\sigma(L|_{Q_1})\subset\mathbb{C}_{-\alpha}^+$; $(d)$  $S_L(\cdot)|_{Q_2}$ (the restriction of $S_L(\cdot)$ on $Q_2$) is  exponentially stable.
\end{enumerate}
\end{definition}
\begin{remark}\label{yu-remark-12-22-1}
    Several remarks on Definition \ref{yu-def-7-2-1} are given in order.
\begin{enumerate}
  \item [$(i)$] The assumption that $A^*$ satisfies (AEDC) in $X$ has been used in \cite{Kunisch-Wang-Yu}
    to characterize the stabilizability property for  the pair $[A,B]$ in the viewpoint of frequency-domain.
  \item[$(ii)$] $L$ satisfies (AEDC) in $X$ if and only if $L^*$ satisfies (AEDC) in $X$
  (see Lemma \ref{yu-lemms-25-7-14-1}).
    \item[$(iii)$] The spaces $Q_1$ and $Q_2$ are respectively called the unstable and stable parts of state space $X$. In  $(ii)$ of Definition \ref{yu-def-7-2-1},
    the condition that $L|_{Q_1}$ is bounded means that there exists $C>0$ so that $|\langle L^*f,\varphi\rangle_X|\leq C\|\varphi\|_X\|f\|_X$ for any $\varphi\in Q_1$ and $f\in D(L^*)$. It implies that $Q_1\subset D(L)$. The condition $\sigma(L|_{Q_1})\subset\mathbb{C}_{-\alpha}^+$ is equivalent to that, for each $\lambda\in \overline{\mathbb{C}_{-\alpha}^-}$ and any $f\in Q_1$, there exists a unique $\varphi\in Q_1$ so that $(\lambda I-L)\varphi=f$.
        The condition $(d)$  means that there exist $\varepsilon>0$ and $C(\varepsilon)>0$ so that $\|S_L(t)f\|_X\leq C(\varepsilon)e^{-\varepsilon t}\|f\|_X$
        for each $f\in Q_2$.
  \item [$(iv)$] If the linear operator $L$ generates a $C_0$-semigroup $S_L(\cdot)$ and  one of the following assumptions holds:
\begin{enumerate}
  \item [$(a)$] $S_L(\cdot)$ is a compact semigroup on $X$.
  \item [$(b)$] $L$ with domain $D(L)$ is normal and $\sigma(L)\cap \mathbb{C}_{-\gamma}^+$ is bounded for each $\gamma\in\mathbb{R}^+$,
\end{enumerate}
     then $L$ satisfies (AEDC) in $X$ (see Lemma \ref{yu-lemma-12-20-1-bb}). Two notes on $(a)$ and $(b)$ are listed as follows:
     \begin{enumerate}
  \item [$(iv)_1$] Throughout this paper, when we say $S_L(\cdot)$ is compact, it means that there exists $t_0\geq 0$ so that $S_L(t)$ is compact for any $t>t_0$, i.e., $S_L(\cdot)$ is immediately or eventually compact (see \cite[Definition 4.23, Section 4, Chapter II]{Engel-Nagel}).
  \item [$(iv)_2$] If $L$ is a self-adjoint operator or a linear partial differential operator on the whole
   space with constant coefficients, then $L$ is normal (see \cite[Theorem 13.24, Chapter 13]{Rudin}).
\end{enumerate}
\end{enumerate}

\end{remark}

    We make the following assumption:
\begin{enumerate}
  \item [$(A_3)$] The operator $A$ satisfies (AEDC) in $X$.
\end{enumerate}
\vskip 5pt
    The main result of this paper is stated as follows.
\begin{theorem}\label{yu-theorem-3-31-1}
    Under Assumptions $(A_1)$, $(A_2)$ and $(A_3)$, the following statements are equivalent:
\begin{enumerate}
  \item [$(i)$] The pair $[\mathcal{A},B]$ is stabilizable in $\mathcal{X}$.
  \item [$(ii)$] The pair $[\mathcal{A},B]$ is stabilizable in $\mathcal{X}$ with bounded feedback law.
  \item[$(iii)$] The pair $[A,B]$ is stabilizable in $X$ with bounded feedback law.
\item[$(iv)$]  There exist $\alpha>0$ and $C(\alpha)>0$ so that
\begin{equation}\label{yu-4-2-2-b}
    \|\varphi\|^2_X\leq C(\alpha)(\|(\lambda I-A^*)\varphi\|_X^2+\|B^*\varphi\|_U^2)\;\;\mbox{for any}\;\;\lambda\in\mathbb{C}_{-\alpha}^+,\;\;\varphi\in X_1^*.
\end{equation}
\end{enumerate}
\end{theorem}
\begin{remark}\label{yu-remark-4-8-1}
    Some comments on Theorem \ref{yu-theorem-3-31-1} are listed in order.
\begin{enumerate}
  \item [$(i)$] The inequality \eqref{yu-4-2-2-b} is a Hautus-type resolvent inequality or a frequency-domain stabilizability test, which was used in \cite{Kunisch-Wang-Yu} to characterize the stabilizability property for the pair $[A,B]$ when $B$ is admissible. Theorem \ref{yu-theorem-3-31-1} shows that under Assumptions $(A_1)$, $(A_2)$ and $(A_3)$, \eqref{yu-4-2-2-b} can still be used to characterize the stabilizability property for the pair $[A,B]$
      even if $B$ may fail to be admissible.
  \item[$(ii)$] By Proposition \ref{yu-proposition-6-27-1} and \cite[Theorem 25]{Trelat-Wang-Xu}, we  conclude that the pair $[\mathcal{A},B]$ is stabilizable in $\mathcal{X}$ if and only if there exist $T>0$, $\delta\in[0,1)$ and $C(T,\delta)>0$ so that
\begin{equation*}\label{yu-4-23-1-bb}
    \|\mathcal{S}^*(T)\varphi\|_{\mathcal{X}}^2\leq C(T,\delta)\int_0^T\|\mathcal{B}^*\mathcal{S}^*(t)\varphi\|_U^2dt
    +\delta\|\varphi\|_{\mathcal{X}}^2\;\;\mbox{for any}\;\;\varphi\in \mathcal{X}^*_1,
\end{equation*}
    which is a weak observability inequality  for the pair $[\mathcal{A}^*,\mathcal{B}^*]$. Here,
    $\mathcal{B}^*$ is the adjoint operator of $B$ w.r.t. the pivot space $\mathcal{X}$, i.e.,
    $\langle Bu,\varphi\rangle_{\mathcal{X}_{-1},\mathcal{X}_1^*}=\langle u, \mathcal{B}^*\varphi\rangle_U$ for any $u\in U$ and $\varphi\in \mathcal{X}_1^*$. Its expression is given in \eqref{yu-3-29-14} below.
\item[$(iii)$] Under Assumptions $(A_1)$ and $(A_2)$, when
$S(\cdot)$ is a compact semigroup (which implies that $(A_3)$ holds by $(iv)$ in Remark \ref{yu-remark-12-22-1}), $(iv)$ is equivalent to that, there exists $\alpha>0$ so that for each $\lambda\in\mathbb{C}_{-\alpha}^+$, if there exists $\varphi\in X_1^*$ so that
  $(\lambda I-A^*)\varphi=0$ and $B^*\varphi=0$, then $\varphi=0$ (see Corollary \ref{yu-proposition-4-22-1}).
  We thus recover partially the results on stabilizability of \cite{Badra-Takahashi, Raymond}. Specifically, Theorem \ref{yu-theorem-3-31-1} explains why in the setting of \cite{Badra-Takahashi, Raymond}, it suffices to search the feedback law  in the family of bounded linear operators.

  \item[$(iv)$] The equivalence between $(i)$ and $(ii)$ (of Theorem \ref{yu-theorem-3-31-1})  is an independent result.  It does not require $A$ to generate an analytic semigroup. For more general results, one can refer to Theorem \ref{yu-theorem-10-18-1}.
  \item[$(v)$] In Theorem \ref{yu-theorem-3-31-1}, the equivalence between  $(i)$ and $(iii)$ provides an answer to the problem
  $(P1)$, and the equivalence between  $(iii)$ and $(iv)$ provides an answer to the problem
  $(P2)$.
\end{enumerate}
\end{remark}

    The following corollary is an application of Theorem \ref{yu-theorem-3-31-1}.
\begin{corollary}\label{yu-corollary-4-22-1}
    Under Assumptions $(A_1)$ and $(A_2)$, if
    the operator $A$ satisfies one of the conditions $(a)$ and $(b)$ in
    $(iv)$ of Remark \ref{yu-remark-12-22-1} (by replacing $L$ by $A$),
    then the following statements are equivalent:
\begin{enumerate}
  \item [$(i)$] The pair $[\mathcal{A},B]$ is rapidly stabilizable in $\mathcal{X}$.
  \item [$(ii)$] The pair $[\mathcal{A},B]$ is rapidly  stabilizable in $\mathcal{X}$ with bounded feedback law.
  \item[$(iii)$] The pair $[A,B]$ is rapidly stabilizable in $X$ with bounded feedback law.
\item[$(iv)$]  For each $\alpha>0$, there exists $C(\alpha)>0$ so that \eqref{yu-4-2-2-b} holds.
\end{enumerate}
\end{corollary}
\begin{remark}\label{yu-remark-7-3-1}
    Several comments on Corollary \ref{yu-corollary-4-22-1} are given in order.
\begin{enumerate}
  \item [$(i)$] By Proposition \ref{yu-proposition-6-27-1} and \cite[Theorem 3.4]{Liu-Wang-Xu-Yu}, we  conclude that
  the pair $[\mathcal{A},B]$ is rapidly stabilizable in $\mathcal{X}$ if and only if for each $\alpha>0$, there exist two positive constants $C(\alpha)$ and $D(\alpha)$ so that
\begin{equation*}\label{yu-25-7-3-1}
     \|\mathcal{S}^*(t)\varphi\|_{\mathcal{X}}^2\leq D(\alpha)\int_0^t\|\mathcal{B}^*\mathcal{S}^*(s)\varphi\|_U^2ds
    +C(\alpha)e^{-\alpha s}\|\varphi\|_{\mathcal{X}}^2\;\;\mbox{for any}\;\;\varphi\in \mathcal{X}^*_1\;\;\mbox{and}\;\;
    t\in\mathbb{R}^+,
\end{equation*}
    where $\mathcal{B}^*$ is the adjoint operator $B$ w.r.t. the space $\mathcal{X}$ (see $(ii)$ in Remark \ref{yu-remark-4-8-1}). Corollary \ref{yu-corollary-4-22-1} presents a frequency-domain characterization (i.e., \eqref{yu-4-2-2-b}) on the rapid stabilizability property for the pair $[\mathcal{A},B]$ and the pair $[A,B]$.
    The most intuitive application is that once the pair $[\mathcal{A},B]$ is null controllable in
    the state space $\mathcal{X}$, the pair $[A,B]$ is rapidly stabilizable in $X$ with bounded feedback law.
    This also explains why many systems require different state spaces when studying controllability and stabilizability. For instance, the heat equation with Dirichlet controls in $\Omega\subset\mathbb{R}^n$ is exactly null controllable
     in $H^{-1}(\Omega)$, but is stabilizable in $L^2(\Omega)$ (see Section \ref{yu-sec-25-7-1}).
  \item [$(ii)$] Similar to $(iii)$ of Remark \ref{yu-remark-4-8-1}, we have that if Assumptions $(A_1)$ and $(A_2)$ are true and
$S(\cdot)$ is a compact semigroup, then $(iv)$ in Corollary \ref{yu-corollary-4-22-1} is equivalent to
the following statement:  for any $\alpha>0$ and $\lambda\in\mathbb{C}_{-\alpha}^+$, if there exists  $\varphi\in X_1^*$ so that
  $(\lambda I-A^*)\varphi=0$ and $B^*\varphi=0$, then $\varphi=0$.
  \item [$(iii)$] Similar to $(iv)$ of Remark \ref{yu-remark-4-8-1}, the equivalence between $(i)$ and $(ii)$ of Corollary \ref{yu-corollary-4-22-1} does not need that $A$ generates an analytic semigroup on $X$ (see Corollary \ref{corollary-yu-12-4-1} for general cases).
\end{enumerate}
\end{remark}
\noindent\textbf{Novelties.}  The novelties of this paper are the followings.

\medskip
   $(i)$ We provide an explanation on the necessity of choosing  different state spaces in the study of controllability and stabilizability of the pair $[A,B]$, as well as their relationship (see $(i)$ of Remark \ref{yu-remark-7-3-1}).

  \medskip $(ii)$ We establish a frequency-domain characterization for the exponential  stabilizability property
  for the pair $[A,B]$ under Assumptions $(A_1)$-$(A_3)$.

  \medskip  $(iii)$ We extend  the main results of \cite{Badra-Takahashi, Fattorini-1966, Fattorini-1967, Raymond} to the case where the $C_0$-semigroup generated by the state operator is not compact. We also extend partially \cite[Theorem 2]{Kunisch-Wang-Yu} to the case where the control operator is not admissible.

\subsection{Comments on related works}
 We now recall the related works of this paper from three aspects.
     The first aspect is the related works on analytic-type control systems (i.e., their state operators generate analytic semigroups). The second aspect is the related works on the frameworks of stabilizability of infinite-dimensional control system with unbounded control operator. The third   aspect concerns
     the choice of the family of feedback laws and the algorithmic procedure to construct the feedback.

  \medskip Analytic-type partial differential equations exist widely in physical processes such as heat conduction and reaction-diffusion. Therefore, it is important and meaningful to study how to stabilize such systems
      with internal or boundary controls. When the state operator in this kind of control system also generates  a compact semigroup, the unstable part (w.r.t. the state operator) of the state space is only finite-dimensional (see Lemma \ref{yu-lemma-6-23-1}) and then the corresponding  stabilization problem  is essentially a stabilization problem for a finite dimensional system. As far as we know, this idea was firstly introduced by H. O. Fattorini in \cite{Fattorini-1966, Fattorini-1967} for the controlled parabolic equation with the bounded control operators (for simplicity, it is called the Fattorini's method later). In the past two to three decades, Fattorini's method has been widely applied to the linear or nonlinear parabolic systems with unbounded control operators (see for instance, \cite{Badra-Mitra-Raymond, Badra-Takahashi, Coron-Trelat-2004, Coron-Trelat-2006, Lhachemi-Prieur-Trelat, Raymond, Raymond-Thevenet, Russell}). In particular, in \cite{Badra-Takahashi}, the authors
      showed that the state spaces should be different for controllability and stabilizability.
      In practice, there are many parabolic systems whose state operators do not generate compact semigroups, for example,  reaction-diffusion equations on whole
      space or half of whole space with constant-valued coefficients. In these cases, even if their state operators satisfy
      (AEDC), the unstable part of the state space is also infinite-dimensional. Hence, the existing results
       (when state operator generates a compact semigroup) cannot be directly applied to these non-compact situations. Recently, in \cite{Kunisch-Wang-Yu},  the authors presented a frequency-domain condition for the stabilizability of infinite-dimensional control systems with admissible control operators under Assumption $(A_3)$ (see \cite[Theorem 2]{Kunisch-Wang-Yu}). The key point for proving \cite[Theorem 2]{Kunisch-Wang-Yu} is
       to use  characterizations of controllability of infinite-dimensional control system with bounded
       state and control operators (see Lemma \ref{yu-lemma-12-9-1}).

  \medskip Feedback stabilizability of infinite-dimensional control systems has a research history of several decades. With the efforts of outstanding researchers, many important and classic  results have been established on this topic (see, for instance, \cite{Alabau-Boussouira, Badra-Mitra-Raymond, Badra-Takahashi, Coron-Trelat-2004, Coron-Trelat-2006, Fattorini-1966, Fattorini-1967, Jacob-Zwart, Lhachemi-Prieur-Trelat, Liu, Raymond, Raymond-Thevenet, Russell, Triggiani-1980, Triggiani}). In practice, we often use boundary controls or pointwise controls to stabilize an infinite dimensional system. Therefore, the feedback stabilization problem with unbounded control operators has always been a focus and a difficult point in feedback stabilizability of infinite-dimensional control systems. However, providing the most appropriate definition for  the feedback stabilizability with unbounded control operators is not an easy task.
      As far as we know, in addition to the Fattorini's framework introduced above, there are at least two other frameworks for defining the feedback stabilizability of infinite-dimensional control systems. The first one is  the framework established by R. F. Curtain, G. Weiss, O. J. Staffans, et al. based on some auxiliary problems such as admissibility of control operator, well-posed and regular linear system and so on (see, for example, \cite{Curtain-Weiss, Staffans-1998, Weiss-Rebarber} and references therein). The second one is the framework established in \cite{Liu-Wang-Xu-Yu, Trelat-Wang-Xu} under the conclusions in the study of LQ problem of infinite-dimensional control system on infinite horizon (see $(i)$ in Remark \ref{yu-remark-11-20-111}). To our knowledge, the relationship between these two frameworks is still an open issue.

  \medskip As one of the characteristics of feedback stabilization problems, we not only need to answer whether the system can be stabilizable or not, but also provide an algorithmic procedure to construct the feedback in some family of linear operators,
     for example, Gramian method, Fattorini's method, Riccati method, back-stepping method and so on. For related works, one can refer to \cite{Badra-Takahashi, Coron-Lv, Fattorini-1966, Flandoli-Lasiecka, Komornik-1997, Vest} and so on.
      When the control operator is unbounded, except for the spectral projection method, the feedback laws designed by these methods are generally unbounded. As far as we know, there is no literature that can provide a precise characterization (or description) of the set formed by these feedback laws, even whether this set is a linear space or not is unknown.

\subsection{Structure of the paper}
    This paper is organized as follows. In Section \ref{yu-sec-1-2-1}, we present some preliminary results. In Section \ref{yu-sec-25-3}, we establish the existence of bounded feedback laws under more general assumptions. In fact, Theorem \ref{yu-theorem-10-18-1} in Section \ref{yu-sec-25-3} covers the equivalence between $(i)$ and $(ii)$ of Theorem \ref{yu-theorem-3-31-1}. In Section \ref{yu-sect-5-12-1}, we prove our main result. Some examples are presented in Section \ref{yu-sec-7-15-100}. Section \ref{yu-sec-6} is the appendix of this paper.

\section{Preliminary results}\label{yu-sec-1-2-1}
\subsection{Stabilizability with bounded control operators}\label{yu-sec-2-1-1}
    In this subsection, we recall characterizations of controllability and stabilizability properties for the pair $[A,B]$ under the assumption that $A$ and $B$ are bounded.
    The following results can be found in  \cite{Takahashi, Triggiani, Tucsnak-Weiss} and are collected in \cite[Lemma 5]{Kunisch-Wang-Yu}.
\begin{lemma}\label{yu-lemma-12-9-1}
    Suppose that $A\in\mathcal{L}(X)$ and $B\in\mathcal{L}(U;X)$. The following statements are equivalent:
\begin{enumerate}
  \item [$(i)$] The pair $[A,B]$ is exactly controllable for some $T>0$, i.e., there exists $T>0$ so that for any $y_0, y_1\in X$, there  exists $u\in L^2(0,T;U)$ so that $y(T;y_0,u)=y_1$.
  \item [$(ii)$] For each $T>0$, there  exists $C(T)>0$ so that
$\|\varphi\|_X^2\leq C(T)\int_0^T\|B^*S^*(t)\varphi\|_U^2dt$ for any $\varphi\in X$.
  \item [$(iii)$] For each $\lambda\in\mathbb{C}$, there  exists $C(\lambda)>0$ so that
  $\|\varphi\|_X^2\leq C(\lambda)(\|(\lambda I-A^*)\varphi\|^2_X+\|B^*\varphi\|_U^2)$ for any $\varphi\in X$.
\end{enumerate}
\end{lemma}
    Suppose that $A\in\mathcal{L}(X)$, $B\in\mathcal{L}(U;X)$ and $[A,B]$ is exactly controllable. By Lemma \ref{yu-lemma-12-9-1} and \cite[Proposition 3]{Ma-Wang-Yu},
    we have that for each $\alpha>0$, there exist $D_1(\alpha)>0$ and $D_2(\alpha)\geq1$ so that

\begin{equation}\label{yu-12-9-2}
    \|S^*(t)\varphi\|_X^2\leq D_1(\alpha)\int_0^t\|B^*S^*(s)\varphi\|^2_Uds
    +D_2(\alpha)e^{-\alpha t}\|\varphi\|_X^2
    \;\;\mbox{for any}\;\;t\in\mathbb{R}^+\;\;\mbox{and}\;\;\varphi\in X.
\end{equation}
    Then it follows from \cite[Theorem 1.1]{Liu-Wang-Xu-Yu} that the pair $[A,B]$ is rapidly stabilizable in $X$. Let $T>0$, $\alpha>0$ and $\varepsilon\geq0$ be fixed arbitrarily. Let $D_1(4\alpha)>0$ and $D_2(4\alpha)\geq 1$ be given so that \eqref{yu-12-9-2} holds by replacing $\alpha$
    by $4\alpha$. We define
\begin{eqnarray*}\label{yu-12-9-3}
    \langle \Lambda_{\alpha,\varepsilon} (t)\varphi,\psi\rangle_X:=D_1(4\alpha)e^{4\alpha T}\int_0^te^{-(4\alpha-\varepsilon) s}
    \langle B^*S^*(-s)\varphi,B^*S^*(-s)\psi\rangle_Uds\nonumber\\
    +D_2(4\alpha)e^{-(4\alpha-\varepsilon) t}\langle S^*(-t)\varphi,S^*(-t)\psi\rangle_X\;\;\mbox{for}\;\;t\in[0,T]
\end{eqnarray*}
    and
$$\langle \Pi_{\alpha,\varepsilon,T}\varphi,\psi\rangle_X:=\int_0^T\langle \Lambda_{\alpha,\varepsilon} (t)\varphi,\psi\rangle_Xdt$$
    for any $\varphi,\psi\in X$. It is obvious that $\Lambda_{\alpha,\varepsilon}(t)$ ($t\in[0,T]$) and $\Pi_{\alpha,\varepsilon,T}$ are in $\mathcal{L}(X)$. By \cite[Lemma 2.1]{Ma-Wang-Yu}, we conclude that they are also invertible.
\begin{lemma}\label{yu-lemma-12-10-120}
    Suppose that  $A\in\mathcal{L}(X)$, $B\in\mathcal{L}(X;U)$ and $[A,B]$ is exactly controllable. Fix $\alpha>0$ arbitrarily.  Let $D_1(\alpha)$ and $D_2(\alpha)$ be two positive constants so that \eqref{yu-12-9-2} holds. Let $T\in((2\alpha)^{-1}\ln (D_2(4\alpha)),+\infty)$, $\widehat{\varepsilon}:=T^{-1}\ln[D_2(4\alpha)]$ and
$K_T:=-TD_1(4\alpha)e^{4\alpha T}B^*\Pi_{\alpha,\widehat{\varepsilon},T}^{-1}$.
    Then the $C_0$-group $S_{K_T}(\cdot)$ generated by $A+BK_T$ satisfies that
\begin{equation}\label{yu-12-20-2-b}
\|S_{K_T}(t)\|_{\mathcal{L}(X)}\leq C(\alpha) e^{-\alpha t}\;\;\mbox{for any}\;\;t\in\mathbb{R}^+,
\end{equation}
    where $C(\alpha)>0$ is a constant.
    Moreover, if $X$ is finite-dimensional, then $S_{K_T}(\cdot)$ is immediately compact.
\end{lemma}
\begin{proof}
    Since $A$ is bounded, $S(\cdot)$ is a $C_0$-group on $X$. By \cite[Theorem 1.2]{Ma-Wang-Yu}, we obtain that
\begin{equation}\label{yu-12-28-1}
    \|S_{K_T}(t)\|_{\mathcal{L}(X)}\leq C(\alpha) e^{-\frac{1}{2}(4\alpha-T^{-1}\ln(D_2(4\alpha)))t}
    \;\;\mbox{for any}\;\;t\in\mathbb{R}^+,
\end{equation}
    where $T\in ((4\alpha)^{-1}\ln (D_2(4\alpha)),+\infty)$.
    Note that  $((2\alpha)^{-1}\ln (D_2(4\alpha)),+\infty)\subset ((4\alpha)^{-1}\ln (D_2(4\alpha)),+\infty)$ and $T^{-1}\ln (D_2(4\alpha))<2\alpha$ if $T\in((2\alpha)^{-1}\ln (D_2(4\alpha)),+\infty)$. According to \eqref{yu-12-28-1},
    the inequality \eqref{yu-12-20-2-b} holds. When $X$ is finite-dimensional, $S_{K_T}(\cdot)$ is compact since any bounded set of finite-dimensional space is pre-compact.
\end{proof}
\begin{remark}
    In Lemma \ref{yu-lemma-12-10-120}, we only present a special way to design the feedback law when $A$ generates a $C_0$-group. This method is based on the weak observability inequality. For related works, one can refer to \cite{Komornik-1997, Urquiza, Vest} and so on.
\end{remark}

\subsection{Proof of Proposition \ref{yu-proposition-6-27-1}}\label{yu-25-sec-1-2-1}
    The proof of Proposition \ref{yu-proposition-6-27-1} requires several preliminary lemmas, which are shown in order.
\begin{lemma}\label{yu-lemma-3-29-1}
    Let $\mathcal{A}$ and $D(\mathcal{A})$ be defined by \eqref{yu-25-6-23-5}, and on $\mathcal{X}$,
\begin{equation*}\label{yu-3-29-1}
    \mathcal{S}(\cdot):=
\begin{cases}
    S(\cdot)&\mbox{if}\;\;\gamma\in[0,\frac{1}{2}),\\
    (\rho_0I-\widetilde{A})^{\frac{1}{2}}S(\cdot)(\rho_0I-\widetilde{A})^{-\frac{1}{2}}&\mbox{if}\;\;\gamma\in[\frac{1}{2},1).
\end{cases}
\end{equation*}
     Then
    $\mathcal{S}(\cdot)$ is an analytic semigroup on $\mathcal{X}$ and its generator is $\mathcal{A}$ with domain $D(\mathcal{A})$.
\end{lemma}
\begin{proof}
    One can directly check that $\mathcal{S}(\cdot)$ is an analytic semigroup on $\mathcal{X}$.
    When $\gamma\in[0,\frac{1}{2})$, it is clear that the generator of $\mathcal{S}(\cdot)$ is $A$ with domain $X_1$. When
    $\gamma\in[\frac{1}{2},1)$, we suppose that
    $\widehat{\mathcal{A}}$ is the generator of $\mathcal{S}(\cdot)$.
    On one hand, according to \eqref{yu-3-28-10}, it is clear that for any $f\in X_{-\frac{1}{2}}$,
    $\lim_{t\to 0^+} t^{-1}(\mathcal{S}(t)-I)f$ exists in $X_{-\frac{1}{2}}$ if and only if $\lim_{t\to0^+} t^{-1}(S(t)-I)(\rho_0I-\widetilde{A})^{-\frac{1}{2}}f$ exists in $X$.
    Meanwhile, we observe that $\lim_{t\to 0^+} t^{-1}(S(t)-I)(\rho_0I-\widetilde{A})^{-\frac{1}{2}}f$ exists in $X$
    if and only if $(\rho_0 I-\widetilde{A})^{-\frac{1}{2}}f\in D(A)(=X_1)$. Hence,
    \begin{equation}\label{wang-9-23-1}
    f\in D(\widehat{\mathcal{A}})\;\;\mbox{if and only if}\;\;(\rho_0 I-\widetilde{A})^{-\frac{1}{2}}f\in X_1.
    \end{equation}
    On the other hand, since
\begin{equation}\label{yu-3-29-4}
    (\rho_0I-\widetilde{A})^{-\frac{1}{2}}X_{\frac{1}{2}}= X_1,
\end{equation}
we obtain from (\ref{wang-9-23-1}) that $D(\widehat{\mathcal{A}})=X_{\frac{1}{2}}$. Moreover, by \eqref{yu-3-28-10} and (\ref{yu-3-29-4}),
 we have that for any $f\in D(\widehat{\mathcal{A}})$ and $g\in X_{-\frac{1}{2}}$,
\begin{eqnarray*}\label{yu-3-29-5}
    &\;&\lim_{t\to0^+}t^{-1}\langle (\mathcal{S}(t)-I)f,g\rangle_{X_{-\frac{1}{2}}}
    =\lim_{t\to0^+}t^{-1}\langle (S(t)-I)(\rho_0I-\widetilde{A})^{-\frac{1}{2}}f,(\rho_0I-\widetilde{A})^{-\frac{1}{2}}g\rangle_X\nonumber\\
    &=&\langle A(\rho_0I-\widetilde{A})^{-\frac{1}{2}}f,(\rho_0I-\widetilde{A})^{-\frac{1}{2}}g\rangle_X
    =\langle (\rho_0I-\widetilde{A})^{\frac{1}{2}}A(\rho_0I-\widetilde{A})^{-\frac{1}{2}}f,g\rangle_{X_{-\frac{1}{2}}}.
\end{eqnarray*}
    Then it follows from \eqref{yu-25-6-23-5} that $\widehat{\mathcal{A}}=\mathcal{A}$ when $\gamma\in[\frac{1}{2},1)$.
\end{proof}
\begin{lemma}\label{yu-lemma-3-28-1}
    The dual space  $X_\gamma^\flat$ of $X_{-\gamma}$ w.r.t. the pivot space $\mathcal{X}$ is given by
\begin{equation}\label{yu-3-28-4}
    X_\gamma^\flat=
\begin{cases}
    X_\gamma^*&\mbox{if}\;\;\gamma\in[0,\frac{1}{2}),\\
    (\rho_0I-\widetilde{A})^{\frac{1}{2}}X_{\gamma-\frac{1}{2}}^*&\mbox{if}\;\;\gamma\in[\frac{1}{2},1).
\end{cases}
\end{equation}
    Here, the norm of $(\rho_0I-\widetilde{A})^{\frac{1}{2}}X_{\gamma-\frac{1}{2}}^*$ is defined as follows:
    $\|f\|_{(\rho_0I-\widetilde{A})^{\frac{1}{2}}X_{\gamma-\frac{1}{2}}^*}=\|(\rho_0 I-A^*)^{\gamma-\frac{1}{2}}    (\rho_0I-\widetilde{A})^{-\frac{1}{2}}f\|_X$ for $f\in (\rho_0I-\widetilde{A})^{\frac{1}{2}}X_{\gamma-\frac{1}{2}}^*$.
\end{lemma}
\begin{proof}
    When $\gamma\in[0,\frac{1}{2})$, the result follows from the fact that $X_{-\gamma}'=X_\gamma^*$ immediately. Thus,
     it suffices to show the result for the case  $\gamma\in[\frac{1}{2},1)$. Indeed,
    by \eqref{yu-3-28-10}, we have that for any $f\in  (\rho_0I-\widetilde{A})^{\frac{1}{2}}X_{\gamma-\frac{1}{2}}^*
    \subset X_{-\frac{1}{2}}$ (i.e., $(\rho_0I-A^*)^{\gamma-\frac{1}{2}}(\rho_0I-\widetilde{A})^{-\frac{1}{2}}f\in X$) and $g\in X_{-\frac{1}{2}}$,
\begin{eqnarray}\label{yu-3-29-10}
    \langle f,g\rangle_{X_{-\frac{1}{2}}}&=&\langle (\rho_0 I-\widetilde{A})^{-\frac{1}{2}}f,
    (\rho_0 I-\widetilde{A})^{-\frac{1}{2}}g\rangle_X\nonumber\\
    &=&\langle (\rho_0I-A^*)^{\frac{1}{2}-\gamma}(\rho_0I-A^*)^{\gamma-\frac{1}{2}}
(\rho_0 I-\widetilde{A})^{-\frac{1}{2}}f,
    (\rho_0 I-\widetilde{A})^{-\frac{1}{2}}g\rangle_X\nonumber\\
    &=&\langle(\rho_0 I-\widetilde{A})^{\frac{1}{2}} (\rho_0I-A^*)^{\gamma-\frac{1}{2}}
(\rho_0 I-\widetilde{A})^{-\frac{1}{2}}f,
    (\rho_0 I-\widetilde{A})^{-\gamma+\frac{1}{2}}g\rangle_{X_{-\frac{1}{2}}}.
\end{eqnarray}
    This implies that
\begin{equation}\label{yu-3-29-11-b}
    \|g\|_{X_{-\gamma}}=\|(\rho_0I-\widetilde{A})^{-\gamma+\frac{1}{2}}g\|_{X_{-\frac{1}{2}}}=\sup_{f\in B_1^{*,\gamma}(0)}\langle f,g\rangle_{X_{-\frac{1}{2}}}
    \;\;\mbox{for any}\;\;g\in X_{-\frac{1}{2}},
\end{equation}
    where $B_1^{*,\gamma}(0):=\{f\in(\rho_0I-\widetilde{A})^{\frac{1}{2}}X_{\gamma-\frac{1}{2}}^*:\|f\|_{ (\rho_0I-\widetilde{A})^{\frac{1}{2}}X_{\gamma-\frac{1}{2}}^*}\leq 1\}$.
    Hence, according to \eqref{yu-3-29-11-b} and the same way used in the proof of \cite[Proposition 2.10.2, Chapter 2]{Tucsnak-Weiss}, $X_{-\gamma}$ is the dual space of
    $(\rho_0I-\widetilde{A})^{\frac{1}{2}}X_{\gamma-\frac{1}{2}}^*$ w.r.t. the pivot space $X_{-\frac{1}{2}}$.
\end{proof}
\begin{lemma}\label{yu-lemma-3-29-3}
    Let $\mathcal{S}^*(\cdot)$ be the dual semigroup of $\mathcal{S}(\cdot)$. Let $\mathcal{A}^*$ and $\mathcal{B}^*$ be the adjoint operators of $\mathcal{A}$ and $B$ w.r.t. the pivot space $\mathcal{X}$, respectively. Then
\begin{equation}\label{yu-3-29-12}
    \mathcal{S}^*(\cdot)=
\begin{cases}
    S^*(\cdot)&\mbox{if}\;\;\gamma\in[0,\frac{1}{2}),\\
(\rho_0I-\widetilde{A})^{\frac{1}{2}}S^*(\cdot)(\rho_0I-\widetilde{A})^{-\frac{1}{2}}
&\mbox{if}\;\;\gamma\in[\frac{1}{2},1),
\end{cases}
\end{equation}
\begin{equation}\label{yu-3-29-13}
    \mathcal{A}^*=
\begin{cases}
    A^*&\mbox{if}\;\;\gamma\in[0,\frac{1}{2}),\\
    (\rho_0I-\widetilde{A})^{\frac{1}{2}}A^*(\rho_0I-\widetilde{A})^{-\frac{1}{2}}&\mbox{if}\;\;\gamma\in[\frac{1}{2},1)
\end{cases}
\end{equation}
   with its domain
\begin{equation}\label{yu-3-29-bbb-1}
    D(\mathcal{A}^*)
    =\begin{cases}
        X_1^*&\mbox{if}\;\;\gamma\in[0,\frac{1}{2}),\\
        (\rho_0I-\widetilde{A})^{\frac{1}{2}}X_1^*&\mbox{if}\;\;\gamma\in[\frac{1}{2},1),
    \end{cases}
\end{equation}
    and
\begin{equation}\label{yu-3-29-14}
    \mathcal{B}^* =
\begin{cases}
    B^*&\mbox{if}\;\;\gamma\in[0,\frac{1}{2}),\\
    B^*(\rho_0I-A^*)^{-\frac{1}{2}}(\rho_0I-\widetilde{A})^{-\frac{1}{2}}&\mbox{if}\;\;\gamma\in[\frac{1}{2},1).
\end{cases}
\end{equation}
    (Here, it should be noted that  $A^*$ and $B^*$ are the adjoint operators of $A$ and $B$ w.r.t. the pivot space $X$, respectively.)
\end{lemma}
 \begin{proof}
    When $\gamma\in[0,\frac{1}{2})$, \eqref{yu-3-29-12}, \eqref{yu-3-29-13} and \eqref{yu-3-29-14} are obvious. Thus, it suffices to show
    them when $\gamma\in[\frac{1}{2},1)$.
\par
    We first prove \eqref{yu-3-29-12}. Since $\gamma\in [\frac{1}{2},1)$,  $\mathcal{X}=X_{-\frac{1}{2}}$. Then
      for any $f,g\in \mathcal{X}$ and $t\in\mathbb{R}^+$, we have that
\begin{eqnarray*}\label{yu-3-29-15}
    \langle \mathcal{S}(t)f,g\rangle_\mathcal{X}
    =\langle (\rho_0I-\widetilde{A})^{-\frac{1}{2}}f,S^*(t)(\rho_0I-\widetilde{A})^{-\frac{1}{2}}g\rangle_X=\langle f,(\rho_0I-\widetilde{A})^{\frac{1}{2}}S^*(t)(\rho_0I-\widetilde{A})^{-\frac{1}{2}}g\rangle_{X_{-\frac{1}{2}}}.
\end{eqnarray*}
    This implies \eqref{yu-3-29-12}. By the same way to show \eqref{yu-3-29-12}, we can obtain \eqref{yu-3-29-13}.
\par
    We next prove \eqref{yu-3-29-14}. For this purpose, we define
$\Phi_\gamma:=(\rho_0I-\widetilde{A})^{\frac{1}{2}}(\rho_0I-A^*)^{\gamma-\frac{1}{2}}
    (\rho_0I-\widetilde{A})^{-\frac{1}{2}}$.
    It is clear that $\Phi_\gamma$ is a linear  isomorphic map from $X_{\gamma}^\flat$ to $X_{-\frac{1}{2}}$, where $X_{\gamma}^\flat$ is given by \eqref{yu-3-28-4}.
    By the definition of dual space w.r.t. the pivot space $X_{-\frac{1}{2}}$ (or, by \eqref{yu-3-29-10} and the density of $X_{-\frac{1}{2}}$ in $X_{-\gamma}$), we have that
$\langle f,g\rangle_{X_{-\gamma},X_\gamma^\flat}=\langle (\rho_0I-\widetilde{A})^{-\gamma+\frac{1}{2}}f,\Phi_\gamma g\rangle_{X_{-\frac{1}{2}}}$ for any $f\in X_{-\gamma}$ and $g\in X_\gamma^\flat$.
    Here, $\langle \cdot,\cdot\rangle_{X_{-\gamma},X_\gamma^\flat}$ is the dual product of $X_{-\gamma}$ and $X_\gamma^\flat$ (see Lemma \ref{yu-lemma-3-28-1}). Therefore, for any $u\in U$ and $g\in X_{\gamma}^\flat$, we get that
\begin{eqnarray*}\label{yu-3-31-1}
    \langle Bu,g\rangle_{X_{-\gamma},X^\flat_\gamma}&=&\langle (\rho_0I-\widetilde{A})^{-\gamma+\frac{1}{2}}Bu,\Phi_\gamma g\rangle_{X_{-\frac{1}{2}}}=
    \langle (\rho_0I-\widetilde{A})^{-\gamma}Bu,(\rho_0I-A^*)^{\gamma-\frac{1}{2}}(\rho_0I-\widetilde{A})^{-\frac{1}{2}}g\rangle_X
    \nonumber\\
    &=& \langle Bu, (\rho_0I-A^*)^{-\frac{1}{2}}(\rho_0I-\widetilde{A})^{-\frac{1}{2}}g\rangle_{X_{-\gamma},X^*_\gamma}=\langle u,B^*(\rho_0I-A^*)^{-\frac{1}{2}}(\rho_0I-\widetilde{A})^{-\frac{1}{2}}g\rangle_U,
\end{eqnarray*}
    which indicates \eqref{yu-3-29-14}. The proof is completed.
    \end{proof}
\par
    We are now in a position to prove Proposition \ref{yu-proposition-6-27-1}.
\begin{proof}[Proof of Proposition \ref{yu-proposition-6-27-1}]
    $(i)$ follows easily from Lemma \ref{yu-lemma-3-29-1}. According to $(ii)$ and
     \cite[Theorem 4.4.3 in Section 4.4 and Proposition 4.2.5 in Section 4.2, Chapter 4]{Tucsnak-Weiss}, $(iii)$ is clear. Thus, it suffices to show $(ii)$.

To this end, let $\mathcal{X}^*_1:=D(\mathcal{A}^*)$ and $\mathcal{X}_{-1}:=D(\mathcal{A}^*)'$ (the dual space of $\mathcal{X}^*_1$ w.r.t. the pivot space $\mathcal{X}$, see the footnote in Definition \ref{yu-definition-10-18-1}). We firstly prove that $B\in \mathcal{L}(U;\mathcal{X}_{-1})$. By
    Lemma \ref{yu-lemma-3-29-3}, it suffices to show that
\begin{equation}\label{yu-3-31-4}
    \mathcal{B}^*\in\mathcal{L}(\mathcal{X}^*_1;U),
\end{equation}
    where $\mathcal{B}^*$ is given by \eqref{yu-3-29-14}. Indeed, when $\gamma\in[0,\frac{1}{2})$, by
    $(A_2)$, Remark \ref{yu-remark-8-00-1} and \eqref{yu-3-29-14}, we have that for any $f\in \mathcal{X}^*_1$,
\begin{eqnarray}\label{yu-3-31-3}
    \|\mathcal{B}^* f\|_U=\|B^*f\|_U\leq \|B^*(\rho_0I-A^*)^{-\gamma}\|_{\mathcal{L}(X;U)}
    \|(\rho_0I-A^*)^{\gamma-1}\|_{\mathcal{L}(X)}\|f\|_{\mathcal{X}^*_1}.
\end{eqnarray}
    When $\gamma\in[\frac{1}{2},1)$,  by
    $(A_2)$, Remark \ref{yu-remark-8-00-1} and \eqref{yu-3-29-14} again, we obtain that for any $f\in \mathcal{X}^*_1$,
\begin{eqnarray*}\label{yu-3-31-5}
    \|\mathcal{B}^* f\|_U&=&\|B^*(\rho_0 I-A^*)^{-\frac{1}{2}}(\rho_0I-\widetilde{A})^{-\frac{1}{2}}f\|_U\nonumber\\
    &\leq&\|B^*(\rho_0 I-A^*)^{-\gamma}\|_{\mathcal{L}(X;U)}\|(\rho_0 I-A^*)^{\gamma-\frac{1}{2}}
    (\rho_0I-\widetilde{A})^{-\frac{1}{2}}f\|_X\nonumber\\
    &=&\|B^*(\rho_0 I-A^*)^{-\gamma}\|_{\mathcal{L}(X;U)}\|
    (\rho_0I-A^*)^{\gamma-\frac{3}{2}}(\rho_0I-A^*)(\rho_0I-\widetilde{A})^{-\frac{1}{2}}f\|_{X}\nonumber\\
    &\leq&\|B^*(\rho_0 I-A^*)^{-\gamma}\|_{\mathcal{L}(X;U)}
    \|(\rho_0I-\mathcal{A}^*)^{\gamma-\frac{3}{2}}\|_{\mathcal{L}(\mathcal{X})}\|f\|_{\mathcal{X}^*_1}.
\end{eqnarray*}
    This, along with \eqref{yu-3-31-3}, implies \eqref{yu-3-31-4}.
\par
    We next prove that for each $T>0$, there  exists $C(T)>0$ so that
\begin{equation}\label{yu-25-7-10-1}
    \int_0^T\|\mathcal{B}^* \mathcal{S}^*(t)\varphi\|_U^2dt\leq C(T)\|\varphi\|_{\mathcal{X}}^2\;\;\mbox{for any}\;\;\varphi\in \mathcal{X}^*_1.
\end{equation}
    If it can be done, then according to the fact that $B\in\mathcal{L}(U;\mathcal{X}_{-1})$ and Lemma \ref{yu-lemma-3-29-3}, $B$ is admissible
    (w.r.t. $\mathcal{A}$) on $\mathcal{X}$. Let $T>0$ be fixed arbitrarily. By \cite[Theorem 6.13, Chapter 2]{Pazy},
    $(A_1)$, $(A_2)$ and Lemma \ref{yu-lemma-3-29-3}, we have that for each $\gamma\in[0,\frac{1}{2})$, there exists $C(\gamma)>0$ so that
\begin{eqnarray*}\label{yu-3-31-6}
    \int_0^T\|\mathcal{B}^* \mathcal{S}^*(t)\varphi\|_U^2dt=\int_0^T\|B^*S^*(t)\varphi\|_U^2dt
    \leq C(\gamma)T^{1-2\gamma}\|B^*(\rho_0I-A^*)^{-\gamma}\|^2_{\mathcal{L}(\mathcal{X};U)}
    \|\varphi\|_{\mathcal{X}}^2
\end{eqnarray*}
   for any $\varphi\in \mathcal{X}$; for each $\gamma\in[\frac{1}{2},1)$, there exists $C(\gamma)>0$ so that
\begin{eqnarray*}\label{yu-3-31-7}
    \int_0^T\|\mathcal{B}^* \mathcal{S}^*(t)\varphi\|_U^2dt
    &\leq&\|B^*(\rho_0I-A^*)^{-\gamma}\|^2_{\mathcal{L}(X;U)}\int_0^T\|(\rho_0I-A^*)^{\gamma-\frac{1}{2}}
    S^*(t)(\rho_0I-\widetilde{A})^{-\frac{1}{2}}\varphi\|_X^2dt\nonumber\\
    &\leq&C(\gamma)T^{2(1-\gamma)}\|B^*(\rho_0I-A^*)^{-\gamma}\|^2_{\mathcal{L}(X;U)}
    \|\varphi\|_{\mathcal{X}}^2\;\;\mbox{for any}\;\;\varphi\in \mathcal{X}.
\end{eqnarray*}
    Therefore, (\ref{yu-25-7-10-1}) holds.
    The proof is completed.
\end{proof}
\section{Stabilizability by bounded feedback laws}\label{yu-sec-25-3}
    In this section, we prove that under Assumption $(A_3)$ and the assumption that $B$ is admissible,
     the feedback law to stabilize the pair $[A,B]$ can be chosen in the family of bounded linear operators.
     For this purpose, we make the following assumptions:
\begin{enumerate}
     \item[$(A_1')$] The linear operator $A:D(A)(:=X_1)\subset X\to X$ generates  a $C_0$-semigroup $S(\cdot)$ on $X$.
\item[$(A_2')$] The operator $B$  is admissible w.r.t. $A$ (i.e., $B\in\mathcal{L}(U;X_{-1})$ and for each $T>0$, there exists  $C(T)>0$ so that \eqref{yu-25-6-22-1} holds).
\end{enumerate}

\begin{lemma}\label{yu-prop-5-22-1}
     Under Assumptions $(A'_1)$ and $(A_2')$, the following statements are true:
\begin{enumerate}
  \item [$(i)$] Suppose that $(A_3)$ holds. The pair $[A,B]$ is stabilizable in $X$ if and only if
   there exist $\alpha>0$ and  $C(\alpha)>0$ so that
  \begin{equation}\label{yu-5-22-1}
    \|\varphi\|_{X}^{2}\leq C(\alpha)(\|(\lambda I-A^*)\varphi\|^2_{X}+\|B^*\varphi\|^2_{U})
    \;\;\mbox{for any}\;\;\lambda\in \mathbb{C}_{-\alpha}^+,\;\varphi\in X^*_1.
  \end{equation}
  \item [$(ii)$]  Suppose that the operator $A$ satisfies one of the assumptions $(a)$ and $(b)$ in $(iv)$ of Remark \ref{yu-remark-12-22-1}. Then the pair $[A,B]$ is rapidly stabilizable in $X$ if and only if
   for each $\alpha>0$, there  exists  $C(\alpha)>0$
     so that \eqref{yu-5-22-1} holds.
\end{enumerate}
\end{lemma}
\begin{proof}
     According to $(ii)$ of Remark \ref{yu-remark-12-22-1}, $(i)$ follows from \cite[Theorem 2]{Kunisch-Wang-Yu}. It suffices to show $(ii)$.

     Note that $A$ satisfies the condition $(a)$/$(b)$ in $(iv)$ of Remark \ref{yu-remark-12-22-1} if and only if $A+\alpha I$ satisfies the condition $(a)$/$(b)$ for any $\alpha>0$. By $(i)$ and $(iv)$ of Remark \ref{yu-remark-12-22-1}, we only need to claim that \emph{the pair $[A,B]$ is rapidly stabilizable in $X$ if and only if for any $\alpha>0$, the pair $[A+\alpha I,B]$
    is stabilizable} in $X$. Indeed, from  \cite[Theorem 3.4]{Liu-Wang-Xu-Yu} it follows that, the pair $[A,B]$ is rapidly stabilizable in $X$ if and only if for each $\alpha>0$, there exist $C_0(\alpha)>0$ and $D_0(\alpha)>0$ so that
\begin{equation}\label{yu-12-11-1}
    \|S^*(t)\varphi\|_X^2\leq C_0(\alpha)\int_0^t\|B^*S^*(s)\varphi\|_U^2ds+D_0(\alpha) e^{-\alpha t}\|\varphi\|_X^2\;\;\mbox{for any}\;\;t\in\mathbb{R}^+,\;\varphi\in X_1^*.
\end{equation}
    On one hand, if $[A,B]$ is rapidly stabilizable in $X$, then for each $\alpha>0$, it follows from \eqref{yu-12-11-1} with $\alpha$ replaced by $3\alpha$ that
\begin{equation}\label{yu-12-11-2}
    \|S^*_\alpha(t)\varphi\|_X^2\leq C_0(3\alpha)e^{2\alpha t}\int_0^t\|B^*S^*_\alpha(s)\varphi\|_U^2ds+D_0(3\alpha) e^{-\alpha t}\|\varphi\|_X^2\;\;\mbox{for any}\;\;t\in\mathbb{R}^+,\;\varphi\in X^*_1.
\end{equation}
    Here and in what follows, $S_\alpha^*(\cdot)$ is the $C_0$-semigroup on $X$ generated by $A^*+\alpha I$. Take $T>0$ so that $\delta_T:=D_0(3\alpha)e^{-\alpha T}<1$. Thus, by \eqref{yu-12-11-2}, we have
\begin{equation}\label{yu-12-11-3}
    \|S^*_\alpha(T)\varphi\|_X^2\leq C_0(3\alpha)e^{2\alpha T}\int_0^T\|B^*S^*_\alpha(s)\varphi\|_U^2ds+\delta_T\|\varphi\|_X^2\;\;\mbox{for any}\;\;\varphi\in X^*_1.
\end{equation}
    This, along with \cite[Theorem 26]{Trelat-Wang-Xu} and \cite[Proposition 3.9]{Liu-Wang-Xu-Yu}, implies that $[A+\alpha I,B]$ is stabilizable in $X$.
    On the other hand, if the pair $[A+\alpha I,B]$ is stabilizable in $X$ for each $\alpha>0$, then by \cite[Theorem 26]{Trelat-Wang-Xu} and \cite[Proposition 3.9]{Liu-Wang-Xu-Yu}, we  conclude that there exist $T>0$, $\delta\in(0,1)$ and $C_1(\delta,T)>0$ so that \eqref{yu-12-11-3} holds by replacing $C_0(3\alpha)e^{2\alpha T}$ and $\delta_T$ by $C_1(\delta,T)$ and $\delta$, respectively. This, together with \cite[Proposition 3]{Ma-Wang-Yu}, yields that there exist $\varepsilon>0$, $C(\alpha,\varepsilon)>0$ and $D(\alpha,\varepsilon)>0$ so that
\begin{equation*}\label{yu-12-11-4}
\|S^*_\alpha(t)\varphi\|_X^2\leq C(\alpha,\varepsilon)\int_0^t\|B^*S^*_\alpha(s)\varphi\|_U^2ds+D(\alpha,\varepsilon) e^{-\varepsilon t}\|\varphi\|_X^2\;\;\mbox{for any}\;\;t\in\mathbb{R}^+,\;\varphi\in X_1^*.
\end{equation*}
     It follows from the latter that \eqref{yu-12-11-1} holds with $D_0(\alpha)=C(\alpha,\varepsilon)$ and $D_0(\alpha)=D(\alpha,\varepsilon)$. Since $\alpha$ is arbitrary, we conclude that the pair $[A,B]$ is rapidly stabilizable in $X$.
     The proof is completed.
\end{proof}
    The main result of this section is stated as follows.

\begin{theorem}\label{yu-theorem-10-18-1}
   Under Assumptions $(A_1')$, $(A_2')$ and $(A_3)$, the following statements are equivalent:
\begin{enumerate}
  \item [$(i)$] The pair $[A,B]$ is  stabilizable in $X$ in the sense of Definition \ref{yu-definition-10-18-1}
  (by replacing $\mathcal{A}$ by $A$).
  \item [$(ii)$]  The pair $[A,B]$  is stabilizable in $X$ with bounded feedback law.
\end{enumerate}
\end{theorem}
\begin{remark}
    Two remarks on Theorem \ref{yu-theorem-10-18-1} are provided in order.
\begin{enumerate}
  \item [$(i)$]
    According to Assumptions $(A_1')$ and $(A_2')$, the stabilizability property for the pair $[A,B]$ in the sense of Definition \ref{yu-definition-10-18-1} (by replacing $\mathcal{A}$ by $A$) is well defined.
  \item [$(ii)$] Theorem \ref{yu-theorem-10-18-1} means that under Assumptions $(A_1')$, $(A_2')$ and $(A_3)$, the feedback operator can be chosen in the family of bounded linear operators.
       However, since
      the $C_0$-semigroup $S(\cdot)$ generated by $A$ is not necessarily compact, the techniques in
       \cite{Badra-Takahashi, Fattorini-1966, Fattorini-1967, Raymond} are no longer applicable.
\end{enumerate}
\end{remark}
\begin{proof}[Proof of Theorem \ref{yu-theorem-10-18-1}]
Since  $(ii)\Rightarrow (i)$ is obvious, we only need to show $(i)\Rightarrow (ii)$.

Suppose that
     $(i)$ is true. It follows from $(i)$ in Lemma \ref{yu-prop-5-22-1} that there exist $\alpha>0$ and $C(\alpha)>0$ so that
\begin{equation}\label{yu-11-11-1}
    \|\varphi\|^2_X\leq C(\alpha)
    \left(\|(\lambda I-A^*)\varphi\|^2_X+\|B^*\varphi\|_U^2\right)
    \;\;\mbox{for any}\;\;\lambda\in\mathbb{C}_{-\alpha}^+,\;\varphi\in X_1^*.
\end{equation}
 Let $\beta\in(0,\alpha)$ be fixed arbitrarily. According to  Assumption $(A_3)$ and $(ii)$ in Remark \ref{yu-remark-12-22-1}, there are two closed subspaces $Q_1:=Q_1(\beta)$ and $Q_2:=Q_2(\beta)$ of $X$ so that $(a_1)$ $X=Q_1\oplus Q_2$;  $(a_2)$ $Q_1$ and $Q_2$ are  invariant subspaces of $S^*(\cdot)$;
    $(a_3)$ $A^*|_{Q_1}$ (the restriction of $A^*$ on $Q_1$) is bounded and satisfies that  $\sigma(A^*|_{Q_1})\subset\mathbb{C}_{-\beta}^+$; $(a_4)$  $S^*(\cdot)|_{Q_2}$
     (the restriction of $S^*(\cdot)$ on $Q_2$) is  exponentially stable.
  From $(a_1)$ it follows that for each $x\in X$, there  exists a unique $(x_1,x_2)\in Q_1\times Q_2$ so that $x=x_1+x_2$. Define a linear operator $P$ from $X$ to $Q_1$ as follows:
\begin{equation}\label{yu-12-12-1}
    Px:=x_1\;\;\mbox{if}\;\;x=x_1+x_2\;\;\mbox{with}\;\;x_1\in Q_1,\;x_2\in Q_2.
\end{equation}
    One can easily check that $P^2=P$,  $Q_1=PX$ and $Q_2=(I-P)X$.
   Since $Q_1$ and $Q_2$ are closed spaces, we have that $P\in \mathcal{L}(X)${\footnote{By the closed graph theorem, it suffices to show that $P$ is a closed operator. Indeed, if $x_n\to x$ and $Px_n\to f$ as $n\to+\infty$, then by the closedness of $Q_1$, we have that $f\in Q_1$. This, along with the closedness of
    $Q_2$, implies that $x-f=\lim_{n\to+\infty}(x_n-Px_n)\in Q_2$. Thus, $f-Px=P(f-x)=0$. It follows that $P$ is a closed operator.},
     $Q_1$ and $Q_2$ are two Hilbert spaces whose norms are inherited from $X$. Then $P$ is a projection operator.}
    By $(a_2)$, we claim that
\begin{equation}\label{yu-12-12-2}
    S^*(\cdot)P=PS^*(\cdot).
\end{equation}
    Indeed, by $(a_2)$, we have that for any $x=x_1+x_2\in X$ with $(x_1,x_2)\in Q_1\times Q_2$,
\begin{equation*}\label{yu-12-12-3}
    S^*(t)Px=S^*(t)x_1=PS^*(t)x_1=PS^*(t)(x_1+x_2)=PS^*(t)x\;\;\mbox{for each}\;\;t\in\mathbb{R}^+.
\end{equation*}
    Hence, \eqref{yu-12-12-2} is true. Let $S_1^*(\cdot):=PS^*(\cdot)$ and $S^*_2(\cdot):=(I-P)S^*(\cdot)$. According to \eqref{yu-12-12-2}, it is clear that $S^*_1(\cdot)$ and $S^*_2(\cdot)$ are two $C_0$-semigroups on $Q_1$ and $Q_2$, respectively. Suppose that the generators of $S_1^*(\cdot)$ (on $Q_1$) and $S^*_2(\cdot)$ (on $Q_2$) are $A^*_1$ (with domain $D(A_1^*)$) and $A_2^*$ (with domain $D(A_2^*)$), respectively. It follows from $(a_4)$ that there exist $\varepsilon>0$ and $C(\beta,\varepsilon)>0$ so that
\begin{equation}\label{yu-12-12-4}
    \|S^*_2(t)\|_{\mathcal{L}(X)}\leq C(\beta,\varepsilon)e^{-\varepsilon t}\;\;\mbox{for each}\;\;t\in\mathbb{R}^+.
\end{equation}
\par
    Let
$H_1:=P^*X$ and $H_2:=\mbox{Ker}(P^*)$.
    Since $P$ is a projection operator, one can directly check that $P^*$ is also a projection operator. Thus, $H_2=(I-P^*)X$ and $X=H_1\oplus H_2$.
    By \eqref{yu-12-12-2}, we have that
\begin{equation}\label{yu-12-12-6}
S(\cdot)P^*=P^*S(\cdot).
\end{equation}
    Further, we claim that
\begin{equation}\label{yu-12-12-6-01}
    H_1=Q_1'\;\;\mbox{and}\;\;H_2=Q_2'
\end{equation}
    in the following sense: for any linear bounded functionals $F(\cdot)$ and $G(\cdot)$ on $Q_1$ and $Q_2$, respectively, there are unique $f\in H_1$ and $g\in H_2$ so that $F(x)=\langle f,x\rangle_X$ for any $x\in Q_1$ and $G(x)=\langle g,x\rangle_X$ for any $x\in Q_2$. Indeed, if $F(\cdot)$ is a linear  bounded functional on $Q_1$, then $F\circ P$ is a linear bounded functional on $X$. Hence, there exists a unique $f\in X$ so that $F(Px)=\langle f, x\rangle_X$ for any $x\in X$,
     which, combined with the fact that $F(Px)=F(P(Px))$ for any $x\in X$, indicates that $\langle f,x\rangle_X=\langle f,Px\rangle_X=\langle P^*f,x\rangle_X$ for any $x\in X$. This implies that $f\in H_1$ and then $H_1$ is the dual space of $Q_1$. By the same way, we can show that $H_2=Q_2'$.

    Define $S_1(\cdot):=P^*S(\cdot)$ and $S_2(\cdot):=(I-P^*)S(\cdot)$.  It is obvious that $S_1(\cdot)$ and $S_2(\cdot)$ are two $C_0$-semigroups on $H_1$ and $H_2$, respectively. Moreover, by \eqref{yu-12-12-6-01} and \eqref{yu-12-12-6}, we get that $S_1(\cdot)$ and $S_2(\cdot)$ are the dual semigroups of $S_1^*(\cdot)$ and $S^*_2(\cdot)$, respectively. Thus,
the generators of $S_1(\cdot)$ (on $H_1$) and $S_2(\cdot)$ (on $H_2$)  are given by
\begin{equation}\label{yu-12-10-11111}
    A_1:=(A_1^*)^*\;\;\mbox{and}\;\;A_2:=(A_2^*)^*
\end{equation}
    in the following sense:
\begin{equation}\label{yu-12-10-1111-2}
\begin{cases}
    \langle A_1f_1,g_1\rangle_{H_1,Q_1}=\langle f_1,A^*_1g_1\rangle_{H_1,Q_1}\;\;\mbox{for any}\;\;f_1\in D(A_1),\;g_1\in D(A_1^*),\\
     \langle A_2f_2,g_2\rangle_{H_2,Q_2}=\langle f_2,A^*_2g_2\rangle_{H_2,Q_2}
     \;\;\mbox{for any}\;\;f_2\in D(A_2),\;g_2\in D(A_2^*).
\end{cases}
\end{equation}
    (Here, for $i=1,2$, $\langle \cdot,\cdot\rangle_{H_i,Q_i}$ is the dual product of $H_i$ and $Q_i$.) Furthermore, by \eqref{yu-12-12-6-01} and the fact that $Q_1$ is a Hilbert space, there exists a linear isomorphic mapping $\Phi: Q_1\to H_1$ so that
\begin{equation}\label{yu-12-13-10}
    \langle f,g\rangle_{H_1,Q_1}=\langle f,\Phi g\rangle_{H_1}\footnote{Here, it should be noted that
    $\langle f,g\rangle_{H_1,Q_1}=\langle f, g\rangle_X$ and $\langle f,\Phi g\rangle_{H_1}= \langle f,\Phi g\rangle_{X}$ for any $f\in H_1$ and $g\in Q_1$. It follows that $(I-\Phi)Q_1\subset
    H_1^\bot$. Moreover, it is obvious that
    $\|\Phi g\|_{H_1}=\|g\|_{Q_1}$  for any $g\in Q_1$
    and
$\Phi^*=\Phi^{-1}$.}
\;\;\mbox{for any}\;\;
    f\in H_1,\;g\in Q_1.
\end{equation}
\par
    The rest of the proof is split  into six steps.  In \emph{Steps 1} and \emph{2}, we establish some decomposition properties
    of $A^*$ and $A$ based on $X=Q_1\oplus Q_2$ and $X=H_1\oplus H_2$, respectively.  In \emph{Step 3}, we introduce the projection of the control
operator and show that it is a bounded control operator (while the initial control operator $B$ may be unbounded). In \emph{Step 4}, we exponentially stabilize the first part
of the system with a bounded feedback law. The last two steps are then devoted to proving that this bounded feedback law stabilizes the whole system.
It is far from
being obvious because the system is triangular in some sense and the feedback law
could destabilize the second part of the system. More precisely, in \emph{Step 5}, we prove that under the feedback law established in \emph{Step 4}, the state of the corresponding closed-loop system  is actually in $L^2(\mathbb{R}^+;X)$. Finally, in \emph{Step 6}, we show the exponential stability of the closed-loop system in the framework of Definition \ref{yu-definition-10-18-1}.
\vskip 5pt
\emph{Step 1. We claim that
\begin{equation}\label{yu-12-12-7}
  Q_1= PX^*_1\subset X_1^*,\;\;A^*P=PA^*\;\;\mbox{in}\;\;X^*_1
\end{equation}
and
\begin{equation}\label{yu-11-11-3}
\begin{cases}
    A^*_1=PA^*\;\;\mbox{in}\;\;Q_1\;\;\mbox{with}\;\;D(A^*_1)=PX^*_1=Q_1,\\
    A^*_2=(I-P)A^*\;\;\mbox{in}\;\;Q_2\;\mbox{with}\;\;D(A^*_2)=(I-P)X^*_1.
\end{cases}
\end{equation}
}

    Firstly, we prove \eqref{yu-12-12-7}. Indeed, by $(a_3)$ and  $(iii)$ in Remark \ref{yu-remark-12-22-1}, we have that $Q_1\subset X^*_1$. Then $Q_1=PX^*_1\subset X_1^*$.  Moreover,  it follows from \eqref{yu-12-12-2} and the boundness of $P$ that for any  $g\in X^*_1$,
\begin{equation*}\label{yu-12-12-9}
    A^*Pg=\lim_{t\to0^+}t^{-1}(S^*(t)Pg-Pg)=P\lim_{t\to0^+}t^{-1}(S^*(t)g-g)=PA^*g\in X.
\end{equation*}
    Thus,  \eqref{yu-12-12-7} holds.

    Secondly, we show the first conclusion in  \eqref{yu-11-11-3}. On one hand,
   for each $f\in PX^*_1$, i.e., $f=Pg$ for some $g\in X^*_1$, we have
\begin{equation}\label{yu-11-11-4}
    \lim_{t\to 0^+}t^{-1}(S_1^*(t)f-f)=P\lim_{t\to 0^+}
    t^{-1}(S^*(t)g-g)=PA^*g\in X.
\end{equation}
    This implies that $f\in D(A^*_1)$ and then $PX_1^*\subset D(A^*_1)$.
    On the other hand, if $f\in D(A^*_1)$, then  $f\in Q_1$ (i.e., $f=Pf$) and
\begin{equation}\label{yu-11-11-5}
    \lim_{t\to 0^+}t^{-1}(S^*(t)Pf-Pf)
    =\lim_{t\to0^+}t^{-1}(PS^*(t)f-f)=\lim_{t\to0^+}t^{-1}(S^*_1(t)f-f)=A^*_1f\in X.
\end{equation}
    It follows from the latter that $f=Pf\in X^*_1$ and then $D(A_1^*)\subset PX^*_1$. Therefore, by the first conclusion in  \eqref{yu-12-12-7} and  \eqref{yu-11-11-4}, we have $D(A^*_1)=PX^*_1=Q_1$
and  $A^*_1f=PA^*f$  for each $f\in Q_1$.
\par
    Finally, by similar arguments used to show \eqref{yu-11-11-4} and \eqref{yu-11-11-5}, we can prove the second conclusion in \eqref{yu-11-11-3}.

\vskip 5pt
\emph{Step 2. We claim that
\begin{equation}\label{yu-12-12-7-bbb}
   H_1=P^*X_1\subset X_1,\;\;AP^*=P^*A\;\;\mbox{in}\;\;X_1
\end{equation}
    and
\begin{equation}\label{yu-11-19-3}
\begin{cases}
    A_1=P^*A\;\;\mbox{in}\;\;H_1\;\;\mbox{with}\;\;D(A_1)=P^*X_1=H_1,\\
    A_2=(I-P^*)A\;\;\mbox{in}\;\;H_2\;\;\mbox{with}\;\;D(A_2)=(I-P^*)X_1.
\end{cases}
\end{equation}
}

   Since $P^*$ is also a projection operator, by \eqref{yu-12-10-11111} and the same way to show \eqref{yu-12-12-7} and
   \eqref{yu-11-11-3}, we can obtain \eqref{yu-12-12-7-bbb} and \eqref{yu-11-19-3}.

\vskip 5pt
    \emph{Step 3.
    We show that $B^*P\in \mathcal{L}(Q_1;U)$ and
\begin{equation}\label{yu-11-15-bb-1}
    (B^*P)^*=(\rho_0I-A_1)P^*(\rho_0I-\widetilde{A})^{-1}B\in \mathcal{L}(U;H_1)
\end{equation}
    in the sense: $\langle u, B^*P\varphi\rangle_U=\langle (\rho_0I-A_1)P^*(\rho_0I-\widetilde{A})^{-1}Bu, \varphi\rangle_{H_1,Q_1}$ for any $\varphi\in Q_1$ and $u\in U$.}

     By \eqref{yu-12-12-7},  \eqref{yu-11-11-3} and $(iii)$ in Remark \ref{yu-remark-12-22-1}, we have that
\begin{equation}\label{yu-1-20-1}
    A_1^*=A^*|_{Q_1}\in\mathcal{L}(Q_1).
\end{equation}
    Then it follows from \eqref{yu-12-12-7}, \eqref{yu-11-11-3} and Assumption $(A_2')$ that for each $\varphi\in Q_1$,
\begin{eqnarray*}\label{yu-11-14-1}
    \|B^*P\varphi\|_U&=&\|B^*(\rho_0I-A^*)^{-1}(\rho_0I-A^*)P\varphi\|_U
    \leq \|B^*(\rho_0I-A^*)^{-1}\|_{\mathcal{L}(X;U)}(\rho_0+\|A_1^*\|_{\mathcal{L}(Q_1)})
    \|\varphi\|_{Q_1}.
\end{eqnarray*}
    This implies that $B^*P\in \mathcal{L}(Q_1;U)$.
    Moreover,  according to Assumption $(A_2')$, \eqref{yu-12-10-11111}, \eqref{yu-12-12-7} and \eqref{yu-11-19-3}, it is clear that for any $\varphi\in Q_1$ and  $u\in U$,
\begin{eqnarray*}\label{yu-11-15-b-2}
    \langle u, B^*P\varphi\rangle_U&=&\langle Bu,P\varphi\rangle_{X_{-1},X^*_{1}}
    =\langle(\rho_0I-\widetilde{A})^{-1}Bu, P(\rho_0I-A^*)\varphi\rangle_X\nonumber\\
    &=& \langle P^*(\rho_0I-\widetilde{A})^{-1}Bu, (\rho_0I-A^*_1)\varphi\rangle_{H_1,Q_1}
    =\langle(\rho_0I-A_1)P^*(\rho_0I-\widetilde{A})^{-1}Bu,\varphi\rangle_{H_1,Q_1},
\end{eqnarray*}
    which indicates  \eqref{yu-11-15-bb-1}.

\vskip 5pt
\emph{Step 4.  We prove that the pair $[A_1,(\rho_0I-A_1)P^*
    (\rho_0I-\widetilde{A})^{-1}B]$ is exactly controllable in $H_1$ and then, for each $\mu>0$, there exists $F:=F(\mu)\in\mathcal{L}(H_1;U)$ so that any solution $z_F(\cdot)$ of the following equation (in $H_1$):
\begin{equation}\label{yu-11-19-8}
    z_t(t)=\left[A_1+(\rho_0I-A_1)P^*
    (\rho_0I-\widetilde{A})^{-1}BF\right]z(t),\;\;t\in\mathbb{R}^+
\end{equation}
    satisfies
\begin{equation}\label{yu-11-19-9}
    \|z_F(t)\|_{H_1}\leq C(\mu)e^{-\mu t}\|z_F(0)\|_{H_1}\;\;\mbox{for any}\;\;t\in\mathbb{R}^+,
\end{equation}
    where $C(\mu)>0$ is a constant.}

    For this purpose, we first claim that there is a constant $C>0$ so that
\begin{equation}\label{yu-11-15-1}
    \|\varphi\|^2_{Q_1}\leq C\left(\|(\lambda I-A^*_1)\varphi\|^2_{Q_1}+\|B^*P\varphi\|_U^2\right)
    \;\;\mbox{for any}\;\;\lambda\in\mathbb{C},\;\;\varphi\in Q_1.
\end{equation}
    Indeed, by \eqref{yu-11-11-1}, \eqref{yu-12-12-7} and \eqref{yu-11-11-3}, we obtain that
\begin{equation}\label{yu-11-13-4}
    \|\varphi\|^2_{Q_1}\leq C(\alpha)
    \left(\|(\lambda I-A^*_1)\varphi\|^2_{Q_1}+\|B^*P\varphi\|_U^2\right)
    \;\;\mbox{for any}\;\;\lambda\in\mathbb{C}_{-\alpha}^+,\;\;\varphi\in Q_1.
\end{equation}
    It follows from \eqref{yu-1-20-1} that $-A^*_1$ with domain $Q_1$ generates an uniformly continuous $C_0$-group on $Q_1$ (which is given by $(S^*_1(\cdot))^{-1}$).
    By $(a_3)$ above, we have that $\sigma(-A^*_1)\subset
    \overline{\mathbb{C}_{\beta}^{-}}$. Thus, it follows from  the fact that $\beta\in(0,\alpha)$, \cite[Corollary 3.12, Section 3, Chapter 4]{Engel-Nagel} and \cite[Theorem 5.3 and Remark 5.4, Section 1.5, Chapter 1]{Pazy} that there exists $C(\beta,\alpha)>0$ so that
$\|(\lambda I+A_1^*)^{-1}\varphi\|_{Q_1}\leq
    C(\alpha,\beta)\|\varphi\|_{Q_1}$ for any $\lambda\in \mathbb{C}_{\frac{\alpha+\beta}{2}}^+$
     and $\varphi\in Q_1$,
    which indicates
$\|\varphi\|_{Q_1}\leq
    C(\alpha,\beta)\|(\lambda I-A_1^*)\varphi\|_{Q_1}$  for any $\lambda\in \mathbb{C}_{-\frac{\alpha+\beta}{2}}^-$ and
    $\varphi\in Q_1$.
    This, along with \eqref{yu-11-13-4}, implies  \eqref{yu-11-15-1}.

    By \eqref{yu-11-15-1}, we have
\begin{equation}\label{yu-12-13-13}
     \|\varphi\|^2_{H_1}\leq C\left(\|(\lambda I-\Phi A^*_1\Phi^{-1})\varphi\|^2_{H_1}+\|B^*P\Phi^{-1}\varphi\|_U^2\right)
    \;\;\mbox{for any}\;\;\lambda\in\mathbb{C},\;\;\varphi\in H_1.
\end{equation}
    (Here, we recall that the linear isomorphic mapping $\Phi: Q_1\to H_1$ satisfies \eqref{yu-12-13-10}.) Let $(\Phi A^*_1\Phi^{-1})^*$ and $(B^*P\Phi^{-1})^*$ be the adjoint operators of $\Phi A^*_1\Phi^{-1}$ and $B^*P\Phi^{-1}$ in the sense:
    $\langle (\Phi A_1^*\Phi^{-1})^*f, g\rangle_{H_1}=\langle f, \Phi A_1^*\Phi^{-1}g\rangle_{H_1}$ for any $f,g\in H_1$
    and
    $\langle (B^*P\Phi^{-1})^*u,\varphi\rangle_{H_1}=
    \langle u, B^*P\Phi^{-1}\varphi\rangle_U$ for any
    $u\in U$ and $\varphi\in H_1$.
On one hand, by \eqref{yu-12-13-10} and \eqref{yu-12-10-1111-2}, we get
\begin{equation*}\label{yu-12-16-1}
    \langle (\Phi A_1^*\Phi^{-1})^*f, g\rangle_{H_1}
    =\langle f,\Phi A_1^*\Phi^{-1}g\rangle_{H_1}=\langle f,A_1^*\Phi^{-1}g\rangle_{H_1, Q_1}
    =\langle A_1f,\Phi^{-1}g\rangle_{H_1, Q_1}=\langle A_1f, g\rangle_{H_1}
\end{equation*}
    for any $f,g\in H_1$. On the other hand, by \eqref{yu-11-15-bb-1}, we have
\begin{eqnarray*}\label{yu-12-16-2}
    \langle (B^*P\Phi^{-1})^*u,\varphi\rangle_{H_1}
    =\langle (\rho_0 I-A_1)P^*(\rho_0I-\widetilde{A})^{-1}Bu,\Phi^{-1}\varphi\rangle_{H_1,Q_1}
    =\langle (\rho_0 I-A_1)P^*(\rho_0I-\widetilde{A})^{-1}Bu,\varphi\rangle_{H_1}
\end{eqnarray*}
    for any $u\in U$ and $\varphi\in H_1$. Therefore, $[(\Phi A^*_1\Phi^{-1})^*,(B^*P\Phi^{-1})^*]
    =[A_1,(\rho_0 I-A_1)P^*(\rho_0I-\widetilde{A})^{-1}B]$, which, combined with \eqref{yu-12-13-13}
    and Lemma \ref{yu-lemma-12-9-1}, indicates that the pair
    $[A_1,(\rho_0 I-A_1)P^*(\rho_0I-\widetilde{A})^{-1}B]$ is exactly controllable in $H_1$. Hence, by Lemma \ref{yu-lemma-12-10-120}, we  conclude that for each $\mu>0$, there exist $F:=F(\mu)\in\mathcal{L}(H_1;U)$ and $C(\mu)>0$ so that any solution $z_F(\cdot)$ of the equation \eqref{yu-11-19-8} in $H_1$  satisfies \eqref{yu-11-19-9}.

\vskip 5pt
    \emph{Step 5. Fix $\mu>0$ arbitrarily. Let $F\in \mathcal{L}(H_1;U)$ be given so that the solution of the equation \eqref{yu-11-19-8} satisfies \eqref{yu-11-19-9}. Define
\begin{equation}\label{yu-11-19-10}
    K=K(\mu):=FP^*.
\end{equation}
    We claim that for each $y_0\in X$, the following integral equation:
\begin{equation}\label{yu-11-19-11}
    y(t)=S(t)y_0+\int_0^t\widetilde{S}(t-s)BKy(s)ds,\;\;t\in\mathbb{R}^+,
\end{equation}
    has a unique solution $y_K(\cdot;y_0)$ in $C([0,+\infty);X)\cap L^2(\mathbb{R}^+;X)$.}

Let $y_0\in X$ be fixed arbitrarily. In what follows, we choose $T>0$ to be such that
\begin{equation}\label{yu-bbb-11-21-b-1}
    \sqrt{2}C(\beta,\varepsilon)e^{-\varepsilon T}<1,
\end{equation}
    where $C(\beta,\varepsilon)>0$ is given so that \eqref{yu-12-12-4} holds.
    Note that by \cite[Corollary 5.5.1, Section 5.5, Chapter 5]{Tucsnak-Weiss}, the equation
    \eqref{yu-11-19-11} has a unique solution $y_K(\cdot;y_0)$ in $C([0,+\infty);X)$.
    Thus, it suffices to show that
\begin{equation}\label{yu-12-16-10}
    \int_{\mathbb{R}^+}\|y_K(t;y_0)\|_X^2dt<+\infty.
\end{equation}
    For this purpose, on one hand, let $z_1(\cdot):=P^*y_K(\cdot;y_0)$
    and $z_2(\cdot):=(I-P^*)y_K(\cdot;y_0)$. By $(iii)$ in Remark \ref{yu-remark-6-22-1}, \eqref{yu-12-12-6}, \eqref{yu-11-19-10} and \eqref{yu-11-19-11}, we have that for each $\varphi\in H_1$ and any $t\in\mathbb{R}^+$,
\begin{eqnarray*}\label{yu-11-20-6}
    &\;&\langle z_1(t),\varphi\rangle_{H_1}=\langle z_1(t),\Phi^{-1}\varphi\rangle_{H_1, Q_1}\nonumber\\
    &=&\langle z_1(t),\Phi^{-1}\varphi\rangle_X
    =\langle S_1(t)P^*y_0,\Phi^{-1}\varphi\rangle_{H_1,Q_1}+
    \left\langle\int_0^t\widetilde{S}(t-s)BKy_K(s)ds,
    \Phi^{-1}\varphi\right\rangle_{X_{-1},X^*_1}\nonumber\\
    &=&\langle S_1(t)P^*y_0,\Phi^{-1}\varphi\rangle_{H_1,Q_1}+\int_0^t\langle
S_1(t-s)P^*(\rho_0I-\widetilde{A})^{-1}BFP^*
y_K(s),(\rho_0I-A^*_1)\Phi^{-1}\varphi\rangle_{H_1,Q_1}ds\nonumber\\
&=&\left\langle S_1(t)P^*y_0+\int_0^t
S_1(t-s)(\rho_0I-A_1)P^*(\rho_0I-\widetilde{A})^{-1}
BFz_1(s)ds,\varphi\right\rangle_{H_1}.
\end{eqnarray*}
    This yields that $z_1(\cdot)$ is the solution of the equation
    \eqref{yu-11-19-8} with the initial data $P^*y_0$. Thus, by \eqref{yu-11-19-9}, we get
\begin{equation}\label{yu-11-20-7}
    \|z_1(t)\|_X=\|z_1(t)\|_{H_1}\leq  C(\mu)e^{-\mu t}\|P^*y_0\|_{H_1}
    =C(\mu)e^{-\mu t}\|P^*y_0\|_X\;\;\mbox{for each}\;\;t\in\mathbb{R}^+.
\end{equation}
     According to \eqref{yu-11-19-10} and \eqref{yu-11-19-11}, it is clear that
\begin{equation}\label{yu-11-21-2}
    z_2(t)=S_2(t)(I-P^*)y_0+(I-P^*)\int_0^t\widetilde{S}(t-s)
    BFz_1(s)ds\;\;\mbox{for each}\;\;t\in\mathbb{R}^+.
\end{equation}
   Moreover, by $(iii)$ in Remark \ref{yu-remark-6-22-1} again, we  conclude that for any $\varphi\in X^*_1$,
\begin{eqnarray*}\label{yu-11-21-1}
    \langle z_2(t),\varphi\rangle_X
    =\langle S_2(t)(I-P^*)y_0,\varphi\rangle_X
    +\left\langle \int_0^tS(t-s)(\rho_0I-\widetilde{A})^{-1}
    BFz_1(s)ds,(\rho_0I-A^*)(I-P)\varphi\right\rangle_{X}.
\end{eqnarray*}
  Here, we note that $(I-P)X^*_1\subset X^*_1$ by \eqref{yu-12-12-7}. Thus,
    for each $n\in\mathbb{N}^+$ and $t\in[nT,(n+1)T]$, we have
\begin{eqnarray*}\label{yu-11-21-3}
    \langle z_2(t),\varphi\rangle_X
    &=&\langle S_2(nT)(I-P^*)y_0,S_2^*(t-nT)\varphi\rangle_X\nonumber\\
    &\;&+\left\langle \int_0^{t}S(t-s)(\rho_0I-\widetilde{A})^{-1}
    BFz_1(s)ds,(\rho_0I-A^*)(I-P)\varphi\right\rangle_{X}\nonumber\\
    &=&\langle S_2(nT)(I-P^*)y_0,S_2^*(t-nT)\varphi\rangle_X\nonumber\\
    &\;&+\left\langle \int_0^{t-nT}S(t-nT-s)(\rho_0I-\widetilde{A})^{-1}
    BFz_1(nT+s)ds, (\rho_0I-A^*)(I-P)\varphi\right\rangle_{X}\nonumber\\
    &\;&+\left\langle \int_0^{nT}S(nT-s)(\rho_0I-\widetilde{A})^{-1}
    BFz_1(s)ds, (I-P)(\rho_0I-A^*)S^*_2(t-nT)\varphi\right\rangle_{X}\nonumber\\
    &=&\langle S_2(t-nT)z_2(nT),\varphi\rangle_X+\left\langle (I-P^*)\int_0^{t-nT}\widetilde{S}(t-nT-s)
    BFz_1(nT+s)ds,\varphi\right\rangle_{X}.
\end{eqnarray*}
    It follows from the density of $X^*_1$ in $X$ that for each $n\in\mathbb{N}^+$,
\begin{equation}\label{yu-11-21-4}
    z_2(t)=S_2(t-nT)z_2(nT)
    +(I-P^*)\int_0^{t-nT}\widetilde{S}(t-nT-s)
    BFz_1(nT+s)ds\;\;\mbox{for any}\;\;t\in[nT,(n+1)T].
\end{equation}
    On the other hand, by Assumption $(A_2')$ and $(i)$ in Remark \ref{yu-remark-6-23-2}, we claim that there exists $C(T)>0$ so that for any $u\in L^2(0,T;U)$,
\begin{equation}\label{yu-11-20-1}
    \left\|\int_0^t\widetilde{S}(t-s)Bu(s)ds\right\|^2_X\leq
    C(T)\int_0^t\|u(s)\|_U^2ds\;\;\mbox{for each}\;\;t\in[0,T].
\end{equation}
   Indeed, for any $t\in[0,T]$, $u\in L^2(0,T;U)$ and $\varphi\in X^*_1$, it follows from $(i)$ in Remark \ref{yu-remark-6-23-2} that
\begin{equation*}\label{yu-11-26-13}
    \left\langle \int_0^t\widetilde{S}(t-s)Bu(s)ds,\varphi\right\rangle_X
    =\int_0^t\langle \widetilde{S}(t-s)Bu(s),\varphi\rangle_{X_{-1},X^*_1}ds
    =\int_0^t\langle u(s),B^*S^*(t-s)\varphi\rangle_Uds.
\end{equation*}
    This, along with Assumption $(A_2')$ and $(i)$ in Remark \ref{yu-remark-6-23-2}, implies that there exists
    $C(T)>0$  (independent of $u$) so that
\begin{eqnarray*}\label{yu-11-26-14}
    \left|\left\langle \int_0^t\widetilde{S}(t-s)Bu(s)ds,\varphi\right\rangle_X\right|
    &\leq&\left(\int_0^t\|u(s)\|_U^2ds\right)^{\frac{1}{2}}
    \left(\int_0^T\|B^*S^*(\sigma)\varphi\|_U^2d\sigma\right)^{\frac{1}{2}}\nonumber\\
    &\leq& (C(T))^{\frac{1}{2}}\left(\int_0^t\|u(s)\|_U^2ds\right)^{\frac{1}{2}}
\|\varphi\|_X,
\end{eqnarray*}
   which, combined with the density of $X^*_1$ in $X$, indicates  \eqref{yu-11-20-1}.
  Then by \eqref{yu-12-12-4}, \eqref{yu-11-21-2}, \eqref{yu-11-21-4} and
    \eqref{yu-11-20-1}, we obtain that
\begin{equation}\label{yu-11-21-5}
    \sup_{t\in[0,T]}\|z_2(t)\|_X^2\leq C(\beta,\varepsilon,T,\mu)
    \|y_0\|_X^2
\end{equation}
    and for each $n\in\mathbb{N}^+$,
\begin{equation}\label{yu-11-21-6}
    \begin{cases}
    \sup_{t\in[nT,(n+1)T]}\|z_2(t)\|_X^2\leq C(\beta,\varepsilon,T,\mu)
    \left(\|z_2(nT)\|_X^2+\int_{nT}^{(n+1)T}\|z_1(s)\|^2_Xds\right),\\
    \|z_2((n+1)T)\|_X^2\leq   2(C(\beta,\varepsilon))^2e^{-2\varepsilon T}\|z_2(nT)\|_X^2+C(T,\mu)\int_{nT}^{(n+1)T}\|z_1(s)\|^2_Xds.
\end{cases}
\end{equation}
     From the second inequality in \eqref{yu-11-21-6} it follows that for any $N\in\mathbb{N}^+$ with $N>2$,
\begin{equation*}\label{yu-11-21-7}
     \sum_{n=2}^{N+1}\|z_2(nT)\|_X^2\leq 2(C(\beta,\varepsilon))^2e^{-2\varepsilon T}
     \sum_{n=1}^N\|z_2(nT)\|_X^2+C(T,\mu)\int_{T}^{(N+1)T}\|z_1(s)\|^2_Xds.
\end{equation*}
     This, together with \eqref{yu-bbb-11-21-b-1} and \eqref{yu-11-21-5}, yields that for each $N\in\mathbb{N}^+$ with $N>2$,
\begin{eqnarray*}\label{yu-11-21-8}
    \sum_{n=1}^N\|z_2(nT)\|_X^2\leq (1-2(C(\beta,\varepsilon))^2e^{-2\varepsilon T})^{-1}\left(C(\beta,\varepsilon,T,\mu)
    \|y_0\|_X^2+C(T,\mu)\int_{0}^{(N+1)T}\|z_1(s)\|^2_Xds\right),
\end{eqnarray*}
     which, combined with \eqref{yu-11-20-7}, indicates
\begin{eqnarray}\label{yu-11-21-9}
    \sum_{n=1}^{+\infty}\|z_2(nT)\|_X^2<+\infty.
\end{eqnarray}
    Hence, by \eqref{yu-11-20-7}, \eqref{yu-11-21-5}, the first inequality in \eqref{yu-11-21-6}
    and \eqref{yu-11-21-9}, we obtain that
\begin{eqnarray*}\label{yu-11-21-10}
    \int_{\mathbb{R}^+}\|z_2(t)\|_X^2dt&=&
    \sum_{n=0}^{+\infty}\int_{nT}^{(n+1)T}\|z_2(t)\|_X^2dt
    \leq T\left(\sup_{t\in[0,T]}\|z_2(t)\|_X^2
    +\sum_{n=1}^{+\infty}\sup_{t\in[nT,(n+1)T]}\|z_2(t)\|_X^2\right)\nonumber\\
    &\leq&TC(\beta,\varepsilon,T,\mu)\left(\|y_0\|_X^2+\sum_{n=1}^{+\infty}\|z_2(nT)\|_X^2
    +\int_{\mathbb{R}^+}\|z_1(s)\|_X^2ds\right)<+\infty.
\end{eqnarray*}
    This, along with \eqref{yu-11-20-7}, implies that
 \eqref{yu-12-16-10} holds.
\vskip 5pt
    \emph{Step 6. We complete the proof.}

    Let $K$ be defined by \eqref{yu-11-19-10} with $\mu>0$. It is obvious that $K\in\mathcal{L}(X;U)$. Since the equation \eqref{yu-11-19-11} is linear, there exists an operator-valued mapping $\Upsilon_K(\cdot):\mathbb{R}^+\to\mathcal{L}(X)$ so that for any $y_0\in X$,
\begin{equation}\label{yu-11-22-1}
    y_K(t;y_0):=\Upsilon_K(t)y_0\;\;\mbox{for each}\;\;t\in\mathbb{R}^+,
\end{equation}
    where $y_K(\cdot;y_0)$ is the solution of the equation \eqref{yu-11-19-11}. On one hand, according to the continuity of $y_K(\cdot;y_0)$, it is clear that $\lim_{t\to0^+}\Upsilon_K(t)y_0=y_0$ in $X$. Moreover, by similar arguments as those
     used for \eqref{yu-11-21-4}, we have that
$y_K(t+s;y_0)=y_K(t;y_K(s;y_0))$ for any $t,s\in\mathbb{R}^+$.
     This implies that $\Upsilon_K(t+s)=\Upsilon_K(t)\Upsilon_K(s)$ for any $t,s\in\mathbb{R}^+$. Therefore, $\Upsilon_K(\cdot)$ is a $C_0$-semigroup on $X$.  Suppose that the generator of $\Upsilon_K(\cdot)$ is $A_K$ with domain $D(A_K)$.  Let $y_0\in D(A_K)$ and $\varphi\in D((A^*)^2)$ be fixed arbitrarily.
     It follows from \eqref{yu-11-19-11} and \eqref{yu-11-22-1} that
\begin{equation*}\label{yu-11-22-3}
    \langle \Upsilon_K(t)y_0,\varphi\rangle_X
    =\langle y_0,S^*(t)\varphi\rangle_X+
    \int_0^t\langle K\Upsilon_K(s)y_0,B^*S^*(t-s)\varphi\rangle_{U}ds
    \;\;\mbox{for all}\;\;t\in\mathbb{R}^+,
\end{equation*}
    which indicates
\begin{eqnarray}\label{yu-11-22-4}
    &\;&\langle \Upsilon_K(t)A_Ky_0,\varphi\rangle_X=\frac{d}{dt}\langle \Upsilon_K(t)y_0,\varphi\rangle_X\nonumber\\&=&\langle S(t)y_0,A^*\varphi\rangle_{X}
    +\langle  K\Upsilon_K(t)y_0,B^*\varphi\rangle_{U}+\int_0^t\langle K\Upsilon_K(s)y_0,B^*S^*(t-s)A^*\varphi\rangle_{U}ds.
\end{eqnarray}
    Passing to the limit for $t\to0^+$ in \eqref{yu-11-22-4}, we get
\begin{equation*}\label{yu-11-22-5}
    \langle A_Ky_0,\varphi\rangle_{X_{-1},X^*_1}=\langle A_Ky_0,\varphi\rangle_X=\langle\widetilde{A}y_0+BKy_0,\varphi\rangle_{X_{-1},X^*_1}\;\;\mbox{for each}\;\;
    y_0\in D(A_K)\;\;\mbox{and}\;\;\varphi\in D((A^*)^2).
\end{equation*}
     Since $y_0$ and $\varphi$ are arbitrary, it follows from the density of $D((A^*)^2)$ in $X_1^*$ that
\begin{equation}\label{yu-11-22-b-5}
    A_Kx=\widetilde{A}x+BKx\;\;\mbox{for all}\;\;x\in D(A_K).
\end{equation}
    On the other hand, by \emph{Step 5}, we have
\begin{equation}\label{yu-11-22-6}
    \int_{\mathbb{R}^+}\|\Upsilon_K(t)y_0\|_X^2dt<+\infty\;\;\mbox{for any}\;\;y_0\in X.
\end{equation}
    Since $\Upsilon_K(\cdot)$ is a $C_0$-semigroup on $X$, it follows from \eqref{yu-11-22-6} and \cite[Lemma 5.1, Section 5.1.2, Chapter 5]{Curtain-Zwart}  that there exist $\gamma>0$ and $C(\gamma)>0$ so that
\begin{equation}\label{yu-11-22-7}
    \|\Upsilon_K(t)\|_{\mathcal{L}(X)}\leq C(\gamma)e^{-\gamma t}\;\;\mbox{for all}\;\;t\in \mathbb{R}^+.
\end{equation}
    Hence, by \eqref{yu-11-22-b-5} and \eqref{yu-11-22-7}, we conclude that the pair $[A,B]$ is stabilizable in $X$ with bounded feedback law.
    In summary, we finish the proof of Theorem~\ref{yu-theorem-10-18-1}.
\end{proof}
\begin{remark}\label{yu-remark-4-14-1}
    Several observations  on the proof of Theorem \ref{yu-theorem-10-18-1} are listed as follows:
\begin{enumerate}
  \item [$(i)$]  The key assumption in the proof of Theorem \ref{yu-theorem-10-18-1} is $(A_3)$. Assumption $(A_3)$ implies that the control system can be split into two parts. Both of them are invariant under the action of the operator $A$. Moreover, the unbounded (admissible) control
operator $B$ is projected onto the first part of the system, which yields a bounded control operator (see \emph{Step 4} in the proof of Theorem \ref{yu-theorem-10-18-1}). Thanks to this fact,
we can design a bounded feedback law. More precisely, the bounded feedback is designed by only using
the first part of the system, i.e., the feedback only depends on the first component of the control system.
Then this feedback is put into the full control system. The most difficult part of the proof consists of showing that the second component (which is impacted by the feedback) remains
exponentially stable.
Similar strategies have been developed in \cite{Badra-Takahashi, Coron-Trelat-2004, Coron-Trelat-2006, Lhachemi-Prieur-Trelat, Russell}, but in these references the first part of the system is finite-dimensional. It should be noted that  we do not have any
Lyapunov functional here. To prove the stability of the full system, we use the exponential stability of the uncontrolled second component and we have to show that the introduction
of the feedback (designed from the first part, but which appears in the dynamics of the second part) does not destabilize the system. This is one of the difficulties in the proof.
\item [$(ii)$]
Since we do not have any Lyapunov functional in general, we do not know whether our bounded feedback control strategy is robust or not to various perturbations, for example, perturbations of
the operators $A$ and $B$, additional semilinearities or noise. This opens interesting new questions regarding the treatment of semilinear systems.
  \item [$(iii)$] In the proof of Theorem \ref{yu-theorem-10-18-1}, we provide an algorithmic procedure to construct the feedback (see \eqref{yu-11-19-11}). In fact, the algorithmic procedure we provided only relies on the structure of the first part of the system mentioned in $(i)$.
      Specifically, the algorithmic procedure  we provided is as follows. Let $P$ be given by \eqref{yu-12-12-1}. According to \emph{Step 4} in the proof of Theorem \ref{yu-theorem-10-18-1}, the pair $[A_1,(\rho_0I-A_1)P^*
    (\rho_0I-\widetilde{A})^{-1}B]$ is exactly controllable in $H_1$ (the expression of $A_1$ can be found in \eqref{yu-11-19-3}). Therefore, for each $\mu>0$, we can use Lemma \ref{yu-lemma-12-10-120} to construct a feedback law $F:=F(\mu)\in \mathcal{L}(H_1;U)$ so that $A_1+(\rho_0I-A_1)P^*
    (\rho_0I-\widetilde{A})^{-1}B F$ is exponentially stable with decay rate $\mu$. Let $K:=K(\mu)$ be defined by \eqref{yu-11-19-10}. By \emph{Step 6} in the proof of Theorem \ref{yu-theorem-10-18-1}, we  conclude that $K$ is the needed feedback law. The algorithmic procedure to construct the feedback in the proof of Theorem \ref{yu-theorem-10-18-1} is briefly described as follows:
\begin{eqnarray*}
    &\;&[A,B]\xrightarrow[Steps\;1-3]{P}
  [P^*A,(\rho_0I-P^*A)P^*(\rho_0I-\widetilde{A})^{-1}B]\;\mbox{in}\;P^*X\;\nonumber\\
  &\;&\xrightarrow[Step\; 4]{F\in \mathcal{L}(P^*X;U)} z_t=(P^*A+(\rho_0I-P^*A)P^*(\rho_0I-\widetilde{A})^{-1}BF)z
  \;\mbox{is stable in}\;P^*X\nonumber\\
  &\;&\xrightarrow[Steps\; 5-6]{K:=FP^*\in \mathcal{L}(X;U)}
  K\;\mbox{is a feedback to stabilize}\;[A,B].
\end{eqnarray*}
\end{enumerate}
\end{remark}
\par
    As an application of Theorem \ref{yu-theorem-10-18-1}, we have the following corollary.
\begin{corollary}\label{corollary-yu-12-4-1}
    Suppose that Assumptions $(A'_1)$ and $(A_2')$  hold. If the operator $A$ satisfies one of the conditions $(a)$ and $(b)$ in
    $(iv)$ of Remark \ref{yu-remark-12-22-1} (by replacing $L$ by $A$),
then the following  statements are true:
\begin{enumerate}
  \item [$(i)$] The pair $[A,B]$ is stabilizable in $X$ if and only if the pair $[A,B]$ is stabilizable in $X$ with bounded feedback law.
  \item [$(ii)$] The pair $[A,B]$ is rapidly stabilizable in $X$ if and only if the pair $[A,B]$ is rapidly stabilizable in $X$ with bounded feedback law.
\end{enumerate}
\end{corollary}
\begin{proof}
    According to Lemma \ref{yu-lemma-12-20-1-bb}, it is clear that if $A$ satisfies one of the assumptions $(a)$ and $(b)$ in $(iv)$ of Remark \ref{yu-remark-12-22-1}, then Assumption $(A_3)$ is true. Thus, $(i)$ follows from Theorem \ref{yu-theorem-10-18-1} directly.

    We next show $(ii)$.
    On one hand, if the pair $[A,B]$ is rapidly stabilizable in $X$ with bounded feedback law, then it is also rapidly stabilizable in $X$. On the other hand, if the pair $[A,B]$ is rapidly stabilizable in $X$, then
     the pair $[A+\alpha I,B]$ is stabilizable in $X$ for each $\alpha>0$ (see \cite{Liu-Wang-Xu-Yu, Trelat-Wang-Xu} or the proof of Lemma \ref{yu-prop-5-22-1}). This, together with $(i)$, yields that the pair $[A+\alpha I,B]$ is stabilizable in $X$ with bounded feedback law for each $\alpha>0$. Hence, for each $\alpha>0$, there exist $\varepsilon:=\varepsilon(\alpha)>0$ and $K:=K(\alpha)\in \mathcal{L}(X;U)$ so that  for each $y_0\in X$, the solution $y_{K,\alpha}(\cdot;y_0)$ of the following equation:
\begin{equation*}\label{yu-12-20-1}
\begin{cases}
    y_t(t)=(A+\alpha I)y(t)+BKy(t),\;\;t\in\mathbb{R}^+,\\
    y(0)=y_0
\end{cases}
\end{equation*}
    satisfies
\begin{equation}\label{yu-12-20-2}
    \|y_{K,\alpha}(t;y_0)\|_X\leq C(K,\alpha,\varepsilon)e^{-\varepsilon t}\|y_0\|_X\;\;\mbox{for any}\;\;t\in\mathbb{R}^+.
\end{equation}
    It is clear that the solution $z(\cdot;y_0)$ of the following equation:
\begin{equation}\label{yu-25-7-22-1}
\begin{cases}
    z_t(t)=Az(t)+BKz(t),\;\;t\in\mathbb{R}^+,\\
    z(0)=y_0
\end{cases}
\end{equation}
    satisfies that  $e^{\alpha t}z(t;y_0)=y_{K,\alpha}(t;y_0)$ for each $t\in\mathbb{R}^+$ (by \cite[Corollary 5.5.1, Section 5.5, Chapter 5]{Tucsnak-Weiss}, the equation \eqref{yu-25-7-22-1} is well-posed).
    This, along with \eqref{yu-12-20-2}, implies that
$\|z(t;y_0)\|_X\leq C(\alpha)e^{-\alpha t}\|y_0\|_X$ for any $t\in\mathbb{R}^+$.
    Since $\alpha$ is arbitrary, it follows that the pair $[A,B]$ is rapidly stabilizable in $X$ with bounded feedback law. The proof is completed.
\end{proof}

\section{Proofs of Theorem \ref{yu-theorem-3-31-1} and Corollary \ref{yu-corollary-4-22-1}}\label{yu-sect-5-12-1}

\begin{proof}[Proof of Theorem \ref{yu-theorem-3-31-1}]
    By Assumption $(A_3)$ and $(ii)$ in Remark \ref{yu-remark-12-22-1}, we observe that
    $A^*$ satisfies the property (AEDC), i.e.,
    (\textbf{CL}): For each $\alpha>0$, there  exist  two closed subspaces $Q_1:=Q_1(\alpha)$ and $Q_2:=Q_2(\alpha)$ of $X$ so that $(a)$ $X=Q_1\oplus Q_2$;  $(b)$ $Q_1$ and $Q_2$ are  invariant subspaces of $S^*(\cdot)$;
    $(c)$ $A^*|_{Q_1}$  is bounded and satisfies that  $\sigma(A^*|_{Q_1})\subset\mathbb{C}_{-\alpha}^+$; $(d)$  $S^*(\cdot)|_{Q_2}$ is  exponentially stable.

     The proof consists of the following five steps.
\vskip 5pt
    \emph{Step 1. Fix $\alpha>0$ arbitrarily. Let $Q_1:=Q_1(\alpha)$ and $Q_2:=Q_2(\alpha)$ be the closed subspaces of $X$ so that  $(a)$-$(d)$ in (\textbf{CL}) hold. Define
\begin{equation}\label{yu-4-2-6}
    W_1:=
\begin{cases}
     Q_1&\mbox{if}\;\;\gamma\in[0,\frac{1}{2}),\\
    (\rho_0I-\widetilde{A})^{\frac{1}{2}}Q_1&\mbox{if}\;\;\gamma\in[\frac{1}{2},1),
\end{cases}
    \;\;\mbox{and}
    \;\;W_2:=\begin{cases}
     Q_2&\mbox{if}\;\;\gamma\in[0,\frac{1}{2}),\\
    (\rho_0I-\widetilde{A})^{\frac{1}{2}}Q_2&\mbox{if}\;\;\gamma\in[\frac{1}{2},1).
\end{cases}
\end{equation}
     We claim that $W_1$ and $W_2$ are two closed subspaces of $\mathcal{X}$, and satisfy that  $(e_1)$ $\mathcal{X}=W_1\oplus W_2$;  $(e_2)$ $W_1$ and $W_2$ are  invariant subspaces of $\mathcal{S}^*(\cdot)$;
    $(e_3)$ $\mathcal{A}^*|_{W_1}$ is bounded and $\sigma(\mathcal{A}^*|_{W_1})\subset\mathbb{C}_{-\alpha}^+$; $(e_4)$  $\mathcal{S}^*(\cdot)|_{W_2}$ is  exponentially stable.}
\par
    Its proof is divided into two cases.
\par
    \emph{Case 1. $\gamma\in[0,\frac{1}{2})$.} By \eqref{yu-25-6-23-5-b}, we have $\mathcal{X}=X$. Thus, the result
     follows from (\textbf{CL}) directly.
\par
    \emph{Case 2. $\gamma\in[\frac{1}{2},1)$.} According to \eqref{yu-4-2-6} and the boundness of $(\rho_0I-\widetilde{A})^{\frac{1}{2}}$ (from $X$ to $\mathcal{X}$), it is clear that $W_1$ and $W_2$ are two closed subspaces of $\mathcal{X}$.

    Firstly, we show $(e_1)$. It is clear that
\begin{equation}\label{yu-4-2-7}
    W_1+W_2=W_1\oplus W_2.
\end{equation}
    Indeed, if $f_1+f_2=0$ for $f_1\in W_1$ and $f_2\in W_2$, then by \eqref{yu-4-2-6}, we have $(\rho_0I-\widetilde{A})^{-\frac{1}{2}}f_1\in Q_1$, $(\rho_0I-\widetilde{A})^{-\frac{1}{2}}f_2\in Q_2$ and $(\rho_0I-\widetilde{A})^{-\frac{1}{2}}f_1+(\rho_0I-\widetilde{A})^{-\frac{1}{2}}f_2=0$. These,
     along with the fact that $Q_1+Q_2=Q_1\oplus Q_2$, imply $f_1=f_2=0$. Thus, \eqref{yu-4-2-7} holds.
    We now prove that
\begin{equation}\label{yu-4-2-8}
    \mathcal{X}\subset W_1+W_2.
\end{equation}
    To this end, for each $f\in \mathcal{X}$, we have $(\rho_0I-\widetilde{A})^{-\frac{1}{2}}f\in X$
    (by \eqref{yu-25-6-23-5-b}). According to $(a)$ in (\textbf{CL}), there exists $(g_1,g_2)\in Q_1\times Q_2$ so that
    $g_1+g_2=(\rho_0I-\widetilde{A})^{-\frac{1}{2}}f$, which indicates $f=(\rho_0I-\widetilde{A})^{\frac{1}{2}}(g_1+g_2)$. This, together with \eqref{yu-4-2-6}, yields that
    $f\in W_1+ W_2$. Since $f$ is arbitrary, \eqref{yu-4-2-8} is true. By \eqref{yu-4-2-7} and \eqref{yu-4-2-8}, we obtain $(e_1)$.

\par
   Secondly, we prove $(e_2)$. Since $Q_1$ is an invariant subspace of $S^*(\cdot)$ (see $(b)$ in (\textbf{CL})),
   it follows from Lemma \ref{yu-lemma-3-29-3} and \eqref{yu-4-2-6} that
$\mathcal{S}^*(\cdot)W_1\subset (\rho_0I-\widetilde{A})^{\frac{1}{2}}S^*(\cdot)Q_1\subset (\rho_0I-\widetilde{A})^{\frac{1}{2}}Q_1=W_1$.
    This implies that $W_1$ is an invariant subspace of $\mathcal{S}^*(\cdot)$. Similarly, we can obtain that
    $W_2$ is also an invariant subspace of $\mathcal{S}^*(\cdot)$.

\par
    Thirdly, we show $(e_3)$. Take $f\in W_1$ arbitrarily. By \eqref{yu-4-2-6}, there  exists $h\in Q_1$ so that $f=(\rho_0I-\widetilde{A})^{\frac{1}{2}}h$.
    This, along with \eqref{yu-3-29-13} and $(c)$ in (\textbf{CL}), implies that
\begin{eqnarray*}\label{yu-4-3-2}
    \|\mathcal{A}^* f\|_{\mathcal{X}}=\|(\rho_0I-\widetilde{A})^{\frac{1}{2}}A^*h\|_{X_{-\frac{1}{2}}}
    =\|A^*h\|_X\leq \|A^*|_{Q_1}\|_{\mathcal{L}(X)}\|h\|_X
    \leq \|A^*|_{Q_1}\|_{\mathcal{L}(X)}\|f\|_{\mathcal{X}}.
\end{eqnarray*}
    Since $f$ is arbitrary, it follows from the latter that $\mathcal{A}^*|_{W_1}$ is bounded. Furthermore, according to
    \eqref{yu-3-29-13}, it is clear that $\sigma(\mathcal{A}^*|_{W_1})=\sigma(A^*|_{Q_1})$, which, combined with
    $(c)$ in (\textbf{CL}), indicates that  $\sigma(\mathcal{A}^*|_{W_1})\subset\mathbb{C}_{-\alpha}^+$.

\par
    Finally, we prove $(e_4)$. Take $f\in W_2$ arbitrarily. By \eqref{yu-4-2-6}, there  exists $g\in Q_2$ so that
    $f=(\rho_0I-\widetilde{A})^{\frac{1}{2}}g$. This, together with \eqref{yu-3-29-12} and $(d)$ in (\textbf{CL}), yields that there  exist $\delta>0$ and $C(\delta)>0$ (which is independent of $f$) so that
$\|\mathcal{S}^*(t)f\|_{\mathcal{X}}=\|(\rho_0I-\widetilde{A})^{\frac{1}{2}}S^*(t)g\|_{X_{-\frac{1}{2}}}
    =\|S^*(t)g\|_X\leq C(\delta)e^{-\delta t}\|g\|_X=C(\delta)e^{-\delta t}\|f\|_{\mathcal{X}}$
    for each $t\in\mathbb{R}^+$. Since $f$ is arbitrary, it follows from the above estimate that $(e_4)$ holds.

\vskip 5pt
    \emph{Step 2. We claim that $(i)\Leftrightarrow (ii)$.}

    According to Proposition \ref{yu-proposition-6-27-1}, \emph{Step 1}, $(ii)$ in Remark \ref{yu-remark-12-22-1} and Theorem \ref{yu-theorem-10-18-1}, the equivalence $(i)\Leftrightarrow (ii)$ is obvious.

\vskip 5pt
    \emph{Step 3. We show that $(i)\Rightarrow (iii)$.}

    Its proof is divided into two cases.
\par
    \emph{Case 1.  $\gamma\in[0,\frac{1}{2})$.} In this case, $\mathcal{X}=X$ (see \eqref{yu-25-6-23-5-b}), $\mathcal{A}^*= A^*$ and $\mathcal{B}^*= B^*$ (see \eqref{yu-3-29-13} and \eqref{yu-3-29-14}). Hence, $(iii)$ follows from
    the conclusion $(i)\Leftrightarrow (ii)$ (see \emph{Step 2}) immediately.

\par
    \emph{Case 2. $\gamma\in[\frac{1}{2},1)$.}  By Proposition \ref{yu-proposition-6-27-1} and \emph{Step 1}, we  apply \cite[Theorem 2]{Kunisch-Wang-Yu} to the pair $[\mathcal{A},B]$ and obtain that the pair $[\mathcal{A},B]$ is stabilizable in $\mathcal{X}$ if and only if there exist $\vartheta>0$ and $C(\vartheta)>0$ so that
\begin{equation}\label{yu-4-7-2}
    \|\varphi\|_{\mathcal{X}}^2\leq C(\vartheta)(\|(\lambda I-\mathcal{A}^*)\varphi\|_{\mathcal{X}}^2+
    \|\mathcal{B}^*\varphi\|_U^2)\;\;\mbox{for any}\;\;\varphi\in \mathcal{X}^*_1,\;\;\lambda\in\mathbb{C}_{-\vartheta}^+,
\end{equation}
    where $\mathcal{A}^*$ and $\mathcal{B}^*$ are the adjoint operators of $\mathcal{A}$ and $B$ w.r.t.
    the pivot space $\mathcal{X}$, respectively (see Lemma \ref{yu-lemma-3-29-3}). Thus, it suffices to show $(iii)$ from \eqref{yu-4-7-2}. To this end, we firstly claim that
\begin{equation}\label{yu-4-14-1}
    \|\varphi\|_X^2\leq C(\vartheta)(\|(\lambda I-A^*)\varphi\|_X^2+\|B^*(\rho_0I-A^*)^{-\frac{1}{2}}\varphi\|_U^2)
    \;\;\mbox{for any}\;\;\varphi\in X_1^*,\;\;\lambda\in\mathbb{C}_{-\vartheta}^+.
\end{equation}
    Indeed, by \eqref{yu-3-29-13}, \eqref{yu-3-29-bbb-1} and \eqref{yu-3-29-14}, we have that for each $\varphi\in X_1^*$,
\begin{equation*}\label{yu-4-14-2}
\begin{cases}
    (\rho_0I-\widetilde{A})^{\frac{1}{2}}\varphi\in \mathcal{X}^*_1,\\
    (\lambda I-A^*)\varphi=(\rho_0I-\widetilde{A})^{-\frac{1}{2}}
    (\lambda I-\mathcal{A}^*)(\rho_0I-\widetilde{A})^{\frac{1}{2}}\varphi
    \;\;\mbox{for any}\;\;\lambda\in\mathbb{C},\\
    B^*(\rho_0I-A^*)^{-\frac{1}{2}}\varphi=\mathcal{B}^* (\rho_0I-\widetilde{A})^{\frac{1}{2}}\varphi.
\end{cases}
\end{equation*}
    These, along with $(ii)$ in Remark \ref{yu-remark-6-22-1}, \eqref{yu-25-6-23-5-b} and \eqref{yu-4-7-2}, imply that
\begin{eqnarray*}\label{yu-4-14-3}
    \|\varphi\|_X^2&=&\|(\rho_0I-\widetilde{A})^{\frac{1}{2}}\varphi\|_{\mathcal{X}}^2
    \leq C(\vartheta)(\|(\lambda I-\mathcal{A}^*)(\rho_0I-\widetilde{A})^{\frac{1}{2}}\varphi\|_{\mathcal{X}}^2
    +\|\mathcal{B}^*(\rho_0I-\widetilde{A})^{\frac{1}{2}}\varphi\|_U^2)\nonumber\\
    &=& C(\vartheta)(\|(\lambda I-A^*)\varphi\|_X^2+\|B^*(\rho_0I-A^*)^{-\frac{1}{2}}\varphi\|_U^2)\;\;
    \mbox{for any}\;\;\varphi\in X^*_1\;\;\mbox{and}\;\;\lambda\in\mathbb{C}_{-\vartheta}^+.
\end{eqnarray*}
   Hence, \eqref{yu-4-14-1} holds.

\par
     Fix $\alpha\in(0,\vartheta)$ arbitrarily. Let $Q_1:=Q_1(\alpha)$ and $Q_2:=Q_2(\alpha)$ be the closed subspaces of $X$ so that $(a)$-$(d)$ in \textbf{(CL)} hold.
      By $(a)$ in \textbf{(CL)} and the same arguments as those in the proof of Theorem \ref{yu-theorem-10-18-1}, we  define a projection operator $P:X\to Q_1$ as follows:
\begin{equation}\label{yu-4-14-4}
    Pf=f_1\;\;\mbox{if}\;\;f=f_1+f_2\;\;\mbox{with}\;\;(f_1,f_2)\in Q_1\times Q_2.
\end{equation}
    Let $S_1^*(\cdot):=PS^*(\cdot)$ and $S^*_2(\cdot):=(I-P)S^*(\cdot)$. It follows from $(b)$ in \textbf{(CL)} that $S_1^*(\cdot)$ and $S_2^*(\cdot)$
      are $C_0$-semigroups on $Q_1$ and $Q_2$, respectively. Suppose that $A_1^*$ with domain $D(A^*_1)$ and $A_2^*$ with domain $D(A_2^*)$ are the generators of $S^*_1(\cdot)$ and $S^*_2(\cdot)$, respectively. By $(c)$ in \textbf{(CL)} and the same way to show \emph{Step 1} in the proof of
    Theorem \ref{yu-theorem-10-18-1}, we obtain that
\begin{equation}\label{yu-4-14-5}
  Q_1= PX_1^*\subset X_1^*,\;\;A^*P=PA^*\;\;\mbox{in}\;\;X_1^*
\end{equation}
and
\begin{equation}\label{yu-4-14-6}
\begin{cases}
    A^*_1=PA^*\;\;\mbox{in}\;\;Q_1\;\;\mbox{with}\;\;D(A^*_1)=PX_1^*=Q_1,\\
    A^*_2=(I-P)A^*\;\;\mbox{in}\;\;Q_2\;\mbox{with}\;\;D(A^*_2)=(I-P)X_1^*.
\end{cases}
\end{equation}
   According to $(d)$ in \textbf{(CL)}, there exist $\varepsilon>0$ and $C(\alpha,\varepsilon)>0$ so that
\begin{equation}\label{yu-4-14-7}
    \|S^*_2(t)\|_{\mathcal{L}(X)}\leq C(\alpha,\varepsilon)e^{-\varepsilon t}\;\;\mbox{for any}\;\;t\in\mathbb{R}^+.
\end{equation}
    Define
$E_1:=P^*X$ and $E_2:=\mbox{Ker}(P^*)$.
   By the similar arguments as those used to show \eqref{yu-12-12-6-01}, we conclude that
    $E_2=(I-P^*)X$, $X=E_1\oplus E_2$, $E_1=Q_1'$ and $E_2=Q_2'$
    in the following sense:  for any linear bounded functionals $F(\cdot)$ and $G(\cdot)$ on $Q_1$ and $Q_2$, respectively, there exist unique $f\in E_1$ and $g\in E_2$ so that $F(x)=\langle f,x\rangle_X$ for any $x\in Q_1$ and $G(x)=\langle g,x\rangle_X$ for any $x\in Q_2$.
    Let $S_1(\cdot):=P^*S(\cdot)$ and $S_2(\cdot):=(I-P^*)S(\cdot)$.  One can easily check that
    $S_1(\cdot)$ and $S_2(\cdot)$ are $C_0$-semigroups on $E_1$ and $E_2$, respectively. Moreover, they are the dual semigroups of $S_1^*(\cdot)$ and $S^*_2(\cdot)$, respectively. Suppose that
     $A_1$ with domain $D(A_1)$ and $A_2$ with domain $D(A_2)$ are the generators of $S_1(\cdot)$ and $S_2(\cdot)$, respectively. Then by the above facts, one can easily get that
    $A_1:=(A_1^*)^*$ and $A_2:=(A_2^*)^*$
    in the following sense:
\begin{equation*}\label{yu-4-14-20}
\begin{cases}
    \langle A_1f_1,g_1\rangle_{E_1,Q_1}=\langle f_1,A^*_1g_1\rangle_{E_1,Q_1}\;\;\mbox{for any}\;\;f_1\in D(A_1),\;g_1\in D(A_1^*),\\
     \langle A_2f_2,g_2\rangle_{E_2,Q_2}=\langle f_2,A^*_2g_2\rangle_{E_2,Q_2}
     \;\;\mbox{for any}\;\;f_2\in D(A_2),\;g_2\in D(A_2^*).
\end{cases}
\end{equation*}
    Meanwhile, by the same way to show \emph{Step 2} and  \emph{Step 3} in the proof of
    Theorem \ref{yu-theorem-10-18-1}, we have that
$E_1=P^*X_1\subset X_1$, $AP^*=P^*A$ in $X_1$,
\begin{equation*}\label{yu-4-15-2}
\begin{cases}
    A_1=P^*A\;\;\mbox{in}\;\;E_1\;\;\mbox{with}\;\;D(A_1)=P^*X_1=E_1,\\
    A_2=(I-P^*)A\;\;\mbox{in}\;\;E_2\;\;\mbox{with}\;\;D(A_2)=(I-P^*)X_1,
\end{cases}
\end{equation*}
    $B^*P\in \mathcal{L}(Q_1;U)$ and
    $(B^*P)^*=(\rho_0I-A_1)P^*(\rho_0I-\widetilde{A})^{-1}B\in \mathcal{L}(U;E_1)$
    in the sense: $\langle u, B^*P\varphi\rangle_U=\langle (\rho_0I-A_1)P^*(\rho_0I-\widetilde{A})^{-1}Bu, \varphi\rangle_{E_1,Q_1}$ for any $\varphi\in Q_1$ and $u\in U$.

\vskip 5pt
    \emph{Sub-step 3.1. We prove that $[A_1,(\rho_0I-A_1)P^*
    (\rho_0I-\widetilde{A})^{-1}B]$ is rapidly stabilizable in $E_1$ with bounded feedback law.}

    According to $(c)$ in \textbf{(CL)} and the similar arguments as those used to show \emph{Step 4} in Theorem \ref{yu-theorem-10-18-1},
    it suffices to prove \eqref{yu-11-13-4}, i.e., there  exists $C(\vartheta)>0$ so that
\begin{equation}\label{yu-4-15-4}
    \|\varphi\|^2_{Q_1}\leq C(\vartheta)
    \left(\|(\lambda I-A^*_1)\varphi\|^2_{Q_1}+\|B^*P\varphi\|_U^2\right)
    \;\;\mbox{for any}\;\;\lambda\in\mathbb{C}_{-\vartheta}^+,\;\;\varphi\in Q_1.
\end{equation}
    To this end, by \eqref{yu-4-14-1},  there  exists $C(\vartheta)>0$ so that
\begin{equation}\label{yu-4-15-5}
    \|(\rho_0I-A^*)^{\frac{1}{2}}\varphi\|_X^2\leq C(\vartheta)(\|(\rho_0I-A^*)^{\frac{1}{2}}(\lambda I-A^*)\varphi\|_X^2+\|B^*\varphi\|_U^2)
    \;\;\mbox{for any}\;\;\varphi\in X_{\frac{3}{2}}^*,\;\;\lambda\in\mathbb{C}_{-\vartheta}^+.
\end{equation}
    Since $P\varphi\in Q_1$ for any $\varphi\in X$, it follows from \eqref{yu-4-14-5} that $A^*(P\varphi)\in Q_1$ and then $(A^*)^2 (P\varphi)\in Q_1$ for any $\varphi\in X$. Hence,
\begin{equation}\label{yu-4-15-6}
    P\varphi\in X_2^*\;\;\mbox{if}\;\;\varphi\in X.
\end{equation}
     This, together with \eqref{yu-4-14-5} and \eqref{yu-4-14-6}, yields that for any $\varphi\in X$,
\begin{eqnarray}\label{yu-4-15-7}
   \|(\rho_0I-A^*)^{\frac{1}{2}}(\lambda I-A^*)P\varphi\|_X
    &\leq&
    \|(\rho_0I-A^*)^{-\frac{1}{2}}\|_{\mathcal{L}(X)}\|(\rho_0I-A^*)(\lambda I-A^*)P\varphi\|_X\nonumber\\
    &\leq&\|(\rho_0I-A^*)^{-\frac{1}{2}}\|_{\mathcal{L}(X)}\|(\rho_0I-A_1^*)(\lambda I-A_1^*)P\varphi\|_X\nonumber\\
    &\leq& \|(\rho_0I-A^*)^{-\frac{1}{2}}\|_{\mathcal{L}(X)}
    \|\rho_0I-A_1^*\|_{\mathcal{L}(Q_1)}\|(\lambda I-A_1^*)P\varphi\|_{Q_1},
\end{eqnarray}
   and
\begin{equation}\label{yu-4-15-8}
    \|P\varphi\|_{Q_1}=\|P\varphi\|_X\leq \|(\rho_0I-A^*)^{-\frac{1}{2}}\|_{\mathcal{L}(X)}
    \|(\rho_0I-A^*)^{\frac{1}{2}}P\varphi\|_X.
\end{equation}
    Therefore, by  \eqref{yu-4-15-5}, \eqref{yu-4-15-6}, \eqref{yu-4-15-7} and \eqref{yu-4-15-8}, we  conclude that there exists $C(\vartheta)>0$ so that \eqref{yu-4-15-4} holds and then  the pair $[A_1,(\rho_0I-A_1)P^*
    (\rho_0I-\widetilde{A})^{-1}B]$ is rapidly stabilizable in $E_1$ with bounded feedback law.

\vskip 5pt
    \emph{Sub-step 3.2. Let $\mu>0$ be fixed arbitrarily. By Sub-step 3.1, there exists $K:=K(\mu)\in \mathcal{L}(E_1;U)$ so that the solution
    $z_K(\cdot)$ of the equation (in $E_1$):
$z_t(t)=\left[A_1+(\rho_0I-A_1)P^*
    (\rho_0I-\widetilde{A})^{-1}BK\right]z(t)$ ($t\in\mathbb{R}^+$)
    satisfies
\begin{equation}\label{yu-4-15-11}
    \|z_K(t)\|_{E_1}\leq C(\mu)e^{-\mu t}\|z_K(0)\|_{E_1}\;\;\mbox{for any}\;\;t\in\mathbb{R}^+.
\end{equation}
Let
\begin{equation}\label{yu-4-15-bb01}
    F:=KP^*.
\end{equation}
    Then for each $y_0\in X$, the solution $y_F(\cdot;y_0)$ of the equation
    \eqref{yu-25-6-30-1}
 is in $C([0,+\infty);X)\cap L^2(\mathbb{R}^+;X)$.}

    To this end, let $y_0\in X$ be fixed arbitrarily. Let $T>0$ be chosen so  that
\begin{equation}\label{yu-4-15-12}
    \sqrt{2}C(\alpha,\varepsilon)e^{-\varepsilon T}<1.
\end{equation}
    Here and in what follows, $\varepsilon>0$ and $C(\alpha,\varepsilon)>0$ are given so that \eqref{yu-4-14-7} holds.
    It is clear that $F\in \mathcal{L}(X;U)$. Then it follows from Lemma~\ref{yu-proposition-25-6-26-2} that
    $y_F(\cdot;y_0)\in C([0,+\infty);X)$.  It suffices to show that $y_F(\cdot;y_0)\in L^2(\mathbb{R}^+;X)$, i.e.,
\begin{equation}\label{yu-4-15-13}
    \int_{\mathbb{R}^+}\|y_F(t;y_0)\|^2_Xdt<+\infty.
\end{equation}
   For this purpose, let $z_1(\cdot):=P^*y_F(\cdot;y_0)$ and $z_2(\cdot):=(I-P^*)y_F(\cdot;y_0)$.
   On one hand, by \eqref{yu-4-15-11} and the same way to prove \eqref{yu-11-20-7}, we get
\begin{equation}\label{yu-4-15-14}
    \|z_1(t)\|_X\leq C(\mu)e^{-\mu t}\|P^*y_0\|_X\;\;\mbox{for any}\;\;t\in\mathbb{R}^+.
\end{equation}
   By \eqref{yu-4-14-5}, $(A_2)$ and the same arguments as those in the proof of \eqref{yu-11-21-4}, we have that for each $n\in\mathbb{N}$,
\begin{equation}\label{yu-4-15-15}
    z_2(t)=S_2(t-nT)z_2(nT)
    +(I-P^*)\int_0^{t-nT}\widetilde{S}(t-nT-s)
    BFz_1(nT+s)ds\;\;\mbox{for all}\;\;t\in[nT,(n+1)T].
\end{equation}
    On the other hand, according to Assumptions $(A_1)$ and $(A_2)$, there exists $C(\gamma)>0$ so that for any $f\in C([0,T];X)$,
\begin{equation}\label{yu-4-24-5}
    \left\|\int_0^T\widetilde{S}(T-s)BFf(s)ds\right\|_X
    \leq C(\gamma)T^{1-\gamma}\|(\rho_0 I-\widetilde{A})^{-\gamma}B\|_{\mathcal{L}(U;X)}
    \|F\|_{\mathcal{L}(X;U)}\sup_{t\in[0,T]}\|f(t)\|_X.
\end{equation}
    Then by \eqref{yu-4-14-7}, \eqref{yu-4-15-14}, \eqref{yu-4-15-15}, \eqref{yu-4-24-5} and the fact that $S_2(\cdot)$ is the dual semigroup of $S_2^*(\cdot)$, we have
\begin{equation}\label{yu-4-15-16}
    \sup_{t\in[0,T]}\|z_2(t)\|^2_X\leq C(\alpha,\varepsilon,\mu,T)\|y_0\|_X^2,
\end{equation}
\begin{equation}\label{yu-4-16-1}
    \sup_{t\in[nT,(n+1)T]}\|z_2(t)\|_X^2
    \leq C(\alpha,\varepsilon,\mu,T)(\|z_2(nT)\|_X^2+e^{-2\mu nT}\|y_0\|_X^2)
\end{equation}
    and
\begin{equation}\label{yu-4-16-2}
    \|z_2((n+1)T)\|_X^2\leq 2(C(\alpha,\varepsilon))^2e^{-2\varepsilon T}\|z_2(nT)\|_X^2
    +C(\mu,T,\mu)e^{-2\mu nT}\|y_0\|^2_X
\end{equation}
    for each $n\in\mathbb{N}^+$. According to \eqref{yu-4-15-12}, \eqref{yu-4-15-16} and \eqref{yu-4-16-2}, it is clear that
$\sum_{n=1}^{+\infty}\|z_2(nT)\|_X^2<+\infty$.
    This, along with \eqref{yu-4-15-16} and \eqref{yu-4-16-1}, implies that
\begin{equation}\label{yu-4-16-4}
    \int_{\mathbb{R}^+}\|z_2(t)\|_X^2dt\leq C(\alpha,\varepsilon,\mu,T)
    \left(\sum_{n=1}^{+\infty}\|z_2(nT)\|_X^2+(1-e^{-2\mu T})^{-1}\|y_0\|_X^2\right)<+\infty.
\end{equation}
    Therefore, by \eqref{yu-4-15-14}, \eqref{yu-4-16-4} and the definitions of $z_1(\cdot)$ and $z_2(\cdot)$, we obtain \eqref{yu-4-15-13}.
\par
    By \emph{Sub-step 3.1} and using the same arguments as those in \emph{Step 6} of
    the proof of Theorem \ref{yu-theorem-10-18-1}, we conclude that if $F$ is given by \eqref{yu-4-15-bb01}, then there exist $\delta>0$ and $C(\delta)>0$ so that for each $y_0\in X$, the  solution $y_F(\cdot;y_0)$ of the closed-loop system \eqref{yu-25-6-30-1} satisfies
$\|y_F(t;y_0)\|_X\leq C(\delta)e^{-\delta t}\|y_0\|_X$ for each $t\in\mathbb{R}^+$.
    Hence, $(iii)$ holds.

\vskip 5pt
    \emph{Step 4. We prove $(iii)\Rightarrow (iv)$.}

    When $\gamma\in[0,\frac{1}{2})$,  the result follows easily from
    \cite[Theorem 2]{Kunisch-Wang-Yu} and
    the fact that $[A,B]$ (in $X$) and  $[\mathcal{A},B]$ (in $\mathcal{X}$) are consistent in this case. It suffices to show $(iii)\Rightarrow (iv)$ for  the case $\gamma\in[\frac{1}{2},1)$.

    Suppose that $(iii)$ holds, i.e., there exist $\beta>0$, $C(\beta)>0$ and $F\in \mathcal{L}(X;U)$ so that for each $y_0\in X$, the solution $y_F(\cdot;y_0)$ of the equation \eqref{yu-25-6-30-1}  satisfies
\begin{equation}\label{yu-25-7-12-1}
    \|y_F(t;y_0)\|_X\leq C(\beta)e^{-\beta t}\|y_0\|_X\;\;\mbox{for each}\;\;t\in\mathbb{R}^+.
\end{equation}
    Let $u_F(\cdot):=Fy_F(\cdot;y_0)$. Then by \cite[Theorem (Ball) on Page 259]{Pazy}, we have that for any $\varphi\in X_1^*$,
\begin{equation}\label{yu-4-9-1-bb}
    \frac{d}{dt}
    \langle y_F(t;y_0),\varphi\rangle_X=\langle y_F(t;y_0),A^*\varphi\rangle_X
    +\langle u_F(t),B^*\varphi\rangle_U\;\;\mbox{for all}\;\;t\in\mathbb{R}^+.
\end{equation}
 According to \eqref{yu-25-7-12-1}, it is clear that for any $\lambda\in\mathbb{C}_{-\beta}^+$,
\begin{eqnarray}\label{yu-4-9-3}
    \xi(\lambda):=\int_{\mathbb{R}^+}e^{-\lambda t}y_F(t;y_0)dt\;\;\mbox{and}\;\;\eta(\lambda):=\int_{\mathbb{R}^+}e^{-\lambda t}u_F(t)dt
\end{eqnarray}
        are well defined, and
\begin{equation}\label{yu-4-9-11-bb}
\begin{cases}
    \|\xi(\lambda)\|_X\leq C(\beta)
    (\mbox{Re}\lambda+\beta)^{-1}\|y_0\|_X,\\
    \|\eta(\lambda)\|_U\leq C(\beta)\|F\|_{\mathcal{L}(X;U)}(\mbox{Re}\lambda+\beta)^{-1}\|y_0\|_X.
\end{cases}
\end{equation}
     These, along with \eqref{yu-4-9-1-bb}, imply that  for each  $\lambda\in \mathbb{C}_{-\beta}^+$ and $\varphi\in X^*_1$, $\int_0^{+\infty}e^{-\lambda t}\frac{d}{dt}
    \langle y_F(t;y_0),\varphi\rangle_Xdt$ is well defined and
\begin{eqnarray}\label{yu-4-9-13}
    &\;&\int_{\mathbb{R}^+}e^{-\lambda t}\frac{d}{dt}
    \langle y_F(t;y_0),\varphi\rangle_Xdt
    =\langle \lambda\xi(\lambda),\varphi\rangle_X
    -\langle y_0,\varphi\rangle_X.
\end{eqnarray}
    Therefore, by \eqref{yu-4-9-1-bb}, \eqref{yu-4-9-3} and \eqref{yu-4-9-13}, we obtain that for any $\varphi\in X_1^*$ and
    $\lambda\in\mathbb{C}_{-\beta}^+$,
$\langle y_0,\varphi\rangle_X=\langle \xi(\lambda),(\bar{\lambda}I-A^*)\varphi\rangle_X
    -\langle\eta(\lambda),B^*\varphi\rangle_U$.
    Since $y_0$ is arbitrary, it follows from the latter equality, \eqref{yu-25-7-12-1} and \eqref{yu-4-9-11-bb} that
$\|\varphi\|_X^2\leq C(\beta,F)(\mbox{Re}\lambda+\beta)^{-2}
    (\|(\lambda I-A^*)\varphi\|_X^2+\|B^*\varphi\|_U^2)$ for any $\varphi\in X_1^*$ and $\lambda\in\mathbb{C}_{-\beta}^+$.
    Hence, \eqref{yu-4-2-2-b} holds by letting $\alpha:=\frac{\beta}{2}$  and $C(\alpha):=
    4C(\beta,F)\beta^{-2}$.
\vskip 5pt
    \emph{Step 5. We prove $(iv)\Rightarrow (i)$.}

    When $\gamma\in[0,\frac{1}{2})$, the result follows easily from
    \cite[Theorem 2]{Kunisch-Wang-Yu}. Thus, it suffices to show $(i)$ for the case $\gamma\in[\frac{1}{2},1)$.

    Suppose that $(iv)$ holds, i.e., there exist $\alpha>0$ and $C(\alpha)>0$ so that \eqref{yu-4-2-2-b} is true.
    By Assumption $(A_3)$ and $(ii)$ of Remark \ref{yu-remark-12-22-1}, we have that $(a)$-$(d)$ in \textbf{(CL)} hold.
    It follows from \eqref{yu-3-29-13} and \eqref{yu-3-29-14} that for any $f\in \mathcal{X}^*_1$,
$(\rho_0I-A^*)^{-\frac{1}{2}}(\rho_0 I-\widetilde{A})^{-\frac{1}{2}}f\in X_{\frac{3}{2}}^*(\subset X_1^*)$,
\begin{eqnarray*}\label{yu-4-16-21}
    \|(\lambda I-A^*)(\rho_0I-A^*)^{-\frac{1}{2}}(\rho_0 I-\widetilde{A})^{-\frac{1}{2}}f\|_X
    &=&\|(\rho_0I-A^*)^{-\frac{1}{2}}(\rho_0I-\widetilde{A})^{-\frac{1}{2}}(\lambda I-\mathcal{A}^*)f\|_X\nonumber\\
    &\leq&\|(\rho_0I-A^*)^{-\frac{1}{2}}\|_{\mathcal{L}(X)}
    \|(\lambda I-\mathcal{A}^*)f\|_{\mathcal{X}}\;\;\mbox{for all}\;\;\lambda\in\mathbb{C}
\end{eqnarray*}
    and
$\|B^*(\rho_0I-A^*)^{-\frac{1}{2}}(\rho_0 I-\widetilde{A})^{-\frac{1}{2}}f\|_U=\|\mathcal{B}^* f\|_U$.
    These, along with \eqref{yu-4-2-2-b} (where $\varphi=(\rho_0I-A^*)^{-\frac{1}{2}}(\rho_0 I-\widetilde{A})^{-\frac{1}{2}}f$, $f\in \mathcal{X}^*_1$), imply that
\begin{equation}\label{yu-4-11-2}
    \|(\rho_0I-A^*)^{-\frac{1}{2}}(\rho_0I-\widetilde{A})^{-\frac{1}{2}}f\|^2_X
    \leq D(\alpha)(\|(\lambda I-\mathcal{A}^*)f\|_{\mathcal{X}}^2+\|\mathcal{B}^*f\|^2_U)
    \;\;\mbox{for any}\;\;\lambda\in \mathbb{C}_{-\alpha}^+,\;\;f\in \mathcal{X}^*_1,
\end{equation}
    where $D(\alpha):=C(\alpha)(1+\|(\rho_0I-A^*)^{-\frac{1}{2}}\|^2_{\mathcal{L}(X)})$. Let $Q_1:=Q_1(\alpha)$ and
    $Q_2:=Q_2(\alpha)$ be given in \textbf{(CL)}, and $W_1:=W_1(\alpha)$ and $W_2:=W_2(\alpha)$ be defined by \eqref{yu-4-2-6}.
    By $(e_1)$ in \emph{Step 1}, we  define
$\mathcal{P}: \mathcal{X}\to W_1$:
$\mathcal{P}g=g_1$ if $g=g_1+g_2$ with $(g_1,g_2)\in W_1\times W_2$.
    By the same arguments as those in the proof of Theorem \ref{yu-theorem-10-18-1}, we have that
   $\mathcal{P}$ is a projection operator, and $W_1$ and $W_2$ are two Hilbert spaces whose norms are inherited from $\mathcal{X}$.

    We claim that
\begin{equation}\label{yu-4-11-4}
    \mathcal{P}=(\rho_0I-\widetilde{A})^{\frac{1}{2}}P(\rho_0I-\widetilde{A})^{-\frac{1}{2}},
\end{equation}
    where $P$ is defined by the same way as \eqref{yu-4-14-4}.
    For this purpose, we take $f\in \mathcal{X}$ arbitrarily. According to the definition of $\mathcal{P}$, there  exists a unique $(f_1,f_2)\in W_1\times W_2$ so  that
    $f=f_1+f_2$ and $\mathcal{P}f=f_1$.
    By $(ii)$ in Remark \ref{yu-remark-6-22-1} and \eqref{yu-4-2-6}, we have
$(\rho_0I-\widetilde{A})^{-\frac{1}{2}}f_1\in Q_1$ and $(\rho_0I-\widetilde{A})^{-\frac{1}{2}}f_2\in Q_2$.
    These imply that
$P(\rho_0I-\widetilde{A})^{-\frac{1}{2}}f=(\rho_0I-\widetilde{A})^{-\frac{1}{2}}f_1$
    and then
$(\rho_0I-\widetilde{A})^{\frac{1}{2}}P(\rho_0I-\widetilde{A})^{-\frac{1}{2}}f=f_1=\mathcal{P}f$.
    Since $f$ is arbitrary, \eqref{yu-4-11-4} follows immediately.
\par
    Define $A^*_1:=PA^*$, $\mathcal{S}^*_2(\cdot):=(I-\mathcal{P})\mathcal{S}^*(\cdot)$ and $\mathcal{A}^*_2:=(I-\mathcal{P})\mathcal{A}^*$.
    It is clear that $\mathcal{S}^*_2(\cdot)$ is a $C_0$-semigroup on $W_2$. By \emph{Step 1} and the same way to show \emph{Step 1} in the proof of  Theorem \ref{yu-theorem-10-18-1}, we obtain that
\begin{equation}\label{yu-4-12-4}
    A_1^*\in\mathcal{L}(Q_1),\;\;PA^*=A^*P\;\;\mbox{in}\;\;X_1^*,
\end{equation}
\begin{equation}\label{yu-4-12-4-b}
    \mathcal{P}\mathcal{A}^*=\mathcal{A}^* \mathcal{P}\;\;\mbox{in}\;\;\mathcal{X}^*_1
\end{equation}
    and
\begin{equation}\label{yu-4-12-5}
    \mathcal{A}^*_2\;\;\mbox{with}\;\;D(\mathcal{A}^*_2):=(I-\mathcal{P})\mathcal{X}^*_1\;\;\mbox{is the generator of}
    \;\;\mathcal{S}^*_2(\cdot).
\end{equation}
    According to $(e_4)$ in \emph{Step 1}, \cite[Theorem 5.3 and Remark 5.4, Section 1.5, Chapter 1]{Pazy} and \eqref{yu-4-12-5}, there  exist $\nu>0$ and  $C(\nu)>0$ so that
$\|f\|_{W_2}\leq C(\nu)\|(\lambda I-\mathcal{A}^*_2)f\|_{W_2}$ for any $\lambda\in\mathbb{C}_{-\nu}^+$ and $f\in D(\mathcal{A}_2^*)$.
    This, together with \eqref{yu-4-12-4-b}, yields that
\begin{eqnarray}\label{yu-4-12-7}
    \|(I-\mathcal{P})f\|_\mathcal{X}&\leq& C(\nu)\|(\lambda I-\mathcal{A}^*)(I-\mathcal{P})f\|_{\mathcal{X}}
    =C(\nu)\|(I-\mathcal{P})(\lambda I-\mathcal{A}^*)f\|_{\mathcal{X}}\nonumber\\
    &\leq& C(\nu,\mathcal{P})\|(\lambda I-\mathcal{A}^*)f\|_{\mathcal{X}}\;\;\mbox{for any}\;\;\lambda\in\mathbb{C}_{-\nu}^+,\;\;f\in \mathcal{X}^*_1.
\end{eqnarray}
    (Here, by the closed graph theorem, we observe that $\mathcal{P}$ is bounded.)
    Moreover, for each $f\in \mathcal{X}^*_1$, i.e., $f=(\rho_0I-\widetilde{A})^{\frac{1}{2}}g$ for some $g\in X_1^*$ (see \eqref{yu-3-29-bbb-1}), it follows from  \eqref{yu-25-6-23-5-b}, \eqref{yu-4-11-4} and \eqref{yu-4-12-4} that
\begin{eqnarray*}\label{yu-4-12-8}
    \|\mathcal{P}f\|_{\mathcal{X}}&=&\|(\rho_0I-\widetilde{A})^{\frac{1}{2}}Pg\|_{\mathcal{X}}=\|Pg\|_X
    =\|(\rho_0I-A^*_1)^{\frac{1}{2}}P(\rho_0I-A^*)^{-\frac{1}{2}}g\|_{Q_1}\nonumber\\
    &\leq& \|P\|_{\mathcal{L}(X)}\|(\rho_0I-A^*_1)^{\frac{1}{2}}\|_{\mathcal{L}(Q_1)}
    \|(\rho_0I-A^*)^{-\frac{1}{2}}(\rho_0I-\widetilde{A})^{-\frac{1}{2}}f\|_X.
\end{eqnarray*}
    This, along with \eqref{yu-4-11-2}, implies that
\begin{equation*}
    \|\mathcal{P}f\|^2_{\mathcal{X}}\leq (D(\alpha))^2 \|P\|^2_{\mathcal{L}(X)}\|(\rho_0I-A^*_1)^{\frac{1}{2}}\|^2_{\mathcal{L}(Q_1)}(\|(\lambda I-\mathcal{A}^*)f\|_X^2+\|\mathcal{B}^*f\|^2_U)
    \;\;\mbox{for any}\;\;\lambda\in \mathbb{C}_{-\alpha}^+,\;\;f\in \mathcal{X}^*_1,
\end{equation*}
    which, combined with \eqref{yu-4-12-7}, indicates that
$\|\varphi\|_{\mathcal{X}}^2\leq C(\beta)(\|(\lambda I-\mathcal{A}^*)\varphi\|_{\mathcal{X}}^2+
    \|\mathcal{B}^*\varphi\|_U^2)$ for any $\varphi\in \mathcal{X}^*_1$ and $\lambda\in\mathbb{C}_{-\beta}^+$,
     where $\beta:=\min\{\nu,\alpha\}$ and $C(\beta):=2((C(\nu,\mathcal{P}))^2+(D(\alpha))^2\|P\|^2_{\mathcal{L}(X)} \|(\rho_0I-A^*_1)^{\frac{1}{2}}\|^2_{\mathcal{L}(Q_1)})$. Therefore, by Proposition \ref{yu-proposition-6-27-1} and \emph{Step 1}, we  apply \cite[Theorem 2]{Kunisch-Wang-Yu} to the pair $[\mathcal{A},B]$ and complete the proof of $(i)$.
    In summary, we finish the proof of Theorem~\ref{yu-theorem-3-31-1}.
\end{proof}
\begin{remark}\label{yu-remark-25-7-13-10}
    The algorithmic procedure to construct the feedback is provided in \emph{Sub-step 3.2}. We find that this algorithmic procedure is consistent with the procedure given in Theorem \ref{yu-theorem-10-18-1} (both of which are based on a controllable system formed by bounded state and control operators). The details can be found in $(iii)$ of Remark \ref{yu-remark-4-14-1}. We give the following table to illustrate the main differences between two frameworks caused by $\gamma$:
    \begin{center}
\begin{tabular}{|c|c|c|c|c|}
\hline
 & SS & SO & OB & PJ\\ \hline
$\gamma\in [0,\frac{1}{2})$ & $X$ &$A$ &$B^*$ & P\\ \hline
$\gamma\in [\frac{1}{2},1)$ &$X_{-\frac{1}{2}}$ & $\widetilde{A}_{-\frac{1}{2}}$ &$B^*(\rho_0I-A^*)^{-\frac{1}{2}}(\rho_0I-\widetilde{A})^{-\frac{1}{2}}$ &$(\rho_0I-\widetilde{A})^{\frac{1}{2}}P(\rho_0I-\widetilde{A})^{-\frac{1}{2}}$\\ \hline
\end{tabular}
\end{center}
    Here, $SS:=$state space, $SO:=$state operator, $OB:=$observe operator, $PJ:=$projection operator from $\mathcal{X}$
    to $W_1$ (see \eqref{yu-4-2-6}) and $P$ is defined by \eqref{yu-4-14-4}.
\end{remark}
\begin{proof}[Proof of Corollary \ref{yu-corollary-4-22-1}]
    We only need to show the following claims:
     $(a)$ the pair $[\mathcal{A},B]$ is rapidly stabilizable in $\mathcal{X}$ if and only if $[\mathcal{A}+\alpha I,B]$ is stabilizable in $\mathcal{X}$ for each $\alpha>0$;
       $(b)$ when $\gamma\in[\frac{1}{2},1)$, $\mathcal{A}+\alpha I=(\rho_0I-\widetilde{A})^{\frac{1}{2}}(A+\alpha I)(\rho_0I-\widetilde{A})^{-\frac{1}{2}}$ for any $\alpha\in\mathbb{R}^+$.
     For $(a)$, we refer to the proof of $(ii)$ of Lemma \ref{yu-prop-5-22-1} with $[A,B]$ (in $X$) replaced by $[\mathcal{A},B]$ (in $\mathcal{X}$). For $(b)$, it follows from \eqref{yu-25-6-23-5} immediately.
\end{proof}

\section{Examples}\label{yu-sec-7-15-100}

\subsection{Parabolic equation on a bounded domain with boundary controls}\label{yu-sec-25-7-1}
\begin{lemma}\label{yu-lemma-4-22-1}
     Suppose that
    $(l_1)$ the spectrum of $L$ with domain $D(L)$ consists of isolated eigenvalues with finite algebraic multiplicity; $(l_2)$  the spectrum of $L$ has no finite cluster point. If $L$ generates an analytic semigroup $S_L(\cdot)$ in $X$, then  $L$ satisfies (AEDC).
\end{lemma}
\begin{proof}
    Since $L$ generates an analytic semigroup, $\mathbb{C}_{-\alpha}^+\cap \sigma(L)$ is bounded in $\mathbb{C}$ for any $\alpha>0$.  This, along with Assumptions $(l_1)$ and $(l_2)$, implies that  the spectrum $\sigma(L)$ of $L$ is countable (or finite or empty) and consists of poles of $R(\cdot):=(\cdot I-L)^{-1}$ with finite algebraic multiplicity only and $\sharp(\sigma(L)\cap \mathbb{C}_{\mu}^+)<+\infty$ for each $\mu\in\mathbb{R}$.
   Thus, by  \cite[Corollary 3.12, Section 3, Chapter IV]{Engel-Nagel} and \cite[Theorem 6.17, Section 6, Chapter 3]{Kato}. we can claim that
    the conclusions in Lemma \ref{yu-lemma-6-23-1} hold. Then by repeating the arguments of the first case in Lemma \ref{yu-lemma-12-20-1-bb},
    we conclude that $L$ satisfies (AEDC).
\end{proof}

   By Theorem \ref{yu-theorem-3-31-1}, Lemma \ref{yu-lemma-4-22-1} and \cite[Theorem 1.6]{Badra-Takahashi}, we have the following corollary.

\begin{corollary}\label{yu-proposition-4-22-1}
    Suppose that Assumptions $(A_1)$ and $(A_2)$ hold, and $A$ satisfies Assumptions
    $(l_1)$ and  $(l_2)$  in Lemma \ref{yu-lemma-4-22-1} (by replacing $L$ by $A$). The following statements are equivalent:
\begin{enumerate}
\item[$(i)$] The  pair $[\mathcal{A},B]$ is stabilizable in $\mathcal{X}$.
\item[$(ii)$] The  pair $[A,B]$ is stabilizable in $X$ with bounded feedback law.
  \item [$(iii)$] There exist $\alpha>0$ and $C(\alpha)>0$ so that \eqref{yu-4-2-2-b} holds.
  \item [$(iv)$] There exists $\alpha>0$ so that the following statement holds: if there  exist $\lambda\in\mathbb{C}_{-\alpha}^+$ and $\varphi\in X_1^*$ so that
  $(\lambda I-A^*)\varphi=0$ and $B^*\varphi=0$, then $\varphi=0$.
\end{enumerate}
\end{corollary}
\begin{remark}
   Corollary \ref{yu-proposition-4-22-1} says that the criterion for stabilizability property for the pair  $[\mathcal{A},B]$ on the state space $\mathcal{X}$ given by \cite{Kunisch-Wang-Yu} within the framework of Definition \ref{yu-definition-10-18-1}  is equivalent to the criterion  given by \cite{Badra-Takahashi, Raymond} for the pair $[A,B]$ in $X$.
    Moreover, Corollary \ref{yu-proposition-4-22-1} theoretically explains why in \cite{Badra-Takahashi, Raymond}, although the authors used bounded feedback laws to stabilize the pair $[A,B]$ in $X$, the requirement for the stabilizability of this system was not strengthened.
\end{remark}
Since $[\mathcal{A},B]$ is rapidly stabilizable in $\mathcal{X}$ if and only if $[\mathcal{A}+\alpha I,B]$ is stabilizable in $\mathcal{X}$ for each $\alpha>0$ (see $(a)$ in the proof of Corollary \ref{yu-corollary-4-22-1}), by the same arguments as those in Corollary \ref{yu-corollary-4-22-1}, Lemma \ref{yu-lemma-4-22-1} and Corollary \ref{yu-proposition-4-22-1}, we have the following result.
\begin{corollary}
    Under assumptions of Corollary \ref{yu-proposition-4-22-1}, the following statements are equivalent:
    \begin{enumerate}
\item[$(i)$] The  pair $[\mathcal{A},B]$ is rapidly  stabilizable in $\mathcal{X}$.
\item[$(ii)$] The  pair $[A,B]$ is rapidly  stabilizable in $X$ with bounded feedback law.
  \item [$(iii)$] For each $\alpha>0$, there exists $C(\alpha)>0$ so that \eqref{yu-4-2-2-b} holds.
  \item [$(iv)$] If there exist $\lambda\in\mathbb{C}$ and $\varphi\in X_1^*$ so that
  $(\lambda I-A^*)\varphi=0$ and $B^*\varphi=0$, then $\varphi=0$.
\end{enumerate}
\end{corollary}

    We next consider  the following controlled heat equation with Dirichlet boundary control:
\begin{equation}\label{yu-4-16-40}
\begin{cases}
    y_t(t,x)=\triangle y(t,x),&(t,x)\in\mathbb{R}^+\times \Omega,\\
    y(t,x)=\chi_{\Gamma}(x)u(t,x),&(t,x)\in\mathbb{R}^+\times \partial\Omega,
\end{cases}
\end{equation}
    where  $\Omega$ is a bounded and open region in $\mathbb{R}^d$ with smooth boundary $\partial\Omega$, $\Gamma$ is
    a nonempty and open subset of $\partial\Omega$, $\chi_\Gamma$ is the characteristic function of $\Gamma$ and
    $u\in L^2(\mathbb{R}^+;L^2(\partial\Omega;\mathbb{C}))$.
\par
    Let
    $\triangle_D$ be the Laplacian operator in $L^2(\Omega;\mathbb{C})$ with homogenous Dirichlet boundary condition, i.e., $\triangle_D:=\triangle$ with domain $D(\triangle_D):=H_0^1(\Omega;\mathbb{C})\cap H^2(\Omega;\mathbb{C})$. It is well known that
    $(\triangle_D,D(\triangle_D))$ is negative, and generates an analytic and compact semigroup on $L^2(\Omega;\mathbb{C})$. By the classical way, $\triangle_D$ can be extended into $H^{-1}(\Omega;\mathbb{C})$ (denoted by $\widehat{\triangle_D}$) and $D(\triangle_D)'$ (denoted by $\widetilde{\triangle_D}$) with domains $H_0^1(\Omega;\mathbb{C})$ and
    $L^2(\Omega;\mathbb{C})$, respectively, where $V'$ means the dual space $V$ w.r.t. the pivot space $L^2(\Omega;\mathbb{C})$.
Introduce the Dirichlet map $\mathcal{D}_\Gamma: L^2(\partial\Omega;\mathbb{C})\to L^2(\Omega;\mathbb{C})$ as follows: $\mathcal{D}_\Gamma v=\varphi$, where $\varphi$ is the unique solution of the following Laplace equation:
\begin{equation*}\label{yu-4-16-41}
\begin{cases}
    \triangle\varphi=0&\mbox{in}\;\;\Omega,\\
    \varphi=\chi_{\Gamma}v&\mbox{on}\;\;\partial\Omega.
\end{cases}
\end{equation*}
    It is well known that $\mathcal{D}_\Gamma\in \mathcal{L}(L^2(\partial\Omega;\mathbb{C});H^{\frac{1}{2}}(\Omega;\mathbb{C}))$. Let
\begin{equation}\label{yu-4-16-42}
    B:=-\widetilde{\triangle_D}\mathcal{D}_\Gamma.
\end{equation}
    Then $B\in\mathcal{L}(L^2(\partial\Omega;\mathbb{C});D(\triangle_D)')$. Indeed, one can show
\begin{equation}\label{yu-4-16-43}
    B\in\mathcal{L}(L^2(\partial\Omega;\mathbb{C});D((-\triangle_D)^{\frac{3}{4}+\varepsilon})')\;\;\mbox{for each}\;\;\varepsilon>0
\end{equation}
    (see \cite[Section 3.1, Chapter 3]{Lasiecka-Triggiani-2000} or \cite[Section 5.1.4.2]{Trelat}).
\begin{remark}\label{yu-25-7-13-1}
  According to \cite[Section 3.1, Chapter 3]{Lasiecka-Triggiani-2000} (see also \cite[Lemma 5.8]{Trelat}) and
    \eqref{yu-4-16-43}, $B$ (defined by \eqref{yu-4-16-42}) is not admissible (w.r.t. $\triangle_D$). This yields that it is inappropriate to use $L^2(\Omega;\mathbb{C})$ as the state space when studying the null controllability of the open-loop system \eqref{yu-4-16-40}. In order to make $B$ admissible,
    the usual approach is to extend the state space $L^2(\Omega;\mathbb{C})$ to $H^{-1}(\Omega;\mathbb{C})$ and extend the state operator $\triangle_D$ from $L^2(\Omega;\mathbb{C})$ to $H^{-1}(\Omega;\mathbb{C})$, i.e., use $\widehat{\triangle_D}$ as the state operator (see, for example, \cite[Section 5.2.5.2, Chapter 5]{Trelat}).
    In this framework, $B$ is admissible (w.r.t. $\widehat{\triangle_D}$) and it has been proved that system \eqref{yu-4-16-40} is null controllable in $H^{-1}(\Omega;\mathbb{C})$, i.e., for any $T>0$ and $y_0\in H^{-1}(\Omega;\mathbb{C})$, there  exists a control $u\in L^2(0,T;L^2(\partial\Omega;\mathbb{C}))$ so that the solution $y(\cdot;y_0,u)$ of the system \eqref{yu-4-16-40} with initial data $y(0)=y_0$ satisfies $y(T;y_0,u)=0$ (see \cite[Proposition 11.5.4]{Tucsnak-Weiss}) .
    Certainly, this conclusion implies that for any $T>0$ and $y_0\in L^2(\Omega;\mathbb{C})$, there  exists a control $u\in L^2(0,T;L^2(\partial\Omega;\mathbb{C}))$ so that the solution $y(\cdot;y_0,u)$ of the system \eqref{yu-4-16-40} with initial data $y(0)=y_0$ satisfies $y(T;y_0,u)=0$, but we cannot guarantee that $y(t;y_0,u)\in L^2(\Omega;\mathbb{C})$ for each $t\in[0,T]$.
\end{remark}
    It is clear that
     the Laplacian operator $\triangle_D$ with domain $D(\triangle_D):=H_0^1(\Omega;\mathbb{C})\cap H^2(\Omega;\mathbb{C})$
     satisfies Assumptions $(l_1)$ and $(l_2)$ in Lemma \ref{yu-lemma-4-22-1}. Thus, by \eqref{yu-4-16-43}, $(i)$ in Remark \ref{yu-remark-7-3-1}, Lemma \ref{yu-lemma-4-22-1}
    and \cite[Proposition 11.5.4, Chapter 11]{Tucsnak-Weiss}, we get the following corollary.
\begin{corollary}
    For each $\alpha>0$, there exist $C(\alpha)>0$ and $F:=F(\alpha)\in \mathcal{L}(L^2(\Omega;\mathbb{C});L^2(\partial\Omega;\mathbb{C}))$ so that any solution
    $y_F(\cdot,\cdot)$ of the following equation:
\begin{equation*}\label{yu-4-23-2}
\begin{cases}
    y_t(t,x)=\triangle y(t,x),&(t,x)\in\mathbb{R}^+\times \Omega,\\
    y(t,x)=\chi_{\Gamma}(x)[Fy(t)](x),&(t,x)\in\mathbb{R}^+\times \partial\Omega,
\end{cases}
\end{equation*}
    satisfies that $\|y_F(t)\|_{L^2(\Omega;\mathbb{C})}\leq C(\alpha)e^{-\alpha t}\|y_F(0)\|_{L^2(\Omega;\mathbb{C})}$ for any $t\in\mathbb{R}^+$.
\end{corollary}

\subsection{Parabolic equation with differential-type control operators}

   Denote by $\mathcal{S}(\mathbb{R}^n;\mathbb{C})$  the $\mathbb{C}$-valued Schwartz space on $\mathbb{R}^n$ ($n\in\mathbb{N}^+$). Its dual space  is the space of tempered distributions (on $\mathbb{R}^n$), which is denoted by $\mathcal{S}'(\mathbb{R}^n;\mathbb{C})$.
     Let $s\in[0,1)$. We introduce the fractional-power Laplacian operator in $\mathbb{R}^n$ as follows:
\begin{equation}\label{yu-4-24-2-wang}
    (-\triangle)^sf(x):=\frac{1}{(2\pi)^n}\int_{\mathbb{R}^n}\int_{\mathbb{R}^n}e^{\mathrm{i}\langle x-y,\xi\rangle}
   |\xi|^{2s}f(y)dyd\xi,\;\;x\in\mathbb{R}^n.
\end{equation}
     It is clear that
    $\mathcal{F}[(-\triangle)^sf]=|\xi|^{2s}\mathcal{F}[f]$ for each $f\in \mathcal{S}(\mathbb{R}^n;\mathbb{C})$, where $\mathcal{F}[f]$ is the Fourier transform of $f$ (see \cite{Zworski}).
    Moreover, by \eqref{yu-4-24-2-wang} and Plancherel's identity, we have that for any $f\in\mathcal{S}(\mathbb{R}^n;\mathbb{C})$,
\begin{equation*}\label{yu-4-25-5}
    \|(-\triangle)^sf\|_{H^{-2s}(\mathbb{R}^n;\mathbb{C})}^2
    =(2\pi)^{-n}\|(1+|\xi|^2)^{-s}|\xi|^{2s}\mathcal{F}[f]\|_{L^2(\mathbb{R}^n;\mathbb{C})}^2
    \leq \|f\|^2_{L^2(\mathbb{R}^n;\mathbb{C})}.
\end{equation*}
    This, along with the density of $\mathcal{S}(\mathbb{R}^n;\mathbb{C})$ in $L^2(\mathbb{R}^n;\mathbb{C})$, implies that
    $(-\triangle)^s$ can be extended uniquely to $L^2(\mathbb{R}^n;\mathbb{C})$ as a bounded linear operator from $L^2(\mathbb{R}^n;\mathbb{C})$ to
    $H^{-2s}(\mathbb{R}^n;\mathbb{C})$. In what follows, we still denote this extended  operator by the same manner.

    Let $X=L^2(\mathbb{R}^n;\mathbb{C})$, $U= L^2(\mathbb{R}^n;\mathbb{C})$ and $\omega$ be a Lebesgue measurable set in $\mathbb{R}^n$.
    Consider the following controlled equation:
\begin{equation}\label{yu-4-24-2}
    y_t(t,x)=\triangle y(t,x)+(-\triangle)^s(a_\omega(x)u(t,x)),\;\;(t,x)\in\mathbb{R}^+\times \mathbb{R}^n,
\end{equation}
    where $u\in L^2(\mathbb{R}^+;U)$ and $a_\omega\in L^\infty(\mathbb{R}^n;\mathbb{R})\cap C^\infty(\mathbb{R}^n;\mathbb{R})$ satisfies
\begin{equation}\label{yu-4-24-1}
    a_\omega(x)\geq c>0\;\;\mbox{for each}\;\;x\in\omega.
\end{equation}
    Let $A=\triangle$ with domain $D(A):=H^2(\mathbb{R}^n;\mathbb{C})$ and
\begin{equation}\label{yu-4-24-4}
    [Bu](x):=(-\triangle)^s(a_\omega(x)u(x)),\;\;x\in\mathbb{R}^n,\;\;u\in U.
\end{equation}
    It is well known that $A$ with domain $H^2(\mathbb{R}^n;\mathbb{C})$ is normal and generates an analytic semigroup
    $e^{\triangle\cdot}$ in $L^2(\mathbb{R}^n;\mathbb{C})$.  Therefore, Assumption $(A_1)$ holds, and by the fact that $\sigma(\triangle)
    \subset(-\infty,0]$, $(b)$ in $(iv)$ of Remark \ref{yu-remark-12-22-1} is true.
    Moreover, since $(-\triangle)^s\in \mathcal{L}(X;H^{-2s}(\mathbb{R}^n;\mathbb{C}))$, we have  $B\in\mathcal{L}(U;H^{-2s}(\mathbb{R}^n;\mathbb{C}))$.
    This, together with  the fact that $X_{-\gamma}=H^{-2\gamma}(\mathbb{R}^n;\mathbb{C})$, yields that Assumption $(A_2)$ with $\gamma=s$  holds, where $X_{-\gamma}$ is defined by \eqref{yu-3-26-5}.
    Let
$\mathcal{X}$ be $L^2(\mathbb{R}^n;\mathbb{C})$  if $s\in[0,\frac{1}{2})$ and $\mathcal{X}$ be $H^{-1}(\mathbb{R}^n;\mathbb{C})$ if $s\in[\frac{1}{2},1)$.
\begin{proposition}\label{yu-proposition-4-27-1}
    The following statements are equivalent:
\begin{enumerate}
  \item [$(i)$] The system \eqref{yu-4-24-2} is rapidly stabilizable in $\mathcal{X}$.
  \item [$(ii)$] The system \eqref{yu-4-24-2} is rapidly stabilizable in $L^2(\mathbb{R}^n;\mathbb{C})$ with bounded feedback law.
  \item [$(iii)$] For each $\alpha>0$, there  exists $C(\alpha)>0$ so that
\begin{equation}\label{yu-4-25-11}
    \|\varphi\|^2_{L^2(\mathbb{R}^n;\mathbb{C})}\leq C(\alpha)(\|(\lambda I-\triangle)\varphi\|_{L^2(\mathbb{R}^n;\mathbb{C})}^2+
   \|a_\omega(-\triangle)^s\varphi\|_{L^2(\mathbb{R}^n;\mathbb{C})}^2)
\end{equation}
    for any $\lambda\in\mathbb{C}_{-\alpha}^+$ and $\varphi\in H^2(\mathbb{R}^n;\mathbb{C})$.
\end{enumerate}
\end{proposition}
\begin{proof}
   Note that $\triangle$ can be extended to $H^{-1}(\mathbb{R}^n;\mathbb{C})$ (with domain $H^1(\mathbb{R}^n;\mathbb{C})$) or $H^{-2}(\mathbb{R}^n;\mathbb{C})$ (with domain $L^2(\mathbb{R}^n;\mathbb{C})$) by a direct way and $X_{-\frac{1}{2}}=H^{-1}(\mathbb{R}^n;\mathbb{C})$. Moreover, it follows from \eqref{yu-4-24-4} and Plancherel's identity
  that for any $u\in L^2(\mathbb{R}^n;\mathbb{C})$ and $\varphi\in H^2(\mathbb{R}^n;\mathbb{C})$,
\begin{eqnarray*}\label{yu-4-25-12}
    \langle Bu,\varphi\rangle_{H^{-2s}(\mathbb{R}^n;\mathbb{C}), H^{2s}(\mathbb{R}^n;\mathbb{C})}
    &=&(2\pi)^{-n}\langle (1+|\xi|^2)^{-s}|\xi|^{2s}\mathcal{F}[a_\omega u],(1+|\xi|^2)^s\mathcal{F}[\varphi]\rangle_{L^2(\mathbb{R}^n;\mathbb{C})}\nonumber\\
    &=&(2\pi)^{-n}\langle\mathcal{F}[a_\omega u], |\xi|^{2s}\mathcal{F}[\varphi]\rangle_{L^2(\mathbb{R}^n;\mathbb{C})}
    =\langle u,a_\omega(-\triangle)^s\varphi\rangle_{L^2(\mathbb{R}^n;\mathbb{C})}.
\end{eqnarray*}
    (Here, we recall that for each $\beta\in\mathbb{R}$ and
    $\varphi\in H^{2\beta}(\mathbb{R}^n;\mathbb{C})$, $\mathcal{F}[(I-\triangle)^\beta\varphi](\xi)
    =(1+|\xi|^2)^\beta\mathcal{F}[\varphi](\xi)$. It is clear that $(I-\triangle)^{\beta}\in\mathcal{L}(H^{2\beta}(\mathbb{R}^n;\mathbb{C});
    L^2(\mathbb{R}^n;\mathbb{C}))$.)
    The above estimates yield that $B^*=a_\omega(-\triangle)^s$. Therefore, by Corollary \ref{yu-corollary-4-22-1}, we finish the proof.
\end{proof}
   \begin{definition}\label{yu-definition-5-9-111}
    Let $E$ be a measurable set in $\mathbb{R}^n$. If there exist $\varepsilon>0$ and $L>0$ so that
    $|E\cap Q_L(x)|\geq \varepsilon L^n$ for each $x\in\mathbb{R}^n$, then $E$ is called the thick set in $\mathbb{R}^n$, where $|\cdot|$ is the Lebesgue measure in $\mathbb{R}^n$ and $Q_L(x)$ is the cube in $\mathbb{R}^n$, centered at $x$ and of side length $L$.
\end{definition}
\begin{lemma}\label{yu-lemma-4-26-1}
    \cite[Theorem 1]{Kovrijkine} If $\omega$ is a thick set in $\mathbb{R}^n$, then for each $R>0$, there  exists $C(R,\omega)>0$ so that
    $\|f\|_{L^2(\mathbb{R}^n;\mathbb{C})}\leq C(R,\omega)\|\chi_\omega f\|_{L^2(\mathbb{R}^n;\mathbb{C})}$ holds for any $f\in L^2(\mathbb{R}^n;\mathbb{C})$ with $\mbox{supp}(\mathcal{F}[f])\subset [-R,R]^n$. Here and in what follows, $\chi_\omega$ is the characteristic function of $\omega$.
\end{lemma}
\begin{theorem}\label{yu-theorem-4-27-1}
    If $\omega$ is a thick set in $\mathbb{R}^n$, then the system \eqref{yu-4-24-2} is rapidly stabilizable with bounded feedback law in $L^2(\mathbb{R}^n;\mathbb{C})$.
\end{theorem}
\begin{proof}
    By Proposition \ref{yu-proposition-4-27-1}, we only need to show that for each $\alpha>0$, there exists $C(\alpha)>0$ so that \eqref{yu-4-25-11} holds. Let $\varphi\in H^2(\mathbb{R}^n;\mathbb{C})$ and $\alpha>0$ be fixed arbitrarily. On one hand, it follows from Plancherel's identity that for each $\lambda\in\mathbb{C}_{-\alpha}^+$,
\begin{equation*}\label{yu-4-27-1}
    \|(\lambda I-\triangle)\varphi\|^2_{L^2(\mathbb{R}^n;\mathbb{C})}=
    (2\pi)^{-n}\|(\lambda+|\xi|^2)\mathcal{F}[\varphi]\|^2_{L^2(\mathbb{R}^n;\mathbb{C})}
    \geq (2\pi)^{-n}\alpha^2\int_{\mathbb{R}^n}\chi_{|\xi|\geq \sqrt{2\alpha}}|\mathcal{F}[\varphi](\xi)|^2d\xi
\end{equation*}
 and
\begin{equation}\label{yu-4-27-3}
    \|\varphi\|^2_{L^2(\mathbb{R}^n;\mathbb{C})}\leq \alpha^{-2}\|(\lambda I-\triangle)\varphi\|^2_{L^2(\mathbb{R}^n;\mathbb{C})}
    +(2\pi)^{-n}\int_{\mathbb{R}^n}\chi_{|\xi|\leq \sqrt{2\alpha}}|\mathcal{F}[\varphi](\xi)|^2d\xi.
\end{equation}
On the other hand,  we have that
\begin{eqnarray}\label{yu-4-27-2}
     \|(\lambda I-\triangle)\varphi\|^2_{L^2(\mathbb{R}^n;\mathbb{C})}&\geq& (2\pi)^{-n}\int_{\mathbb{R}^n}(|\xi|^2-\alpha)^2(1+|\xi|^2)^{-2s}
     |\xi|^{2s}|\mathcal{F}[\varphi](\xi)|^2d\xi\nonumber\\
    &\geq& (2\pi)^{-n}C(s,\alpha)\int_{\mathbb{R}^n}\chi_{|\xi|\geq \sqrt{2\alpha}}|\xi|^{2s}|\mathcal{F}[\varphi](\xi)|^2d\xi
\end{eqnarray}
 for each $\lambda\in\mathbb{C}_{-\alpha}^+$. Here, $C(s,\alpha):=\left(\inf_{|\xi|\geq \sqrt{2\alpha}}\frac{|\xi|^2-\alpha}{(1+|\xi|^2)^s}\right)^2>0$.
    Then by Lemma \ref{yu-lemma-4-26-1}, \eqref{yu-4-24-1} and \eqref{yu-4-27-2}, we observe that
\begin{eqnarray*}\label{yu-4-27-4}
    &\;&\int_{\mathbb{R}^n}\chi_{|\xi|\leq \sqrt{2\alpha}}|\mathcal{F}[\varphi](\xi)|^2d\xi\nonumber\\
    &\leq& \int_{\mathbb{R}^n}\chi_{|\xi|\leq \sqrt{2\alpha}}|\xi|^{2s}|\mathcal{F}[\varphi](\xi)|^2d\xi
    \leq C(\omega,\alpha)\int_{\mathbb{R}^n}\chi_{\omega}|\mathcal{F}^{-1}[\chi_{|\xi|\leq \sqrt{2\alpha}}|\xi|^{2s}\mathcal{F}[\varphi]](x)|^2dx\nonumber\\
    &\leq& 2C(\omega,\alpha)\int_{\mathbb{R}^n}\chi_\omega |(-\triangle)^s\varphi(x)|^2dx
    +2C(\omega,\alpha)(2\pi)^{-n}\int_{\mathbb{R}^n}\chi_{|\xi|\geq \sqrt{2\alpha}}|\xi|^{2s}|\mathcal{F}[\varphi](\xi)|^2d\xi\nonumber\\
    &\leq&2C(\omega,\alpha)c^{-2}\|a_\omega(-\triangle)^s\varphi\|^2_{L^2(\mathbb{R}^n;\mathbb{C})}
    +2C(\omega,\alpha) C(s,\alpha)^{-1}\|(\lambda I-\triangle)\varphi\|^2_{L^2(\mathbb{R}^n;\mathbb{C})}
\end{eqnarray*}
  for each $\lambda\in\mathbb{C}_{-\alpha}^+$.  Since $\varphi$ and $\alpha$ are arbitrary, it follows from the latter inequality and \eqref{yu-4-27-3}
   that for each $\alpha>0$, there exists $C(\alpha)>0$ so that \eqref{yu-4-25-11} holds. The proof is completed.
\end{proof}

\section{Appendices}\label{yu-sec-6}

\subsection{Equivalence of $L$ and $L^*$ on (AEDC)}\label{yu-sec-7-14-1}
\begin{lemma}\label{yu-lemms-25-7-14-1}
     The operator $L$ with domain $D(L)(\subset X)$ satisfies (AEDC) in $X$  if and only if $L^*$ with domain $D(L^*)$ satisfies (AEDC) in $X$.
\end{lemma}
\begin{proof}
    Since $(L^*)^*=L$, it suffices to show that if  $L$ with domain $D(L)$ satisfies (AEDC) in $X$, then $L^*$ with domain $D(L^*)$ satisfies (AEDC) in $X$.
\par
    Suppose that  $L$ with domain $D(L)$ satisfies (AEDC) in $X$. According to Definition \ref{yu-def-7-2-1},
    $L$ with domain $D(L)$ generates a $C_0$-semigroup $S_L(\cdot)$ on $X$, and for each $\alpha>0$, there are two closed subspaces $Q_1:=Q_1(\alpha)$ and $Q_2:=Q_2(\alpha)$ of $X$ so that $(a)$ $X=Q_1\oplus Q_2$;  $(b)$ $Q_1$ and $Q_2$ are  invariant subspaces of $S_L(\cdot)$;
    $(c)$ $L|_{Q_1}$ is bounded and satisfies that  $\sigma(L|_{Q_1})\subset\mathbb{C}_{-\alpha}^+$; $(d)$  $S_L(\cdot)|_{Q_2}$ is  exponentially stable.
     By \cite[Corollary 10.6, Section 1.10, Chapter 1]{Pazy}, we have that $L^*$ with domain $D(L^*)$ generates a $C_0$-semigroup
    $S^*_L(\cdot)$ on $X$ and
$S^*_L(\cdot)=(S_L(\cdot))^*$.

    Let $\alpha>0$ be fixed arbitrarily. We claim that there exist two closed subspaces $H_1$ and $H_2$ of $X$ so that $(a')$ $X=H_1\oplus H_2$;  $(b')$ $H_1$ and $H_2$ are  invariant subspaces of $S^*_L(\cdot)$;
    $(c')$ $L^*|_{H_1}$ is bounded and satisfies that  $\sigma(L^*|_{H_1})\subset\mathbb{C}_{-\alpha}^+$; $(d')$  $S^*_L(\cdot)|_{H_2}$ is  exponentially stable. If it can be done, then by Definition \ref{yu-def-7-2-1} and the arbitrariness of $\alpha$, we have that $L^*$ satisfies (AEDC) in $X$.
    To this end,  we define
$P_Lx=x_1$ if $x=x_1+x_2$ with $(x_1,x_2)\in Q_1\times Q_2$.
    By $(a)$, the closedness of $Q_1$ and $Q_2$, we can easily check that  $P_L$ is a projection operator (from $X$ to $Q_1$).
    Then $P^*_L$ is also a projection operator and
\begin{equation}\label{yu-25-7-14-3}
    \mbox{Ker}(P^*_L)=(I-P^*_L)X.
\end{equation}
    Let
\begin{equation}\label{yu-25-7-14-4}
    H_1:=P^*_LX\;\;\mbox{and}\;\;H_2:=(I-P_L^*)X.
\end{equation}
In the rest of the proof, we aim to show that $H_1$ and $H_2$ satisfy $(a')$-$(d')$.
Firstly, by \eqref{yu-25-7-14-3} and \eqref{yu-25-7-14-4}, we have $X=H_1\oplus H_2$, i.e., $(a')$ holds.
Secondly, for each
 $t\in\mathbb{R}^+$ and $h\in H_1$, it follows from $(b)$ and \eqref{yu-25-7-14-4} that
$\langle (I-P^*_L)S^*_L(t)h,g\rangle_X=\langle (I-P^*_L)h,S_L(t)g\rangle_X=0$ for any $g\in X$.
This implies that $(I-P^*_L)S^*_L(t)h=0$, i.e., $S^*_L(t)h\in H_1$, which indicates that $H_1$ is an invariant subspace of $S_L^*(\cdot)$.
Similarly, $H_2$ is also an invariant subspace of $S_L^*(\cdot)$. Hence, $(b')$ is true.
Thirdly, we observe the following facts:
\begin{enumerate}
  \item [(O1)] The boundness of $L|_{Q_1}$  is equivalent to that there  exists $C>0$ so that $|\langle L^*f,h\rangle_X|\leq C\|f\|_X\|h\|_X$ for any $f\in D(L^*)$ and $h\in Q_1$.
  \item [(O2)] $H_1=Q'_1$ and $H_2=Q_2'$ in the following sense:  for any linear bounded functionals $F(\cdot)$ and $G(\cdot)$ on $Q_1$ and $Q_2$, respectively, there are unique $f\in H_1$ and $g\in H_2$ so that $F(x)=\langle f,x\rangle_X$ for any $x\in Q_1$ and $G(x)=\langle g,x\rangle_X$ for any $x\in Q_2$.
\end{enumerate}
    (Indeed, (O1) follows from $(iii)$ in Remark \ref{yu-remark-12-22-1}, and (O2) follows from the same arguments as those to show \eqref{yu-12-12-6-01}.) By $(c)$ and (O1), we have that
\begin{equation}\label{yu-25-7-14-6}
    Q_1=P_LD(L)\subset D(L)\;\;\mbox{and}\;\;\|LP_L\|_{\mathcal{L}(X)}\leq C\|P_L\|_{\mathcal{L}(X)}.
\end{equation}
     It follows from the first conclusion in (\ref{yu-25-7-14-6}) and $(b)$ that $(I-P_L)D(L)\subset D(L)$ and
     $LP_L=P_LL$ on $D(L)$. These, along with the inequality in \eqref{yu-25-7-14-6},  imply that
\begin{eqnarray*}\label{yu-25-7-14-7}
    |\langle f, Lh\rangle_X|\leq|\langle f, LP_Lh\rangle_X|+
    |\langle f,(I-P_L)Lh\rangle_X
    =|\langle f, LP_Lh\rangle_X|\leq C\|P_L\|_{\mathcal{L}(X)}\|f\|_X\|h\|_X
\end{eqnarray*}
   for any $f\in H_1$ and $h \in D(L)$, which indicates that $L^*|_{H_1}$ is bounded and $H_1=P^*_LD(L^*)\subset D(L^*)$. These, together with the first conclusion in \eqref{yu-25-7-14-6}
    again and (O2), yield that  for any $f\in H_1$ and $g\in Q_1$,
\begin{equation*}\label{yu-25-7-14-8}
    \langle L^*f,g\rangle_{H_1,Q_1}=\langle L^*f, g\rangle_X=\langle f, Lg\rangle_X=\langle f,Lg\rangle_{H_1, Q_1}.
\end{equation*}
    From the latter it follows that $(L^*|_{H_1})^*=L|_{Q_1}$ (in the sense:
    $\langle L^*|_{H_1}f,g\rangle_{H_1,Q_1}=\langle f,L|_{Q_1}g\rangle_{H_1,Q_1}$ for any $f\in H_1$ and $g\in Q_1$). Thus, by $(c)$, we get $\sigma(L^*|_{H_1})\subset \mathbb{C}_{-\alpha}^+$. Finally, we prove $(d')$. Indeed, by $(b)$, we have that for any $t\in\mathbb{R}^+$, $f\in H_2$ and $g\in X$,
\begin{equation*}\label{yu-25-7-14-9}
    \langle S^*_L(t)f,g\rangle_X=\langle f, S_L(t)((I-P_L)g+P_Lg)\rangle_X=\langle f, S_L(t)(I-P_L)g\rangle_X.
\end{equation*}
    This, along with $(d)$, implies that  there exist $\varepsilon>0$ and $C(\varepsilon)>0$ so that for any $f\in H_2$ and $g\in X$,
\begin{equation*}\label{yu-25-7-14-10}
    |\langle S^*_L(t)f,g\rangle_X|\leq C(\varepsilon)\|I-P_L\|_{\mathcal{L}(X)} e^{-\varepsilon t}\|f\|_X\|g\|_X
    \;\;\mbox{for each}\;\;t\in\mathbb{R}^+,
\end{equation*}
    which indicates $(d')$.
    The proof is completed.
\end{proof}
\subsection{Proof of  the claim in $(iv)$ of Remark \ref{yu-remark-12-22-1}}\label{yu-sec-6.1}

     The following lemma  can be deduced from  \cite[Corollary 3.2 in Chapter V and Corollary 3.12 in Chapter IV]{Engel-Nagel} and \cite[Theorem 6.17, Section 6, Chapter 3]{Kato}.
\begin{lemma}\label{yu-lemma-6-23-1}
   Suppose that $L$ with domain $D(L)$ generates a compact semigroup $S_L(\cdot)$ in $X$, i.e., $(a)$ in $(iv)$ of Remark \ref{yu-remark-12-22-1} holds. Then for each $\alpha>0$, there is  a bounded closed region $M$ in $\mathbb{C}_{-\alpha}^+$, whose boundary (denoted by $\partial M$) is a rectifiable and simple curve with the counter-clockwise direction, so that $\sigma(L)\cap \partial M=\emptyset$ and the  Kato projection operator $\mathcal{P}(\alpha)\in\mathcal{L}(X)$ defined by
$\mathcal{P}(\alpha):=\frac{1}{2\pi \mathrm{i}}\int_{\partial M}(\lambda I-L)^{-1}d\lambda\footnote{Here,  $\mathrm{i}$ represents the unit imaginary number.}$
     satisfies
$(q_1)$ $\mathcal{P}(\alpha)D(L)\subset D(L)$, $L\mathcal{P}(\alpha)=\mathcal{P}(\alpha)L$;
  $(q_2)$ $\sigma((\mathcal{P}(\alpha)L)|_{\mathcal{P}(\alpha)X})\subset \mathbb{C}_{-\alpha}^+$;
  $(q_3)$ $\mathcal{P}(\alpha)X$ is a finite dimensional space;
  $(q_4)$ for each $\beta\in(0,\alpha)$, there  exists $C(\alpha,\beta)>0$ so that $\|(I-\mathcal{P}(\alpha))S_L(t)\|_{\mathcal{L}(X)}\leq C(\alpha,\beta)e^{-\beta t}$ for any $t\in\mathbb{R}^+$.
\end{lemma}
We now prove the claim in $(iv)$ of Remark \ref{yu-remark-12-22-1}.
  \begin{lemma}\label{yu-lemma-12-20-1-bb}
      Suppose that $L$ with domain $D(L)$ generates a $C_0$-semigroup $S_L(\cdot)$ in $X$.
       If it satisfies one of the conditions $(a)$ and $(b)$ in $(iv)$ of Remark \ref{yu-remark-12-22-1}, then  $L$ satisfies (AEDC) in $X$.
\end{lemma}
\begin{proof}
    The proof is divided into two cases.
\vskip 5pt
    \emph{Case 1. $S_L(\cdot)$ is a compact semigroup on $X$.} For each $\alpha>0$, we fix it. According to  Lemma \ref{yu-lemma-6-23-1},
    there exists  a bounded closed region $M$ in $\mathbb{C}_{-\alpha}^+$, whose boundary (denoted by $\partial M$) is a rectifiable and simple curve with the counter-clockwise direction so that $\sigma(L)\cap \partial M=\emptyset$ and the following $(v_1)-(v_4)$ hold:
 $(v_1)$ $\mathcal{P}(\alpha)D(L)\subset D(L)$, $L\mathcal{P}(\alpha)=\mathcal{P}(\alpha)L$;
  $(v_2)$ $\sigma((\mathcal{P}(\alpha)L)|_{\mathcal{P}(\alpha)X})\subset \mathbb{C}_{-\alpha}^+$;
  $(v_3)$ $\mathcal{P}(\alpha)X$ is a finite dimensional space;
  $(v_4)$ for each $\beta\in(0,\alpha)$, there  exists $C(\alpha,\beta)>0$ so that $\|(I-\mathcal{P}(\alpha))S_L(t)\|_{\mathcal{L}(X)}\leq C(\alpha,\beta)e^{-\beta t}$ for any $t\in\mathbb{R}^+$.
  Here,  $\mathcal{P}(\alpha)$ is the Kato projection operator which is given in Lemma \ref{yu-lemma-6-23-1}.

  Let $Q_1:=\mathcal{P}(\alpha)X$ and $Q_2:=(I-\mathcal{P}(\alpha))X$. We aim to show $(a)$-$(d)$ in $(ii)$ of Definition \ref{yu-def-7-2-1}.
  Since $(a)$ and $(b)$ are obvious, it suffices to check $(c)$ and $(d)$. Indeed, by $(v_3)$, we suppose that  $N:=\mbox{Dim}(Q_1)$. Let $\Phi$ be a linear isomorphic mapping from $Q_1$ to $\mathbb{C}^N$. It is clear that $L_N:=\Phi L\Phi^{-1}$ is bounded (which is in $\mathbb{C}^{N\times N}$) and $L|_{Q_1}=\Phi^{-1} L_N\Phi$. These imply $(c)$.
  From $(v_1)$ it follows that $\mathcal{P}(\alpha)S_L(\cdot)
  =S_L(\cdot)\mathcal{P}(\alpha)$. This, along with $(v_4)$, implies $(d)$.
  In summary, $L$ satisfies (AEDC) in $X$.
\vskip 5pt
    \emph{Case 2. $L$ with domain $D(L)$ is normal and $\sigma(L)\cap \mathbb{C}_{-\gamma}^+$ is bounded for each $\gamma\in\mathbb{R}^+$.} Let
    $\mathfrak{M}^{L}$ be the (unique) spectral measure (or spectral decomposition) of $L$ (the existence and uniqueness of $\mathfrak{M}^{L}$ can be found in \cite[Theorem 13.33, Chapter 13]{Rudin}, which is called the spectral theorem). For each $\alpha>0$, we fix it. It is well known that $\mathfrak{M}^{L}(\mathbb{C}_{-\alpha}^+)$ is self-adjoint and  orthogonal on $X$. Since $\sigma(L)\cap \mathbb{C}_{-\gamma}^+$ is bounded for each $\gamma>0$,
     it follows from the spectral theorem that $\mathfrak{M}^{L}(\mathbb{C}_{-\alpha}^+)X\subset D(L)$, $L\mathfrak{M}^{L}(\mathbb{C}_{-\alpha}^+)$ is bounded and $L\mathfrak{M}^{L}(\mathbb{C}_{-\alpha}^+)x=\mathfrak{M}^{L}(\mathbb{C}_{-\alpha}^+)Lx$ for any $x\in D(L)$. These imply that $\mathfrak{M}^{L}(\mathbb{C}_{-\alpha}^+)X$ is  an invariant subspace of $L$ (also $L\mathfrak{M}^{L}(\mathbb{C}_{-\alpha}^+)$) and $L\mathfrak{M}^{L}({\mathbb{C}_{-\alpha}^+})\in\mathcal{L}(\mathfrak{M}^{L}
    ({\mathbb{C}_{-\alpha}^+})X)$.
    Hence, $\mathfrak{M}^{L}(\mathbb{C}^+_{-\alpha})X$ is an invariant subspace of $S_L(\cdot)$ and
     $\mathfrak{M}^{L}(\mathbb{C}^+_{-\alpha})S_L(\cdot)(=S_L(\cdot)\mathfrak{M}^{L}
     ({\mathbb{C}_{-\alpha}^+}))$ is the $C_0$-group on $X$, which is generated by  $L\mathfrak{M}^{L}({\mathbb{C}_{-\alpha}^+})$.

    Let
$Q_1:=\mathfrak{M}^{L}({\mathbb{C}_{-\alpha}^+})X$, $Q_2:=(I-\mathfrak{M}^{L}({\mathbb{C}_{-\alpha}^+}))X$
    and
\begin{equation}\label{yu-12-10-11}
    V(\cdot):=(I-\mathfrak{M}^{L}({\mathbb{C}_{-\alpha}^+}))S_L(\cdot).
\end{equation}
We are devoted to proving $(a)$-$(d)$ in $(ii)$ of Definition \ref{yu-def-7-2-1} one by one. Firstly, since $\mathfrak{M}^{L}({\mathbb{C}_{-\alpha}^+})$ is a projection operator,
    it is obvious that $X=Q_1\oplus Q_2$. Secondly, since $Q_1$ is an invariant subspace of $S_L(\cdot)$, it follows from $(a)$ that $Q_2$ is also an invariant subspace of $S_L(\cdot)$. Thirdly, by the definition of $Q_1$, we have that
    $L|_{Q_1}=\mathfrak{M}^{L}({\mathbb{C}_{-\alpha}^+})L$ on $Q_1$ (here, we note that
    $Q_1\in D(L)$ and $\mathfrak{M}^{L}({\mathbb{C}_{-\alpha}^+})Lx=L\mathfrak{M}^{L}({\mathbb{C}_{-\alpha}^+})x$
    for any $x\in D(L)$), which, combined with the boundness of $\mathfrak{M}^{L}({\mathbb{C}_{-\alpha}^+})L$, indicates that $L|_{Q_1}$ is bounded.
    The conclusion that $\sigma(L|_{Q_1})\subset C_{-\alpha}^+$ follows directly from the spectral theorem. Finally,
    it follows from the invariance of $S_L(\cdot)$ on $Q_2$ that $V(\cdot)$ defined by \eqref{yu-12-10-11} is $C_0$-semigroup on $Q_2$.
    Suppose that the generator of $V(\cdot)$ (on $Q_2$) is $J$.
    We claim that
\begin{equation}\label{yu-12-12-13}
    J=(I-\mathfrak{M}^{L}({\mathbb{C}_{-\alpha}^+}))L\;\;\mbox{and}\;\;
    D(J)=(I-\mathfrak{M}^{L}({\mathbb{C}_{-\alpha}^+}))D(L).
\end{equation}
    Indeed, on one hand, if $f\in (I-\mathfrak{M}^{L}({\mathbb{C}_{-\alpha}^+}))D(L)$, i.e., $f=(I-\mathfrak{M}^{L}({\mathbb{C}_{-\alpha}^+}))g$ for some $g\in D(L)$,
    then
\begin{equation}\label{yu-12-12-14}
    \lim_{t\to 0^+}t^{-1}(V(t)f-f)=(I-\mathfrak{M}^{L}({\mathbb{C}_{-\alpha}^+}))\lim_{t\to 0^+}
    t^{-1}(S_L(t)g-g)=(I-\mathfrak{M}^{L}({\mathbb{C}_{-\alpha}^+}))Lg\in X.
\end{equation}
    It follows that $f\in D(J)$ and then $(I-\mathfrak{M}^{L}({\mathbb{C}_{-\alpha}^+})D(L)\subset D(J)$.
    On the other hand, if $f\in D(J)$, then  $f\in Q_2$ (i.e., $f=(I-\mathfrak{M}^{L}({\mathbb{C}_{-\alpha}^+}))f$) and
\begin{eqnarray*}\label{yu-12-12-15}
    &\;&\lim_{t\to 0^+}t^{-1}(S_L(t)(I-\mathfrak{M}^{L}({\mathbb{C}_{-\alpha}^+}))f
    -(I-\mathfrak{M}^{L}({\mathbb{C}_{-\alpha}^+}))f)\nonumber\\
    &=&\lim_{t\to0^+}t^{-1}((I-\mathfrak{M}^{L}({\mathbb{C}_{-\alpha}^+}))S_L(t)f-f)
    =\lim_{t\to0^+}t^{-1}(V(t)f-f)=Jf\in X.
\end{eqnarray*}
    It follows that $f=(I-\mathfrak{M}^{L}({\mathbb{C}_{-\alpha}^+}))f\in D(L)$ and then $D(J)\subset (I-\mathfrak{M}^{L}({\mathbb{C}_{-\alpha}^+}))D(L)$. Therefore, $D(J)=(I-\mathfrak{M}^{L}({\mathbb{C}_{-\alpha}^+}))D(L)$. This, together with
     \eqref{yu-12-12-14}, $L\mathfrak{M}^{L}({\mathbb{C}_{-\alpha}^+})=\mathfrak{M}^{L}({\mathbb{C}_{-\alpha}^+})L$ on $D(L)$ and the fact $(I-\mathfrak{M}^{L}({\mathbb{C}_{-\alpha}^+}))^2=I-\mathfrak{M}^{L}({\mathbb{C}_{-\alpha}^+})$, yields that $J=(I-\mathfrak{M}^{L}({\mathbb{C}_{-\alpha}^+}))L$ and \eqref{yu-12-12-13} holds.
     By \eqref{yu-12-12-13}, the fact that $\mathfrak{M}^{L}({\mathbb{C}_{-\alpha}^+})$  is self-adjoint
     and the spectral theorem, we have that $J$ is a normal operator and $\sigma(J)\subset \overline{\mathbb{C}_{-\alpha}^-}$. These, along
     with \cite[Corollary 3.14, Section 3, Chapter IV and Lemma 1.9, Section 1, Chapter V]{Engel-Nagel}, imply that for any $\beta\in(0,\alpha)$, there exists $C(\alpha,\beta)>0$ so that
$\|V(t)\|_{\mathcal{L}(X)}\leq C(\alpha,\beta)e^{-\beta t}$ for any $t\in\mathbb{R}^+$,
    which, combined with \eqref{yu-12-10-11} and the definition of $Q_2$, indicates $(d)$.
    The proof is completed.
\end{proof}

\end{document}